\documentclass[11pt,reqno]{amsart}

\usepackage[usenames,dvipsnames]{xcolor}
\usepackage{fancyhdr}
\usepackage{amsmath,amsfonts,amsbsy,amsgen,amscd,amssymb,amsthm}
\usepackage{url}

\usepackage{mathtools}

\mathtoolsset{showonlyrefs}

\usepackage[font=small,margin=0.25in,labelfont={bf},labelsep={colon}]{caption}

\usepackage{subcaption}

\usepackage{tikz}
\usepackage{microtype}
\usepackage{enumitem}

\definecolor{dark-gray}{gray}{0.3}
\definecolor{dkgray}{rgb}{.4,.4,.4}
\definecolor{dkblue}{rgb}{0,0,.5}
\definecolor{medblue}{rgb}{0,0,.75}
\definecolor{rust}{rgb}{0.5,0.1,0.1}

\usepackage[colorlinks=true]{hyperref}

\hypersetup{urlcolor=rust}
\hypersetup{citecolor=dkblue}
\hypersetup{linkcolor=dkblue}

\usepackage{setspace}

\usepackage{graphicx}
\usepackage{booktabs,longtable,tabu} 
\usepackage{multirow} 

\usepackage{float}

\usepackage[full]{textcomp}

\usepackage[scaled=.98,sups,osf]{XCharter}
\usepackage[scaled=1.04,varqu,varl]{inconsolata}
\usepackage[type1]{cabin}
\usepackage[charter,vvarbb,scaled=1.07]{newtxmath}
\usepackage[cal=boondoxo]{mathalfa}
\usepackage[T1]{fontenc}

\usepackage{bm} 

\graphicspath{{figures/}}

\newtheorem{theorem}{Theorem}[section]
\newtheorem{lemma}[theorem]{Lemma}

\newtheorem{proposition}[theorem]{Proposition}

\newtheorem{corollary}[theorem]{Corollary}

\theoremstyle{definition}

\newtheorem{remark}[theorem]{Remark}
\newtheorem{assumption}[theorem]{Assumption}

\numberwithin{equation}{section} 
\numberwithin{figure}{section}
\numberwithin{table}{section}

\floatstyle{plaintop}
\newfloat{recipe}{thp}{lor}
\floatname{recipe}{Recipe}
\numberwithin{recipe}{section}

\renewcommand{\phi}{\varphi}

\newcommand{\diff}{\mathrm{d}}
\newcommand{\idiff}{\, \diff}

\newcommand{\triplenorm}[1]{{\left\vert\kern-0.25ex\left\vert\kern-0.25ex\left\vert #1
    \right\vert\kern-0.25ex\right\vert\kern-0.25ex\right\vert}}

\def\XXint#1#2#3{{%
\setbox0=\hbox{$#1{#2#3}{\int}$}
\vcenter{\hbox{$#2#3$}}\kern-.5\wd0}}

\renewcommand{\leq}{\leqslant}
\renewcommand{\geq}{\geqslant}
\renewcommand{\subset}{\subseteq}

\usepackage{mathrsfs}
\allowdisplaybreaks
\hypersetup{bookmarksdepth=2}
\title[Global regularity and singularity formation for periodic gCLM]{The Periodic Generalized Constantin--Lax--Majda Equation II:\\ Global Regularity and Singularity Formation}
\author{Thomas Y. Hou}
\address{Applied and Computational Mathematics,
California Institute of Technology, Pasadena, CA, USA}
\email{hou@cms.caltech.edu}

\author{Xiang Qin}
\address{Applied and Computational Mathematics,
California Institute of Technology, Pasadena, CA, USA}
\email{xqin2@caltech.edu}

\author{Peicong Song}
\address{Applied and Computational Mathematics,
California Institute of Technology, Pasadena, CA, USA}
\email{psong2@caltech.edu}
\date{}
\begin{document}
\begin{abstract}
We study global regularity and finite-time singularity formation for the
periodic generalized Constantin--Lax--Majda equation. A finite
inverse-vorticity moment gives blowup for every \(0<a<1\), including the
smooth datum \(\omega_0=-\sin x\). Finite-dimensional selection
also realizes every focusing whole-line or finite-period companion profile
by smooth periodic blowup. For one critical parameter
in \((0.689,0.690)\), a computer-assisted coercivity verification yields
blowup with inverse spatial scale asymptotic to a positive multiple of
\(\sqrt{\log(1/(T-t))}\), where \(T\) is the blowup time.
These three profile constructions have type-I vorticity growth.
We also prove global well-posedness for \(a\geq1/2\) and odd
one-sided-sign data satisfying an origin bound and a condition at interior zeros.
\end{abstract}
\maketitle
\tableofcontents

\section{Introduction}
\label{sb:paper:introduction}

We consider the generalized Constantin--Lax--Majda equation on the circle,
\begin{equation}\label{sb:paper:equation}
 \omega_t+a u\omega_x=u_x\omega,\qquad
 u_x=H_{\mathbb T}\omega,\qquad u(0,t)=0,
 \qquad \mathbb T=\mathbb R/(2\pi\mathbb Z),
\end{equation}
where \(a\) measures the strength of transport relative to vorticity
stretching and \(H_{\mathbb T}\) denotes the periodic Hilbert transform.
Our purpose is to identify endpoint conditions for global regularity or
finite-time blowup, and to connect stationary solutions of the rescaled
equation with singular solutions arising from smooth periodic data.
An inverse-vorticity moment gives a direct blowup criterion throughout
\(0<a<1\). The three profile constructions describe the singularity more
precisely: they have type-I vorticity growth, while their physical spatial
scales are a power of the remaining time, an inverse square root of its
logarithm, or a fixed positive length.

The companion paper~\cite{Profiles} constructs physical steady states and
rescaled profiles by variational and fixed-point methods, respectively.
The selected rescaled profiles have three
possible geometries: a whole-line profile with a positive spatial rate,
a whole-line profile with zero spatial rate, and a finite-period profile.

Passing from these profiles to smooth periodic dynamics requires an
additional argument. A whole-line profile must be localized to a proper
arc of the circle. A finite-period profile can have a noninteger vanishing
order at the identified endpoints of its period, which we call the
\emph{seam}; it must be smoothed there. At a critical whole-line profile,
the stationary spatial rate vanishes, so the variation of the periodic
kernel determines the concentration scale. We address all three problems.
For the focusing and finite-period profiles, the normalized generator has
a negative essential growth bound. A finite-dimensional choice of smooth
initial perturbations then controls its remaining spectral directions.
For one selected critical profile, a verified coercive estimate gives the
stronger control needed to compute the logarithmic concentration law.

\subsection{Smooth periodic singularity formation}

An inverse-vorticity moment first gives a criterion for finite-time
blowup throughout \(0<a<1\).

\begin{theorem}[Blowup from a finite inverse-vorticity moment]
\label{sb:paper:inverse-blowup}
Let \(0<a<1\), and let \(\omega_0\in C^\infty(\mathbb T)\) be real,
odd, and strictly negative on \((0,\pi)\). If
\[
 \int_0^\pi(-\omega_0(x))^{-a}\,\idiff x<\infty,
\]
then the maximal lifespan \(T_*\) is finite, and the solution satisfies
\(\limsup_{t\uparrow T_*}\|\omega(t)\|_\infty=\infty\).
In particular, this holds for \(\omega_0(x)=-\sin x\) at every
\(0<a<1\).
\end{theorem}

Theorem~\ref{sb:inverse:blowup} proves this criterion through a conserved
inverse-vorticity measure and its increasing first spatial moment.
We next construct solutions with prescribed asymptotic profiles and
concentration scales.

The profile estimates in \cite[Corollary~5.3 and Proposition~6.2]{Profiles}
give two explicit sufficient parameter bounds:
\begin{equation}\label{sb:paper:parameter-thresholds}
 \begin{aligned}
 a_{\mathrm{foc}}&=\frac{17+16\log2}{53-16\log2}
       =0.670259908298743\ldots,\\
 a_{\mathrm{per}}&=\frac{64}{176-9\pi^2}
       =0.7341675585\ldots .
 \end{aligned}
\end{equation}
The first constant follows from a sharp first-moment inequality.
Strictness of that inequality at a fixed point includes
\(a=a_{\mathrm{foc}}\) in the focusing range. The second constant is the
upper bound on whole-line profiles obtained from the fixed-point moment
identities.

In the statements below, \(\bar\Omega\) denotes an odd stationary vorticity profile
of the normalized equation, \(\bar c_\omega<0\) is its amplitude rate, and
\(\lambda(t)=C_l(t)^{-1}\) is the physical spatial scale. These objects
are defined in Section~\ref{sb:common:framework}.

A rescaled profile is called \emph{noncritical} if it is either
focusing, with \(\bar L=\infty\) and \(\bar c_l>0\), or finite-period,
with \(\bar L<\infty\) and \(\bar c_l=0\).

\begin{theorem}[Smooth periodic blowup from noncritical profiles]
\label{sb:paper:main-blowup}
Let \(0<a<1\), and let \(\bar\Omega\) be any normalized noncritical profile
constructed in \cite{Profiles}. There are smooth odd
periodic initial data whose solution of \eqref{sb:paper:equation} has a
finite maximal time \(T\), remains smooth for \(0\leq t<T\), and satisfies
\begin{equation}\label{sb:paper:main-limits}
 (T-t)\|\omega(t)\|_{L^\infty(\mathbb T)}
       \longrightarrow-\frac{\|\bar\Omega\|_\infty}{\bar c_\omega}.
\end{equation}
The spatial behavior is determined by the profile geometry:
\begin{enumerate}[label=\textup{(\roman*)}]
\item If \(\bar\Omega\) is focusing, then for some \(\lambda_*>0\),
\begin{equation}\label{sb:paper:main-focusing}
 \lambda(t)\sim\lambda_*(T-t)^{-\bar c_l/\bar c_\omega},\qquad
 (T-t)\omega(\lambda(t)x,t)\longrightarrow-\frac{\bar\Omega(x)}{\bar c_\omega}
       \quad\hbox{in }C^4_{\mathrm{loc}}(\mathbb R).
\end{equation}
\item If \(\bar\Omega\) has finite half-period \(\bar L\), then
\begin{equation}\label{sb:paper:main-finite}
 \lambda=\frac\pi{\bar L},\qquad
 (T-t)\omega(y,t)\longrightarrow
 -\frac1{\bar c_\omega} \bar\Omega\!\left(\frac{\bar L}\pi y\right)
       \quad\hbox{in }C^1(\mathbb T).
\end{equation}
\end{enumerate}
In particular, such blowup exists for every
\(a\in(0,a_{\mathrm{foc}}]\cup[a_{\mathrm{per}},1)\): the profile can be
chosen focusing in the first interval and finite-period in the second.
At \(a=1/2\), an explicit profile gives \(-\bar c_l/\bar c_\omega=1/3\).
\end{theorem}

\begin{proof}
These profiles satisfy Assumption~\ref{sb:focus:profile-input}
in the focusing case and Proposition~\ref{sb:finite:profile-input}
in the finite-period case. Theorems~\ref{sb:focus:theorem} and
\ref{sb:finite:theorem} give the corresponding smooth periodic solutions,
and Lemma~\ref{sb:common:physical} gives \eqref{sb:paper:main-limits}.
The parameter assertions follow from
Corollary~\ref{sb:focus:parameter-range} and
\eqref{sb:finite:parameter-range}; the explicit example is
Corollary~\ref{sb:focus:explicit}.
\end{proof}

Both constructions are computer-free. They produce at least one choice
of finitely many initial coefficients while allowing unstable eigenvalues
and neutral Jordan blocks. Thus the conclusion is realization by smooth
data through finite-dimensional selection; strict stability of every
profile is not an input.

Realizing a critical whole-line profile with its concentration law requires
a separate argument. For one selected profile, a verified coercive estimate
gives the following logarithmic spatial scale.

\begin{theorem}[A selected critical parameter]
\label{sb:paper:main-critical-theorem}
For one parameter \(a_{\mathrm c}\in(0.689,0.690)\), there are smooth odd
periodic initial data whose solution blows up at a finite time \(T\), and
a positive even Schwartz function \(f\), such that, with
\(\bar\Omega=-xf\), \(b=H\bar\Omega(0)>0\), and \(\bar c_\omega=-(1-a_{\mathrm c})b\),
\begin{equation}\label{sb:paper:main-critical}
 \begin{gathered}
 C_l(t)^2=-\frac{a_{\mathrm c}\mathsf b_2}{6\bar c_\omega}
       \log\frac1{T-t}
       +O\!\left(\log\log\frac1{T-t}\right),
 \qquad \mathsf b_2=\frac2\pi\int_0^\infty x^2f(x)\,\idiff x>0,\\
 (T-t)\omega(C_l(t)^{-1}x,t)\longrightarrow-\frac{\bar\Omega(x)}{\bar c_\omega}
       \quad\hbox{locally uniformly on }\mathbb R.
 \end{gathered}
\end{equation}
The solution remains smooth for \(0\leq t<T\) and satisfies
\eqref{sb:paper:main-limits}.
\end{theorem}

\begin{proof}
Theorem~\ref{sb:crit:transfer} applies to the critical profile and coercive
estimate proved in Appendix~\ref{sb:sc:introduction}.
\end{proof}

We use \(a_*\) for the critical parameter obtained by analytical
continuation in \cite[Theorem~6.4]{Profiles} and \(a_{\mathrm c}\) for the profile
constructed together with a quantitative coercivity estimate here.
The equality of these profiles and uniqueness of a critical parameter
remain open. The bounds in \eqref{sb:paper:parameter-thresholds} give
explicit sufficient parameter ranges for the two noncritical constructions.
The quantitative input
for \eqref{sb:paper:main-critical} is stated in
Assumption~\ref{sb:dyn:profile-input}.

\subsection{Global well-posedness and the origin condition}

The positive transport velocity for one-sided-sign data moves interior
particles away from the origin. Global regularity depends on controlling
the growth of the derivative of this flow at the origin. The natural
endpoint quantity is
\begin{equation}\label{sb:paper:endpoint-quantity}
 m_a(\omega_0)=\operatorname*{ess\,sup}_{0<x<\pi}
                    \frac{|\omega_0(x)|}{x^{1/a}}.
\end{equation}
The next theorem concerns the original periodic equation and is independent
of the profile construction.

\begin{theorem}[Global well-posedness]\label{sb:paper:main-global}
Let \(a\geq1/2\), and suppose \(\omega_0\in H^1(\mathbb T)\) is real,
odd, and nonpositive on \((0,\pi)\). Assume
\begin{equation}\label{sb:paper:zero-compatibility}
 H_{\mathbb T}\omega_0(z)\leq0
 \qquad\hbox{whenever }0<z<\pi\hbox{ and }\omega_0(z)=0.
\end{equation}
If either \(m_a(\omega_0)<\infty\) or
\begin{equation}\label{sb:paper:weighted-endpoint}
 \int_0^\pi|\omega_0(x)|\cot^{1+1/a}(x/2)\,\mathrm{d}x<\infty,
\end{equation}
then the unique \(H^1\) solution is global and preserves oddness and its
one-sided sign. Every additional Sobolev regularity of the initial data
persists on finite time intervals. In particular, smooth initial data
under these hypotheses give a global smooth solution.
\end{theorem}

The proof is given in Theorem~\ref{sb:gwp:theorem}, after the local theory
and continuation criterion in Section~\ref{sb:gwp:local-section}.
Condition~\eqref{sb:paper:zero-compatibility} is automatic for strictly
negative data on \((0,\pi)\). For smooth odd data,
\(m_a(\omega_0)<\infty\) holds for \(a\geq1\); for
\(1/2\leq a<1\), it holds if the origin slope vanishes.
The origin condition is dynamically meaningful. Differentiating
\eqref{sb:paper:equation} at the symmetry point gives
\begin{equation}\label{sb:paper:origin-slope}
 \frac{\idiff}{\idiff t}\omega_x(0,t)
      =(1-a)u_x(0,t)\omega_x(0,t),
\end{equation}
so the zero-slope class is invariant. For \(1/2\leq a<1\), smooth odd
data in this class satisfy the origin bound because they are \(O(x^3)\).
The blowup data in Theorems~\ref{sb:paper:main-blowup} and
\ref{sb:paper:main-critical-theorem} have a nonzero origin slope, and
hence an infinite value of \eqref{sb:paper:endpoint-quantity} when \(a<1\).
At \(a=1\), Remark~\ref{sb:gwp:Chen-comparison} compares the interior-zero
condition with the weighted derivative hypothesis in Chen's global
regularity theorem~\cite[Theorem 2]{Chen2025Regularity}.

\begin{remark}[Comparison with the steady-state vanishing order]
\label{sb:inverse:vanishing-order}
The physical steady states in \cite[Proposition~2.6]{Profiles} have
vorticity vanishing order \(1/a\) at both endpoints.
To compare the two dynamical criteria, fix \(1/2\leq a<1\) and suppose
smooth odd data, strictly negative on \((0,\pi)\), have a common
endpoint vanishing order \(q>0\), in the sense that
\[
 -\omega_0(x)\asymp\operatorname{dist}(x,\{0,\pi\})^q
 \quad\hbox{near }\{0,\pi\}.
\]
The inverse-vorticity moment in Theorem~\ref{sb:paper:inverse-blowup}
is finite exactly when \(q<1/a\), whereas
\(m_a(\omega_0)<\infty\) exactly when \(q\geq1/a\).
Theorems~\ref{sb:paper:inverse-blowup} and \ref{sb:paper:main-global}
therefore give blowup in the first case and global regularity in the second.
Thus, within this endpoint class, the physical steady-state order
separates these two behaviors.
\end{remark}

\subsection{Proof strategy}

The inverse-vorticity density \((-\omega)^{-a}\) is transported as a
measure by the characteristic flow. Its total mass is conserved and its
first spatial moment is bounded. Positivity of the periodic velocity
kernel on a fixed interior interval, together with H\"older's inequality,
forces that moment to increase at a uniform positive rate when \(a<1\).
This proves Theorem~\ref{sb:paper:inverse-blowup}.

For a focusing profile, the logarithmic coordinate \(z=\log x\)
turns transport into a vector field whose speed has positive limits at
both ends. Two normalization identities cancel the first and third origin
terms of the nonlocal linearization. A weighted energy then gives negative
essential growth, with only finitely many nondecaying spectral directions.
Localization introduces a forcing term that decays as the rescaled period
grows. We control this forcing at a lower Sobolev level and close a bounded
higher energy. A topological exit argument selects the initial values of
the remaining finite-dimensional components.

For a finite-period profile, the coordinate adapted to the stationary
velocity gives the constant transport operator
\(-\partial_s-\delta\). The profile vanishes at the seam with an exponent
larger than one. Its ordinary high Sobolev norms can be infinite, so the
energy uses derivatives along this transport field. The corresponding weighted spaces allow smooth seam
cutoffs to approach the profile at every fixed conormal order. The same
finite-dimensional selection then gives a smooth periodic solution whose
normalized profile converges in \(C^1\) on the entire circle.

At the selected critical profile, strict coercivity controls the normalized
shape. The leading periodic Hilbert correction is the spatial constant
\(-q_0L^{-2}\), where \(q_0=\pi^2\mathsf b_2/12\).
Adding a constant \(d\) to the Hilbert transform changes the velocity
by \(d x\) and the rates \((c_l,c_\omega)\) by \((-ad,-d)\).
Consequently both coefficients \(c_lx+aU\) and \(c_\omega+H_L\Omega\)
in the normalized shape equation are unchanged. The correction remains
in the scale equation, giving
\[
 c_l=kL^{-2}+O(L^{-4}),\qquad
 \|\Omega_{\mathrm{ext}}-\bar\Omega\|_{Y_m}=O(L^{-4}),\qquad
 k=aq_0>0,
\]
as proved in Lemmas~\ref{sb:crit:periodic-source}--
\ref{sb:crit:periodic-estimates} and Theorem~\ref{sb:crit:transfer}.
Here \(\Omega_{\mathrm{ext}}\) is the zero extension of one cell and
\(Y_m=Y_\sigma^m\), for a fixed \(\sigma>6\), is the weighted space
defined in \eqref{sb:crit:polynomial-space}.
Since \((L^2)'=2k+O(L^{-2})\), this separation yields the logarithmic
concentration in \eqref{sb:paper:main-critical}. Proposition
\ref{sb:crit:second-law} also determines the next coefficient: with
\(s=\log(1/(T-t))\),
\[
 C_l(t)^2=-\frac{2k}{\pi^2\bar c_\omega}s
   +\frac1{\pi^2}\left(-\frac{k}{a\bar c_\omega}
                  +\frac{\varkappa}{k}\right)\log s+O(1),
\]
where \(\varkappa\) is defined by the linear corrector in
\eqref{sb:crit:corrector}. The computer-assisted part supplies the
quantitative profile and coercivity bounds used in this analysis.

For global regularity, the main estimate is a corrected moment at an
exponent just below the endpoint \(1+1/a\). A positive Fourier multiplier
controls its quadratic form, and a small multiple of the total mass
absorbs the contribution from the opposite endpoint. The resulting bound
controls the accumulated velocity gradient at the origin. A nonnegative
Lagrangian force and the condition at interior zeros then control vorticity
amplification along the flow, closing the continuation argument.

\subsection{Related work}

The gCLM equation isolates the competition between transport and vortex
stretching in the vorticity formulation of the three-dimensional Euler
equations.  The continuation criterion of Beale, Kato, and Majda
\cite{BKM1984} identifies the time integral of the maximum vorticity as a
quantity controlling classical continuation, while the geometric criteria
of Constantin, Fefferman, and Majda \cite{CFM1996} explain how coherence of
the vorticity direction can weaken stretching.  The role of convection and
dynamic depletion was further investigated by Okamoto and Ohkitani
\cite{OO2005}, Hou and Li \cite{HouLi2006,HouLi2008}, and Hou and Lei
\cite{HouLei2009}.  Although the gCLM model has no incompressibility
constraint or vorticity direction, it retains a nonlocal stretching rate
and permits the relative strength of transport to vary through one
parameter.  In two dimensions, the double-exponential small-scale creation
of Kiselev and \v Sver\'ak \cite{KS2014} provides a complementary example
in which transport produces rapid growth while the vorticity remains
bounded.

Self-similar constructions for Euler provide a second source of motivation.
Elgindi \cite{Elgindi2021} constructed finite-time singularities with
\(C^{1,\alpha}\) velocity, and Elgindi, Ghoul, and Masmoudi
\cite{EGM2021Euler} established stability of the corresponding mechanism.
For smooth data in a domain with boundary, the numerical scenario of Luo
and Hou \cite{LH2014PNAS,LH2014MMS} led to the analytic and computational
work of Chen and Hou \cite{CH2022Analysis,CH2025Numerics};
\cite{CH2025PNAS} presents the resulting boundary singularity construction.
The one-dimensional Hou--Luo model has also been treated through integral
estimates \cite{Choi2017} and, more recently, through a purely analytic
fixed-point construction of smooth profiles \cite{HQWW2025HL}.
Shao, Wei, Zhang, and Zhang \cite{SWZZ2026} construct self-similar profiles
for higher-dimensional axisymmetric Euler without swirl and prove a
finite-codimensional stability theorem.  These works distinguish two
steps that are also essential here: constructing a stationary solution of
the rescaled equation and realizing it as the asymptotic profile of an
evolving solution in the desired spatial domain and regularity class.

Constantin, Lax, and Majda \cite{CLM1985} introduced the model without
transport and obtained finite-time singularity by an explicit solution
formula.  De Gregorio \cite{DG1990,DG1996} added transport with the same
coefficient as stretching.  Okamoto, Sakajo, and Wunsch \cite{OSW2008}
introduced the parameterized family considered here, and Wunsch
\cite{Wunsch2011} developed further aspects of its analysis.
Periodic equilibria and traveling waves were studied by Okamoto, Sakajo,
and Wunsch \cite{OSW2014}.  The related nonlocal scalar transport equations
of C\'ordoba, C\'ordoba, and Fontelos \cite{CCF2005}, and Castro and
C\'ordoba \cite{CastroCordoba2010}, show how transport can also reinforce
singularity formation; Castro's thesis \cite{Castro2010Thesis} develops
these nonlocal fluid models in a broader setting.

The construction of smooth self-similar gCLM profiles began with the
perturbative analysis of Elgindi and Jeong \cite{EJ2020} near the explicit
CLM solution.  Elgindi, Ghoul, and Masmoudi \cite{EGM2021Transport} proved
stability in a corresponding small-parameter regime.  Chen, Hou, and
Huang \cite{CHH2021} developed a dynamic-rescaling and weighted-energy
method that also yields smooth blowup for the De Gregorio equation on
the real line.  The explicit profile at \(a=1/2\), found independently in
\cite{Chen2020,LSS2021}, provides another parameter around which stable
smooth blowup can be analyzed.  Huang, Tong, and Wei \cite{HTW2023}
subsequently constructed infinitely many compactly supported De Gregorio
profiles through a compact self-adjoint operator.  Exactly self-similar
solutions of lower H\"older regularity, for a different range of parameters
and regularity exponents, were constructed in \cite{EJ2020,Zheng2023}.

Huang, Qin, Wang, and Wei \cite{HQWW2024}, abbreviated HQWW, gave a unified
fixed-point construction for every \(a\leq1\).  Their normalized variable
\(f=-\bar\Omega/x\) is nonnegative, nonincreasing, and convex as a function of
\(x^2\).  These shape properties yield compactness and preserve a convex
set under the nonlinear reconstruction map.  The resulting profiles are
smooth in the interior of their supports, and their spatial scaling rate
distinguishes compact support, Gaussian logarithmic decay, and algebraic
decay. The focusing profiles used here belong to this family; the sharper
first-moment estimate in \cite{Profiles} gives the focusing bound in
\eqref{sb:paper:parameter-thresholds}.
For \(a>1\), Hou, Qin, and Song \cite{HyperConvective} construct a stable
local branch of compactly supported real-line profiles near \(a=3/2\).
Within this neighborhood, they obtain spatially expanding smooth-data
blowup for \(a<3/2\) and global decay for \(a>3/2\). They also prove
global smoothness for odd compactly supported smooth data under sign and
interior-zero hypotheses when \(a\geq2\), or when \(a\geq1\) and the
initial slope at the origin vanishes.

On the circle, the domain geometry introduces phenomena that are absent
from the whole-line profile equation.  Jia, Stewart, and \v Sver\'ak
\cite{JSS2019} analyzed the stability of sinusoidal De Gregorio equilibria.
Lei, Liu, and Ren \cite{LLR2020} established global regularity for
sign-definite vorticity and developed weighted estimates near an
equilibrium. Chen \cite{Chen2021} proved asymptotically self-similar
periodic blowup for \(a<1\) sufficiently close to \(1\), and global decay
for \(a>1\) sufficiently close to \(1\).
Chen \cite{Chen2025Regularity} obtained further global regularity and
blowup criteria in an odd, one-sided-sign class. The inverse-moment
criterion in Theorem~\ref{sb:paper:inverse-blowup} gives qualitative
periodic blowup throughout \(0<a<1\).
The numerical study \cite{LSS2021} identified different focusing
and finite-period blowup regimes and a transition near
\(a=0.6890665\).  This numerical value guides the critical construction;
the rigorous parameter conclusions follow from the analytic estimates
and, for the selected critical dynamics, the verified bounds stated in
this paper.

The recent periodic analysis of Guo and Jiu \cite{GuoJiu2025} treats
stability and instability near the first excited sinusoidal state.
Silantyev, Lushnikov, Siegel, and Ambrose \cite{SLAS2025} derive exact
periodic pole solutions for dissipative gCLM equations.
At \(a=0\), Huang, Qin, and Wang \cite{HQW2025Multiscale} constructed
asymptotically self-similar blowups with two collapsing spatial scales,
using the explicit CLM solution formula and its complex pole dynamics.
Huang, Tong, and Wang \cite{HTW2026Singular} study a different normalization
with derivative degeneracy at the origin.  Their work combines explicit
singular profiles and fixed-point constructions of traveling waves with
numerical investigations of additional one-scale and two-scale regimes.
Our normalization instead retains a nonzero vorticity slope at the origin.
Recent spectral analyses of Xu~\cite{Xu2026CLMSpectrum,Xu2026GCMLSpectrum}
study the dependence of the CLM and gCLM linearized spectra on the
operator realization. In the present constructions, the essential
semigroup bounds are proved in the weighted spaces used to close
the nonlinear evolution.

Dissipative variants have a related, but distinct, regularity theory.
Schochet \cite{Schochet1986} found explicit viscous-model solutions.
Transport equations with nonlocal drift and dissipation were studied by
Dong \cite{Dong2008}, Li and Rodrigo \cite{LiRodrigo2008}, and Silvestre
and Vicol \cite{SV2016}; Kiselev \cite{Kiselev2010} surveys the surrounding
active-scalar theory.  Chen \cite{Chen2020} and Ambrose, Lushnikov, Siegel,
and Silantyev \cite{ALSS2024} address singularity formation and global
existence for dissipative gCLM equations, with the latter emphasizing the
distinction between the real line and periodic domains.  These works
underscore the need to retain the periodic Biot--Savart law throughout
the reconstruction and dynamical arguments.

\subsection{Conventions and organization}

We write \((\omega,u)\) for the physical variables, \((\Omega,U)\) for the
rescaled variables, and \((\bar\Omega,\bar U)\) for stationary profiles.
The rescaled half-period is \(L\), and its stationary
value is denoted by \(\bar L\). We write \(H=H_\infty=H_{\mathbb R}\)
for the whole-line transform. The stationary rates are
\(\bar c_l,\bar c_\omega\), and \(f=-\bar\Omega/x\) is the scalar profile;
its focusing decay exponent is \(\rho\), and its finite-endpoint
vanishing exponent is \(p_L\).
Auxiliary weights, coordinates, and regularity exponents are defined
locally in the section where they are used.

The remainder of the paper is organized as follows.
Section~\ref{sb:common:framework} fixes the Hilbert-transform convention,
normalization, and physical reconstruction formulas.
Section~\ref{sb:gwp:local-section} establishes the local theory, and
Section~\ref{sb:gwp:section} proves global well-posedness and the
inverse-moment blowup criterion.
Sections~\ref{sb:focus:section}, \ref{sb:paper:critical}, and
\ref{sb:finite:section} treat focusing, critical, and fixed-period blowup,
respectively. Each construction states the profile hypotheses it uses.
Appendix~\ref{sb:sc:introduction}
constructs the selected critical profile and proves its quantitative
coercivity estimate. Appendix~\ref{sb:repro:appendix} describes the
computer-assisted verification and gives the location of its code and data.

\section{Periodic conventions and dynamic rescaling}
\label{sb:common:framework}

This section fixes the rescaling used in the three profile
constructions. We first derive the normalized equation and its two origin
conditions, then record the reconstruction lemma that converts normalized
convergence into physical blowup.

Throughout the paper the physical circle is
\(\mathbb T=\mathbb R/(2\pi\mathbb Z)\), and the equation is
\begin{equation}\label{sb:common:pde}
 \omega_t+a u\omega_y=u_y\omega,\qquad
 u_y=H_{\mathbb T}\omega,\qquad u(0,t)=0.
\end{equation}
We work with odd vorticity, so that its mean is zero and its normalized
velocity is odd. Our sign convention is
\begin{equation}\label{sb:common:hilbert}
 H_{\mathbb T}v(y)=\frac1{2\pi}\operatorname{p.v.}
  \int_{-\pi}^{\pi}v(\eta)\cot\frac{y-\eta}{2}\,\idiff\eta,
 \qquad
 Hv(x)=\frac1\pi\operatorname{p.v.}\int_{\mathbb R}
                  \frac{v(\eta)}{x-\eta}\,\idiff\eta .
\end{equation}
On a rescaled cell of half-period \(L\), set
\[
 H_Lv(x)=\frac1{2L}\operatorname{p.v.}\int_{-L}^{L}
        v(\eta)\cot\frac{\pi(x-\eta)}{2L}\,\idiff\eta,
 \qquad H_\infty=H.
\]
The parameter \(a\) remains fixed during every physical evolution.

For positive functions \(C_l,C_\omega\), define
\begin{equation}\label{sb:common:rescaling}
 \Omega(x,\tau)=C_\omega(\tau)^{-1}
       \omega(C_l(\tau)^{-1}x,t(\tau)),\qquad
 L=\pi C_l,\qquad
 t'=C_\omega^{-1},\quad
 \frac{C_l'}{C_l}=c_l,\quad
 \frac{C_\omega'}{C_\omega}=-c_\omega .
\end{equation}
Primes on the scale functions denote differentiation in \(\tau\).
For profiles and one-variable kernels, primes denote spatial derivatives.
The rescaled velocity is
\[
 U(x,\tau)=C_lC_\omega^{-1}u(C_l^{-1}x,t(\tau))
          =\int_0^xH_L\Omega(\eta,\tau)\,\idiff\eta.
\]
Direct differentiation gives
\begin{equation}\label{sb:common:normalized}
 \Omega_\tau+(c_lx+aU)\Omega_x
                  =(c_\omega+H_L\Omega)\Omega .
\end{equation}
The changing cell in this identity is obtained from a solution of
\eqref{sb:common:pde} on a fixed physical circle.

\begin{lemma}[Origin normalization]\label{sb:common:gauges}
For a smooth odd solution, the conditions
\(\Omega_x(0)=-1\), \(\Omega_{xxx}(0)=6\) are preserved by
\begin{align}
 d_L(v)&=-\frac{3-a}{12}(H_Lv)''(0),\notag\\
 \ell_{l,L}(v)&=d_L(v)-aH_Lv(0),\qquad
 \ell_{\omega,L}(v)=d_L(v)-H_Lv(0),\label{sb:common:two-gauges}\\
 c_l&=\ell_{l,L}(\Omega),\qquad
 c_\omega=\ell_{\omega,L}(\Omega).\notag
\end{align}
At fixed \(L\), one can instead set \(c_l=0\) and preserve just
\(\Omega_x(0)=-1\) by choosing
\begin{equation}\label{sb:common:one-gauge}
 c_\omega=-(1-a)H_L\Omega(0).
\end{equation}
\end{lemma}
\begin{proof}
Writing \(b_\Omega=H_L\Omega(0)\) and
\(j_\Omega=(H_L\Omega)''(0)\), parity and
\eqref{sb:common:normalized} give
\[
 \partial_\tau\Omega_x(0)
 =\bigl(c_\omega-c_l+(1-a)b_\Omega\bigr)\Omega_x(0)
\]
and
\[
 \partial_\tau\Omega_{xxx}(0)
 =\bigl(c_\omega-3c_l+(1-3a)b_\Omega\bigr)\Omega_{xxx}(0)
                    +(3-a)j_\Omega\Omega_x(0).
\]
Substituting the prescribed derivatives proves both formulas.
\end{proof}

A stationary vorticity profile \(\bar\Omega\), its velocity \(\bar U\), and its constant
rates satisfy
\begin{equation}\label{sb:common:stationary}
 (\bar c_lx+a\bar U)\bar\Omega'=(\bar c_\omega+\bar U')\bar\Omega,
 \qquad \bar U'=H_{\bar L}\bar\Omega,\qquad \bar U(0)=0.
\end{equation}
The amplitude rate satisfies \(\bar c_\omega<0\). The three geometries are
\(\bar L=\infty,\bar c_l>0\) (focusing),
\(\bar L=\infty,\bar c_l=0\) (critical), and
\(\bar L<\infty,\bar c_l=0\) (finite period).
In the last case the spatial scale is fixed throughout the construction.

For the chosen normalization, let \(\ell_l,\ell_\omega\) denote its
linear rate functionals. Define the perturbation velocity and the
transport and stretching coefficients by
\[
 \begin{gathered}
 U_h(x)=\int_0^xH_{\bar L}h(y)\,\idiff y,\\
 B=\bar c_lx+a\bar U,\qquad Q=\bar c_\omega+\bar U'.
 \end{gathered}
\]
The exact linear and quadratic expressions are
\begin{align}
 \mathcal Lh={}&-Bh'+Qh+(H_{\bar L}h+\ell_\omega(h))\bar\Omega
                    -(aU_h+\ell_l(h)x)\bar\Omega',\label{sb:common:linear}\\
 \mathcal N(h)={}&(H_{\bar L}h+\ell_\omega(h))h
                    -(aU_h+\ell_l(h)x)h'.\label{sb:common:quadratic}
\end{align}
We also write \(f=-\bar\Omega/x\) and
\(b=H_{\bar L}\bar\Omega(0)\). The origin equation gives
\begin{equation}\label{sb:common:rate-normalization}
 \bar c_\omega=\bar c_l+(a-1)b.
\end{equation}
The origin normalization does not impose \(\bar c_\omega=-1\).
The alternative tuple
\[
 \left(\frac{\bar\Omega}{-\bar c_\omega},
       \frac{\bar U}{-\bar c_\omega},
       \frac{\bar c_l}{-\bar c_\omega},-1\right)
\]
has unit amplitude decay rate and the same half-period, but different
origin coefficients. The physical focusing exponent is
\(-\bar c_l/\bar c_\omega\) in either convention. At critical or finite-period
profiles, \(\bar c_l=0\), so \(\bar c_\omega=-(1-a)b\).
The tangent conditions are \(h'(0)=h'''(0)=0\) for the whole-line
profiles and \(h'(0)=0\) at fixed finite period. Each case below specifies
the space, closed operator domain, and estimates used for these expressions.
The estimates below provide the dynamical control needed to transfer
the stationary profiles to smooth periodic solutions.

\begin{lemma}[Physical reconstruction]\label{sb:common:physical}
Suppose a normalized solution is smooth for every finite \(\tau\),
\(c_\omega(\tau)\to\bar c_\omega<0\), and the rescaled vorticity converges
uniformly to a nonzero profile \(\bar\Omega\). For whole-line profiles, uniform
convergence here refers to the zero extension of one period; convergence
of the actual periodic composition is asserted on compact intervals.
Then the terminal physical time \(T\) is finite and
\begin{equation}\label{sb:common:type-I}
 (T-t)C_\omega\longrightarrow(-\bar c_\omega)^{-1},\qquad
 (T-t)\|\omega(t)\|_\infty\longrightarrow(-\bar c_\omega)^{-1}\|\bar\Omega\|_\infty.
\end{equation}
If additionally
\(c_l-\bar c_l,c_\omega-\bar c_\omega\in L^1(0,\infty)\), then for some
\(A>0\),
\begin{equation}\label{sb:common:power-law}
 C_l(t)^{-1}\sim A(T-t)^{-\bar c_l/\bar c_\omega}.
\end{equation}
For the fixed-period construction, \(C_l=\bar L/\pi\) exactly.
\end{lemma}
\begin{proof}
Set \(A_\omega=C_\omega^{-1}\); then
\(A_\omega'/A_\omega=c_\omega\to\bar c_\omega\). Consequently
\(T=t(0)+\int_0^\infty A_\omega(\tau)\,\idiff\tau<\infty\).
For every sufficiently small \(\epsilon>0\), the ratio
\(A_\omega(\tau+r)/A_\omega(\tau)\) lies between
\(e^{(\bar c_\omega-\epsilon)r}\) and
\(e^{(\bar c_\omega+\epsilon)r}\) for all \(r\geq0\) and large \(\tau\).
Integration and then \(\epsilon\downarrow0\) prove the first limit.
Use \(\omega(C_l^{-1}x,t)=C_\omega\Omega(x,\tau)\) for the second.
Under the integrable-error hypothesis,
\(C_\omega\sim A_1e^{-\bar c_\omega\tau}\) and
\(C_l\sim A_2e^{\bar c_l\tau}\), which imply
\eqref{sb:common:power-law}. The nonzero amplitude limit excludes smooth
continuation through \(T\). Smoothness at every earlier physical time
follows from the assumed normalized continuation and the finite positive
scales there.
\end{proof}

\section{Local solutions and continuation}\label{sb:gwp:local-section}

We establish the \(H^1\) solution theory and a bounded-vorticity
continuation criterion. The derivative equation provides both the local
energy estimate and continuous dependence at this endpoint. A logarithmic
estimate for the Hilbert transform then reduces continuation to control
of \(\|\omega\|_\infty\).

In this section and the global regularity argument, \(x\) is the
physical coordinate on \(\mathbb T=\mathbb R/(2\pi\mathbb Z)\).
We use
\[
 H_{\mathbb T}f(x)=\frac1{2\pi}\operatorname{p.v.}
   \int_{-\pi}^{\pi}f(y)\cot\frac{x-y}{2}\,\mathrm{d}y,
 \qquad \widehat{H_{\mathbb T}f}(n)=-i\operatorname{sgn}(n)\widehat f(n).
\]
The velocity is fixed uniquely by \(u_x=H_{\mathbb T}\omega\) and
\(u(0)=0\). In particular,
\begin{equation}\label{sb:gwp:velocity-bound}
 \|u\|_{H^{s+1}}\leq C_s\|\omega\|_{H^s}\qquad(s\geq0).
\end{equation}
Indeed, the nonzero Fourier coefficients of \(u\) are
\(-\widehat\omega(n)/|n|\), and the value condition controls its
constant coefficient by Cauchy--Schwarz.

\begin{proposition}[Local solutions]\label{sb:gwp:local}
For every real \(a\) and real \(\omega_0\in H^1(\mathbb T)\),
the equation
\begin{equation}\label{sb:gwp:equation}
 \omega_t+a u\omega_x=u_x\omega,
 \qquad u_x=H_{\mathbb T}\omega,\qquad u(0)=0
\end{equation}
has a unique maximal solution
\[
 \omega\in C([0,T_*);H^1)\cap C^1([0,T_*);L^2).
\]
Its mean is constant. Odd initial data give odd solutions and odd
velocities; in that case \(u(t,\pi)=0\).
The solution depends continuously on its initial data in
\(C([0,T];H^1)\) for every \(T<T_*\).
If \(\omega_0\in H^s\), \(s>1\), then this regularity persists
on every compact subinterval of \([0,T_*)\).
\end{proposition}

\begin{proof}
The Hilbert transform is bounded on every \(H^s\), and
\(H^1(\mathbb T)\) embeds into \(C^{0,1/2}\).
Writing \(v=\omega_x\), we obtain the differentiated equation
\begin{equation}\label{sb:gwp:derivative-equation}
 v_t+a u v_x=(1-a)u_xv+\omega H_{\mathbb T}v.
\end{equation}
Multiplication by \(\omega\) and \(v\), respectively, gives
\begin{align}
 \frac12\frac{\mathrm{d}}{\mathrm{d}t}\|\omega\|_2^2
   &=(1+a/2)\int_{\mathbb T}u_x\omega^2\,\mathrm{d}x,
                                                        \label{sb:gwp:L2-energy}\\
 \frac12\frac{\mathrm{d}}{\mathrm{d}t}\|v\|_2^2
   &=(1-a/2)\int_{\mathbb T}u_xv^2\,\mathrm{d}x
       +\int_{\mathbb T}\omega vH_{\mathbb T}v\,\mathrm{d}x.
                                                        \label{sb:gwp:derivative-energy}
\end{align}
Consequently
\begin{equation}\label{sb:gwp:H1-energy}
 \frac{\mathrm{d}}{\mathrm{d}t}\|\omega\|_{H^1}
 \leq C_a(\|H_{\mathbb T}\omega\|_\infty+\|\omega\|_\infty)
                      \|\omega\|_{H^1}
 \leq C_a\|\omega\|_{H^1}^2.
\end{equation}

For clarity, these estimates also construct the solution. Project the
equation onto Fourier modes \(|n|\leq N\). The resulting finite
ordinary differential equation obeys the same energy bounds, because
the projection is orthogonal and commutes with differentiation.
On a time interval depending only on \(a\) and
\(\|\omega_0\|_{H^1}\), the approximations are uniformly bounded
in \(L^\infty H^1\), and their time derivatives in \(L^\infty L^2\).
Fourier truncation and the Arzel\`a--Ascoli theorem give a subsequence
converging in \(C L^2\). Interpolation gives convergence in
\(C H^\sigma\) for every \(\sigma<1\), hence uniform convergence
of \(\omega\) and \(u_x\) when \(\sigma>1/2\).
The limit solves \eqref{sb:gwp:equation} distributionally.

The limit actually lies in \(C H^1\). Its derivative \(v\in
L^\infty L^2\) solves \eqref{sb:gwp:derivative-equation} with
an \(L^\infty L^2\) right-hand side and velocity \(u\in C C^1\).
Changing variables by its \(C^1\) flow gives the Duhamel formula
in \(L^2\), and therefore \(v\in C L^2\), including at the initial
time. This change of variables is valid for distributional solutions:
spatial mollification gives a commutator bounded by
\(C\|u_x\|_\infty\|v\|_2\), which tends to zero by density of
smooth functions. The equation now gives \(\omega_t\in C L^2\).

For uniqueness, subtract two solutions, using the first velocity in
the transport term. If \(w\) is their difference, then
\(\|u_1-u_2\|_\infty\leq C\|w\|_2\), and multiplication by
\(w\) gives
\[
 \frac{\mathrm{d}}{\mathrm{d}t}\|w\|_2^2
 \leq C_a(\|\omega_1\|_{H^1}+\|\omega_2\|_{H^1})\|w\|_2^2.
\]
Gronwall's inequality proves uniqueness. Restarting the construction
defines the maximal interval. Skew-adjointness of \(H_{\mathbb T}\)
shows \(\frac{\mathrm{d}}{\mathrm{d}t}\int_{\mathbb T}\omega\,\idiff x
=(1+a)\int_{\mathbb T}\omega H_{\mathbb T}\omega\,\idiff x=0\).
Reflection and uniqueness prove oddness; an odd periodic velocity
vanishes at both \(0\) and \(\pi\).

We record continuous dependence at the endpoint \(H^1\).
For convergent \(H^1\) data, the uniform local bound and the difference
estimate give \(\omega_n\to\omega\) in \(C L^2\), hence uniformly,
and \(u_n\to u\) in \(C C^1\). Their flows \(\eta_n\) converge
in \(C C^1\); the pullbacks \(P_n(t)v=v\circ\eta_n(t)\) and their
inverses are uniformly bounded on \(L^2\) and converge strongly,
uniformly in time, as follows first on smooth functions and then by density.
For \(V_n=P_n\omega_{nx}\) and \(V=P\omega_x\), where \(P(t)v=v\circ\eta(t)\), the derivative equation becomes the
bounded \(L^2\) ordinary differential equation
\[
 (V_n)_t=\mathcal A_n(t)V_n,\qquad
 \mathcal A_n=(1-a)M_{u_{nx}\circ\eta_n}
               +M_{\omega_n\circ\eta_n}P_nH_{\mathbb T}P_n^{-1},
\]
where \(M_b\) denotes multiplication by \(b\).
These operators are uniformly bounded and converge strongly to
\(\mathcal A(t)\), uniformly in time on each fixed vector.
This convergence is uniform on the compact \(L^2\) set
\(\{P(t)\omega_x(t):0\leq t\leq T\}\).
The integral equation and Gronwall therefore give \(V_n\to V\)
in \(C L^2\), and the inverse pullbacks give
\(\omega_{nx}\to\omega_x\) in \(C L^2\).
This proves local continuous dependence in \(H^1\); finitely many
restarts give it on every compact interval of the limiting solution.

Here is also a direct higher regularity estimate, to specify the
regularity assertion for noninteger \(s\):
\begin{equation}\label{sb:gwp:higher-energy}
 \frac{\mathrm{d}}{\mathrm{d}t}\|\omega\|_{H^s}
       \leq C_{a,s}\|\omega\|_{H^1}\|\omega\|_{H^s},\qquad s>1.
\end{equation}
The product term has this bound by convolution of Fourier
coefficients and \(\sum_n|\widehat\omega(n)|\leq C\|\omega\|_{H^1}\).
For the transport commutator, with \(\langle n\rangle=(1+n^2)^{1/2}\),
\[
 \frac{|k|}{|\ell|}
  |\langle k+\ell\rangle^s-\langle k\rangle^s|
 \leq C_s(\langle k\rangle^s+\langle\ell\rangle^s),\qquad\ell\ne0.
\]
The mean-value theorem and Young's inequality prove this scalar
inequality. Using \(|\widehat u(\ell)|=|\widehat\omega(\ell)|/|\ell|\),
the same convolution estimate bounds the commutator; the remaining
transport term is integrated by parts. This proves
\eqref{sb:gwp:higher-energy}, first on the approximations and then
by passage to the limit.
\end{proof}

\begin{proposition}[Continuation criterion]\label{sb:gwp:continuation}
If \(T_*<\infty\) and
\(\sup_{t<T_*}\|\omega(t)\|_\infty<\infty\), then the solution
extends beyond \(T_*\).
\end{proposition}

\begin{proof}
The elementary logarithmic estimate
\begin{equation}\label{sb:gwp:log-Hilbert}
 \|H_{\mathbb T}f\|_\infty
 \leq C\{1+\|f\|_\infty\log(e+\|f\|_{H^1})\}
\end{equation}
follows by subtracting \(f(x)\) in the principal value integral.
On \(|y|<\delta\), use
\(|f(x-y)-f(x)|\leq\|f'\|_2|y|^{1/2}\). Together with the
complementary integral this gives the bound
\[
 C\|f'\|_2\delta^{1/2}+C\|f\|_\infty\log(2/\delta).
\]
Choose \(\delta=(e+\|f\|_{H^1})^{-2}\).
Combining \eqref{sb:gwp:log-Hilbert} with
\eqref{sb:gwp:H1-energy}, the quantity
\(Y=\log(e+\|\omega\|_{H^1})\) satisfies
\(Y'\leq C_{a,M}(1+Y)\) if \(\|\omega\|_\infty\leq M\).
Thus \(\|\omega\|_{H^1}\) stays bounded up to \(T_*\), and the
uniform local existence interval from Proposition~\ref{sb:gwp:local}
permits a restart across \(T_*\).
\end{proof}

\section{Endpoint conditions for global regularity and blowup}
\label{sb:gwp:section}

The global argument has two components. A corrected moment controls the
accumulated velocity gradient at the origin, while a positive Lagrangian
force controls vorticity amplification along the remaining trajectories.
The endpoint and interior-zero hypotheses are stated below.
After completing this argument, we use the same flow and positive
velocity kernel to prove the inverse-vorticity blowup criterion in
Subsection~\ref{sb:inverse:subsection}.
Throughout this section the physical half-period is \(\pi\), and
the equation and Hilbert transform are those of
Section~\ref{sb:gwp:local-section}. A period-\(\pi\) formulation is
converted to these conventions by
\(\Omega(x,t)=\omega(x/2,t)\), \(U(x,t)=2u(x/2,t)\):
\(U_x=H_{\mathbb T}\Omega\) and
\(\Omega_t+aU\Omega_x=U_x\Omega\), with unchanged time and parameter.

\begin{theorem}[Global solutions]\label{sb:gwp:theorem}
Let \(a\geq1/2\), and put \(\beta_c=1+1/a\).
Suppose \(\omega_0\in H^1(\mathbb T)\) is real, odd, and
\(\omega_0\leq0\) on \([0,\pi]\). At every interior zero assume
\begin{equation}\label{sb:gwp:CZ}
 H_{\mathbb T}\omega_0(z)\leq0
 \qquad(0<z<\pi,\ \omega_0(z)=0).\qquad\text{(CZ)}
\end{equation}
Assume in addition at least one of the two endpoint conditions
\begin{align}
 m_a(\omega_0)&:=\operatorname*{ess\,sup}_{0<x<\pi}
                 \frac{|\omega_0(x)|}{x^{1/a}}<\infty,
                                                        \label{sb:gwp:trace}\\
 Q_{\beta_c}(0)&:=\int_0^\pi|\omega_0(x)|
                           \cot^{\beta_c}(x/2)\,\mathrm{d}x<\infty.
                                                        \label{sb:gwp:weighted-endpoint}
\end{align}
Then the maximal solution from Proposition~\ref{sb:gwp:local} is
global, remains odd, and preserves its sign on \([0,\pi]\).
It belongs to \(C([0,\infty);H^1)\cap C^1([0,\infty);L^2)\).
Every additional \(H^s\) regularity, \(s>1\), persists on finite
time intervals. In particular, smooth data give a global smooth solution.
\end{theorem}

The point values in \eqref{sb:gwp:CZ} are well-defined because
\(H_{\mathbb T}\omega_0\in H^1\subset C^0\).
Strict negativity on \((0,\pi)\) makes \eqref{sb:gwp:CZ} vacuous.
For smooth odd data, \eqref{sb:gwp:trace} is automatic if \(a\geq1\);
if \(1/2\leq a<1\), it holds when \(\omega_0'(0)=0\).

\begin{remark}[The De Gregorio comparison]\label{sb:gwp:Chen-comparison}
At \(a=1\), Chen's theorem~\cite[Theorem~2]{Chen2025Regularity}
assumes the same \(H^1\), oddness, and sign conditions, together with
\(\omega_0/x\in L^\infty\) and
\[
 \int_0^\pi\frac{|\omega_{0x}(x)|^2}{|\omega_0(x)|}\sin x\,\mathrm{d}x<\infty.
\]
It does not impose \eqref{sb:gwp:CZ}. The present proof instead uses
the stated zero-set compatibility and does not require this additional
weighted derivative integral. The hypotheses of the two results are
therefore different.
\end{remark}

\subsection{Flow identities and endpoint control}

For \(a>0\), let \(\eta_t=a u(t,\eta)\), \(\eta(0,\theta)=\theta\), and
\(J=\eta_\theta\). On each compact interval of existence,
\(\eta\) is a \(C^1\) diffeomorphism and
\begin{equation}\label{sb:gwp:flow-identity}
 J_t=a u_x(t,\eta)J,\qquad
 \omega(t,\eta(t,\theta))=\omega_0(\theta)J(t,\theta)^{1/a}.
\end{equation}
Thus signs and zeros are transported. Oddness fixes both endpoints.
For odd \(\omega\) nonpositive on \((0,\pi)\), integrating the odd Hilbert
kernel gives
\begin{equation}\label{sb:gwp:positive-velocity}
 u(x)=-\frac1\pi\int_0^\pi\omega(y)
        \log\left|\frac{\sin((x+y)/2)}{\sin((x-y)/2)}\right|\,\mathrm{d}y
       \geq0,\qquad 0\leq x\leq\pi.
\end{equation}
Indeed, the difference of the squared sines in this ratio is
\(\sin x\sin y\geq0\). Since \(a>0\), all interior trajectories
move away from the origin.

The endpoint quotient has a useful cancellation:
\[
 (\partial_t+a u\partial_x)(\omega/x^{1/a})
       =(u_x-u/x)(\omega/x^{1/a}),\qquad
 |u_x(x)-u(x)/x|\leq Cx^{1/2}\|\omega\|_{H^1}.
\]
The second estimate follows by averaging
\(u_x(x)-u_x(y)\) over \(0<y<x\) and using
\(\|u_{xx}\|_2\leq\|\omega_x\|_2\).
Consequently, under \eqref{sb:gwp:trace},
\begin{equation}\label{sb:gwp:trace-propagation}
 |\omega(t,x)|\leq m_a(\omega_0)x^{1/a}
     \exp\left(Cx^{1/2}\int_0^t\|\omega(s)\|_{H^1}\,\mathrm{d}s\right).
\end{equation}
Here monotonicity bounds every earlier point of the trajectory by
its terminal coordinate \(x\).

We will also use a qualitative fact requiring no pointwise trace.
For each fixed \(b>0\), finiteness of
\(\int|\omega_0|\cot^b(x/2)\,\mathrm{d}x\) propagates on compact regular
time intervals. Indeed, \(J\) is bounded above and below there,
\(\eta(t,\theta)\) is comparable to \(\theta\), and
\(\pi-\eta(t,\theta)\) is comparable to \(\pi-\theta\).
Apply \eqref{sb:gwp:flow-identity} and change variables in the integral.
The same argument gives uniform integrability of these weighted
densities on compact regular time intervals.

\subsection{A positive quadratic form}

We first prove the positivity needed for the corrected moment estimate in
Lemma~\ref{sb:gwp:corrected-moment}.

\begin{lemma}\label{sb:gwp:positive-form}
Put
\[
 \ell(s)=\log\left|\frac{s+1}{s-1}\right|,
 \qquad \Phi_m(s)=(s^m+1)\ell(s)
                     -(s^m-1)\frac{2s}{s^2-1}.
\]
For \(0<m\leq4\), the quadratic form
\begin{equation}\label{sb:gwp:positive-form-equation}
 \int_0^\pi\!\int_0^\pi
 v(x)v(y)\cot^m(y/2)
        \Phi_m\!\left(\frac{\cot(x/2)}{\cot(y/2)}\right)\mathrm{d}x\,\mathrm{d}y
 \geq0
\end{equation}
whenever the logarithmic transform used below belongs to \(L^2\).
\end{lemma}

\begin{proof}
Set \(x=2\arctan e^r\) and
\[
 G_m(r)=\frac{2v(2\arctan e^r)e^{(1-m/2)r}}{1+e^{2r}},\qquad
 W_m(z)=2\cosh(mz/2)\log|\coth(z/2)|
                      -2\frac{\sinh(mz/2)}{\sinh z}.
\]
The left side of \eqref{sb:gwp:positive-form-equation} is
\(\int G_m(r)G_m(t)W_m(t-r)\,\mathrm{d}r\,\mathrm{d}t\).
Here \(W_m\) is even and integrable. We give its transform to make
the positivity assertion self-contained. Write \(m=2q\), \(0<q\leq2\).
For \(z>0\), the geometric series for the logarithm and \(1/\sinh z\)
give
\[
 W_{2q}(z)=4\sum_{r=1,3,5,\ldots}e^{-rz}
                      \{r^{-1}\cosh(qz)-\sinh(qz)\}.
\]
The \(r=1\) summand is \(4e^{-(1+q)z}\). For \(r\geq3\), the
\(L^1(0,\infty)\) norm of the displayed summand is at most
\(4(1+q)/(r^2-q^2)\). Thus the paired series converges in \(L^1\)
and can be transformed term by term.
For \(\zeta=q+i\xi\), the paired odd-integer sum is
\[
 \widehat W_{2q}(\xi)
 =4\operatorname{Re}\sum_{n\in2\mathbb Z+1}
                  \frac{1/n-1}{n-\zeta}
 =2\pi\operatorname{Re}
       \left((\zeta^{-1}-1)\tan(\pi\zeta/2)\right).
\]
In the last equality we used
\(\operatorname{p.v.}\sum_{n\in2\mathbb Z+1}(n-\zeta)^{-1}
=(\pi/2)\tan(\pi\zeta/2)\), the partial-fraction expansion of
the tangent, and the vanishing paired sum of \(1/n\).
Equivalently, this partial-fraction identity follows by applying
the residue theorem to the tangent on expanding rectangles and
pairing the poles at positive and negative odd integers.
Taking the real part yields
\begin{equation}\label{sb:gwp:positive-symbol}
 \widehat W_{2q}(\xi)=
 \frac{2\pi\{\xi\sinh(\pi\xi)
              -(\xi^2+q^2-q)\sin(\pi q)\}}
      {(\xi^2+q^2)(\cosh(\pi\xi)+\cos(\pi q))}.
\end{equation}
The removable value at \((q,\xi)=(1,0)\) is \(4\).
For \(0<q<1\), its numerator divided by \(2\pi\) is at least
\(\xi^2(\pi-\sin\pi q)+q(1-q)\sin\pi q>0\).
For \(1\leq q\leq2\), both terms are nonnegative because
\(\sin\pi q\leq0\) and \(q^2-q\geq0\).
The denominator is positive away from the removable point.
Plancherel's identity now proves the lemma, by approximation for
arbitrary \(G_m\in L^2\).
\end{proof}

\subsection{Corrected moments below the endpoint}

\begin{lemma}\label{sb:gwp:corrected-moment}
Fix \(a\geq1/2\), and put \(\beta_*=1+1/(2a)\).
There exist \(C_a,\kappa_a>0\), independent of
\(\beta\in[\beta_*,\beta_c)\), such that every regular odd
nonpositive solution with finite \(Q_\beta\) obeys
\begin{equation}\label{sb:gwp:corrected-inequality}
 \mathcal E_\beta'\leq C_a\varepsilon\mathcal E_\beta^2,
 \qquad \mathcal E_\beta=Q_\beta+\kappa_a\varepsilon M,
 \qquad \varepsilon=\beta_c-\beta,
\end{equation}
where \(Q_\beta=\int_0^\pi|\omega|\cot^\beta(x/2)\,\mathrm{d}x\)
and \(M=\int_0^\pi|\omega|\,\mathrm{d}x\).
\end{lemma}

\begin{proof}
For \(A=\cot(x/2)\), \(B=\cot(y/2)\), integration by parts
and symmetrization give
\begin{align}
 -Q_\beta'
 &=\frac1{2\pi}\sum_{m\in\{\beta-1,\beta+1\}}
       \int_0^\pi\!\int_0^\pi
       \omega(x)\omega(y)B^m P_m(A/B)\,\mathrm{d}x\,\mathrm{d}y,
                                                    \label{sb:gwp:symmetry}\\
 P_m(s)&=\frac{a\beta}{2}(s^m+1)\ell(s)
                 -\frac{1+a}{2}(s^m-1)\frac{2s}{s^2-1}.
 \notag
\end{align}
For example, this identity follows immediately from
\(A'=-(1+A^2)/2\),
\(\partial_x\ell(A/B)=B(1+A^2)/(A^2-B^2)\), and
\eqref{sb:gwp:positive-velocity}.
Since \(1+a-a\beta=a\varepsilon\),
\[
 P_m=\frac{a\beta}{2}\Phi_m
             -\frac{a\varepsilon}{2}\Psi_m,
 \qquad \Psi_m(s)=(s^m-1)\frac{2s}{s^2-1}.
\]
The parameter restriction gives \(0<\beta-1<\beta+1<4\).
Both \(\Phi_m\) forms are nonnegative by
Lemma~\ref{sb:gwp:positive-form}. Its \(L^2\) condition follows
from \(Q_\beta<\infty\) and boundedness of \(\omega\): directly,
\(\|G_m\|_2^2\leq C\int|\omega|^2(1+\cot^\beta(x/2))\,\mathrm{d}x\).
It follows that
\[
 Q_\beta'\leq C_a\varepsilon\iint|\omega(x)\omega(y)|K_\beta(A,B)\,\mathrm{d}x\,\mathrm{d}y,
 \quad K_\beta(A,B)=2AB\frac{A^{\beta+1}-B^{\beta+1}
                                  +A^{\beta-1}-B^{\beta-1}}{A^2-B^2}.
\]
The continuous value is used at \(A=B\). Uniformly in the stated
\(\beta\) interval,
\begin{equation}\label{sb:gwp:kernel-bound}
 K_\beta(A,B)\leq C_a\left(A^\beta B+AB^\beta
       +\mathbf1_{AB<1}\frac{1-AB}{A+B}\right).
\end{equation}
To verify this, suppose \(A=rB\geq B\). The two quotient terms
are at most \(C_a(A^\beta B+AB^\beta)\) and
\(C_aA^{\beta-2}B\), respectively. If \(A\geq1\), the latter
is bounded by \(C_aA^\beta B\). If \(A<1\) and \(AB\leq1/2\),
it is at most \(C_a\), while \((1-AB)/(A+B)\geq1/4\).
In the remaining case \(1/2<AB<1\), both \(A,B\) are bounded
away from zero and the first two terms suffice.

Define
\[
 h=\int_0^\pi|\omega|\cot(x/2)\,\mathrm{d}x,\qquad
 D=\iint_{x+y>\pi}|\omega(x)\omega(y)|
                            \bigl(-\cot((x+y)/2)\bigr)\,\mathrm{d}x\,\mathrm{d}y.
\]
The addition formula \(\cot((x+y)/2)=(AB-1)/(A+B)\) turns
\eqref{sb:gwp:kernel-bound} into
\(Q_\beta'\leq C_a\varepsilon(hQ_\beta+D)\).
We claim \(hQ_\beta\leq C_a(Q_\beta^2+D)\).
H\"older gives \(h\leq Q_\beta^{1/\beta}M^{1-1/\beta}\).
The claim is immediate if \(h\leq Q_\beta\); otherwise
\(Q_\beta<M\) and \(h<M\).
The part \(h_n\) with \(A\geq1\) is at most \(Q_\beta\).
For \(0<A<1\), \(B>0\), one has
\[
 A\leq C_a\left((AB)^\beta
                 +\mathbf1_{AB<1}\frac{1-AB}{A+B}\right).
\]
Indeed, use the first term if \(AB\geq1/2\); if \(AB\leq1/2\),
the second is at least \(A/3\). Integration gives
\(h_f M\leq C_a(Q_\beta^2+D)\), proving the claim.
Thus
\begin{equation}\label{sb:gwp:Q-estimate}
 Q_\beta'\leq C_a\varepsilon(Q_\beta^2+D).
\end{equation}

The mass supplies precisely the missing negative term:
\[
 M'=\frac{1+a}{2\pi}\iint|\omega(x)\omega(y)|
                  \cot((x+y)/2)\,\mathrm{d}x\,\mathrm{d}y
       \leq C_aQ_\beta^2-c_aD.
\]
Here \((AB-1)/(A+B)\leq(AB)^\beta\) for \(AB\geq1\), and
the other contribution is \(-D\). Add \(\kappa_a\varepsilon\)
times this inequality to \eqref{sb:gwp:Q-estimate}, choosing
\(\kappa_a\) to absorb \(D\), to obtain
\eqref{sb:gwp:corrected-inequality}.

We specify the endpoint justification at \(H^1\) regularity.
All quantities above are finite on compact regular time intervals:
\(h<\infty\) by the \(1/2\)-H\"older endpoint bound, and
\(D<\infty\) since its only corner singularity is bounded by
\(C\sqrt{\pi-x}\sqrt{\pi-y}/(2\pi-x-y)\).
Insert smooth cutoffs near the two endpoints in the moment identity.
The origin cutoff error is bounded by a constant times
\(\int_0^{2\delta}|\omega|\cot^\beta(x/2)\,\mathrm{d}x\), because
\(|u(x)|\leq Cx\). This tends uniformly to zero on compact regular
time intervals by the flow-based uniform integrability proved above.
The other endpoint is easier. The velocity logarithm is nonnegative,
so its double integral is absolutely integrable by the same moment
bound. For the Hilbert term, diagonal truncation converges uniformly
on \(H^1\) inputs, with error at most \(C\delta^{1/2}\|\omega_x\|_2\).
After symmetrization, \eqref{sb:gwp:kernel-bound} dominates the
remaining integrals. The \(\Phi_m\) terms also converge by
\(G_m\in L^2\), \(W_m\in L^1\).
These observations justify \eqref{sb:gwp:symmetry} and the locally
absolutely continuous moment inequalities, including under
\eqref{sb:gwp:weighted-endpoint} without a pointwise trace assumption.
\end{proof}

\subsection{A positive force and vorticity control}

\begin{lemma}\label{sb:gwp:force}
For real \(v\in H^2(\mathbb T)\), the distribution
\(\mathscr F_v=-vv_{xx}-H_{\mathbb T}(vH_{\mathbb T}v_{xx})\)
has a continuous nonnegative representative. More precisely, if
\(v=\sum c_ne^{inx}\) and \(\Phi_k=\sum_{n\geq k}c_ne^{inx}\), then
\begin{equation}\label{sb:gwp:force-squares}
 \mathscr F_v=2\sum_{k\geq1}(2k-1)|\Phi_k|^2\geq0.
\end{equation}
For the solution, set \(\varrho=J^{-1/(2a)}\), a negative power of the
flow Jacobian. Then
\begin{equation}\label{sb:gwp:rho-data}
 \varrho(0,\theta)=1,\qquad
 \varrho_t(0,\theta)=-\tfrac12u_{0x}(\theta),\qquad
 \omega(t,\eta(t,\theta))=\frac{\omega_0(\theta)}{\varrho(t,\theta)^2},
\end{equation}
and the transformed Jacobian satisfies
\begin{equation}\label{sb:gwp:rho-equation}
 \varrho_{tt}=\frac{\omega_0(\theta)^2}{4\varrho^3}
             +\frac a2\mathscr F_u(t,\eta)\varrho\geq0.
\end{equation}
\end{lemma}

\begin{proof}
For a trigonometric polynomial, the zero Fourier coefficient of
\(\mathscr F_v\) is \(2\sum_{n\geq1}n^2|c_n|^2\), and its
coefficient at frequency \(j>0\) is
\(2\sum_{n\geq1}n^2c_{n+j}\overline{c_n}\).
These are the coefficients of \eqref{sb:gwp:force-squares}, since
\(\sum_{k=1}^n(2k-1)=n^2\).
For \(v\in H^2\), Cauchy--Schwarz gives
\(|\Phi_k|\leq Ck^{-3/2}\|v\|_{H^2}\). The squared series
converges uniformly, and the polynomial approximations converge
both distributionally and uniformly in this identity.

Oddness gives mean zero for \(\omega\). The periodic Cotlar identity is
therefore
\[
 H_{\mathbb T}(\omega H_{\mathbb T}\omega)
       =\tfrac12\{(H_{\mathbb T}\omega)^2-\omega^2\}.
\]
To see this directly, \(\omega+iH_{\mathbb T}\omega\) has only
positive Fourier modes; apply the Hilbert transform relation between
the real and imaginary parts of its square. For nonzero mean the
usual constant correction would be required; it is absent here.
Since \(H_{\mathbb T}u_{xx}=-\omega_x\),
\(H_{\mathbb T}(u\omega_x)=\mathscr F_u+uu_{xx}\). Thus
\[
 (\partial_t+a u\partial_x)u_x
       =\tfrac12(u_x^2-\omega^2)-a\mathscr F_u.
\]
Its right-hand side is continuous. The transport Duhamel formula
therefore gives this identity pointwise along the flow even for
\(H^1\) solutions. Differentiate
\(\varrho_t=-u_x(t,\eta)\varrho/2\) to obtain
\eqref{sb:gwp:rho-equation} and use \eqref{sb:gwp:flow-identity}.
\end{proof}

\begin{lemma}\label{sb:gwp:no-collapse}
Suppose a nontrivial solution with the sign condition and
\eqref{sb:gwp:CZ} exists on \([0,T)\), \(T<\infty\). If
\begin{equation}\label{sb:gwp:origin-integrable}
 \int_0^T u_x(t,0)\,\mathrm{d}t<\infty,
\end{equation}
then \(\inf_{t<T,\,0\leq\theta\leq\pi}\varrho(t,\theta)>0\).
In particular \(\sup_{t<T}\|\omega(t)\|_\infty<\infty\).
\end{lemma}

\begin{proof}
For \(\omega_0(\theta)\ne0\), let
\(E_\theta=\varrho_t^2+\omega_0(\theta)^2/(4\varrho^2)\).
Equation \eqref{sb:gwp:rho-equation} gives
\(E_\theta'=a\mathscr F_u\varrho\varrho_t\).
Since \(\varrho\) is convex, its decreasing interval, if present,
starts at zero. On that interval \(E_\theta\leq E_\theta(0)\).
After its minimum \(\varrho\) increases. Consequently
\begin{equation}\label{sb:gwp:angular-barrier}
 \varrho(t,\theta)^2\geq
 \frac{\omega_0(\theta)^2}{u_{0x}(\theta)^2+\omega_0(\theta)^2}
 \quad\hbox{if }\omega_0(\theta)\ne0.
\end{equation}
At an interior zero, \eqref{sb:gwp:CZ} and convexity give
\(\varrho(t,\theta)\geq1\). At the right endpoint,
\[
 u_{0x}(\pi)=\frac1\pi\int_0^\pi\omega_0(y)\tan(y/2)\,\mathrm{d}y<0,
\]
so \(\varrho(t,\pi)\geq1\). At the origin,
\begin{equation}\label{sb:gwp:origin-barrier}
 u_x(t,0)=-\frac1\pi\int_0^\pi\omega(t,y)\cot(y/2)\,\mathrm{d}y\geq0,
 \qquad
 \varrho(t,0)=\exp\left(-\frac12\int_0^t u_x(s,0)\,\mathrm{d}s\right).
\end{equation}
Thus \eqref{sb:gwp:origin-integrable} bounds \(\varrho(t,0)\)
away from zero.

These pointwise barriers give a uniform barrier without assuming
isolated zeros. Otherwise choose \(t_n\uparrow T\), \(\theta_n\),
and minimizing times \(s_n\in[0,t_n]\) with
\(\varrho(s_n,\theta_n)\to0\).
Compact regular-time positivity implies \(s_n\to T\).
Pass to \(\theta_n\to\theta_*\). Convexity gives, for \(t\leq s_n\),
\[
 -C_0\leq\varrho_t(t,\theta_n)\leq0,\qquad
 \varrho(t,\theta_n)\leq\varrho(s_n,\theta_n)+C_0(s_n-t),
 \quad C_0=\tfrac12\|u_{0x}\|_\infty.
\]
For each fixed \(t<T\), continuity in the label yields
\(\varrho(t,\theta_*)\leq C_0(T-t)\).
Letting \(t\uparrow T\) contradicts the appropriate barrier above,
whether \(\theta_*\) is a nonzero label, an interior zero, or an
endpoint. Finally \eqref{sb:gwp:rho-data} bounds \(\omega\) uniformly.
\end{proof}

\subsection{Completion of the global argument}

\begin{proof}[Proof of Theorem~\ref{sb:gwp:theorem}]
The zero solution is immediate. Suppose the nontrivial solution has
a finite maximal time \(T_*\).
First assume \eqref{sb:gwp:trace}, and set \(m=m_a(\omega_0)>0\).
For any regular \(T_0<T_*\), \eqref{sb:gwp:trace-propagation} and
\(\cot(x/2)=(2/x)(1+O(x^2))\) give
\begin{equation}\label{sb:gwp:residue}
 \limsup_{\varepsilon\downarrow0}
   \varepsilon\mathcal E_{\beta_c-\varepsilon}(T_0)
       \leq 2^{\beta_c}m.
\end{equation}
Indeed, the singular part of the integral is
\(2^{\beta_c-\varepsilon}m\int_0^d x^{-1+\varepsilon}\mathrm{d}x\);
the factor from trace propagation differs from one by
\(O_{T_0}(x^{1/2})\), whose integral stays bounded as
\(\varepsilon\downarrow0\). The remaining interval and corrected
mass term make no contribution to the limit in
\eqref{sb:gwp:residue}.

Choose \(K_a>2^{\beta_c}\). For each \(T_0\) choose a sufficiently
small \(\varepsilon>0\) with \(\beta\geq\beta_*\) and
\(\varepsilon\mathcal E_\beta(T_0)\leq K_am\).
The scalar differential inequality \eqref{sb:gwp:corrected-inequality}
gives
\[
 \mathcal E_\beta(t)\leq
 \frac{\mathcal E_\beta(T_0)}
 {1-C_a\varepsilon\mathcal E_\beta(T_0)(t-T_0)}.
\]
Thus the moment remains bounded for a forward time at least
\((2C_aK_am)^{-1}\), independently of \(T_0\).
Take \(T_0\) close enough to \(T_*\) that this interval covers
\([T_0,T_*)\). For the resulting fixed positive \(\varepsilon\),
both \(Q_\beta\) and \(M\) stay bounded up to \(T_*\).

If instead \eqref{sb:gwp:weighted-endpoint} holds, start at zero.
Splitting where \(\cot(x/2)\) is above and below one gives
\(Q_{\beta_c-\varepsilon}(0)\leq Q_{\beta_c}(0)+M(0)\).
Choose \(\varepsilon>0\) so small that
\(C_a\varepsilon\mathcal E_\beta(0)T_*<1/2\).
The same comparison bounds \(Q_\beta\) and \(M\) on \([0,T_*)\).

In either case,
\(h\leq Q_\beta^{1/\beta}M^{1-1/\beta}\) and
\eqref{sb:gwp:origin-barrier} imply
\(\int_0^{T_*}u_x(t,0)\,\mathrm{d}t<\infty\).
Lemma~\ref{sb:gwp:no-collapse} bounds \(\omega\) in \(L^\infty\),
and Proposition~\ref{sb:gwp:continuation} extends the solution
beyond \(T_*\), a contradiction. Higher regularity follows from
\eqref{sb:gwp:higher-energy} on each finite interval.
\end{proof}

\subsection{Blowup from a finite inverse-vorticity moment}
\label{sb:inverse:subsection}

The flow identity \eqref{sb:gwp:flow-identity} and the positive velocity
formula \eqref{sb:gwp:positive-velocity} hold for every \(a>0\).
They also give a blowup criterion.
The inverse-vorticity measure has constant mass, whereas its first
spatial moment grows at a uniform positive rate.

\begin{theorem}[Inverse-vorticity criterion]\label{sb:inverse:blowup}
Let \(0<a<1\), and let \(\omega_0\in C^\infty(\mathbb T)\) be real and odd.
Suppose
\[
 \omega_0(x)<0\quad(0<x<\pi),\qquad
 M_{-a}:=\int_0^\pi(-\omega_0(x))^{-a}\,\idiff x<\infty.
\]
Then the maximal lifespan \(T_*\) is finite, and the solution satisfies
\[
 \limsup_{t\uparrow T_*}\|\omega(t)\|_{L^\infty(\mathbb T)}=\infty.
\]
In particular, the datum \(\omega_0(x)=-\sin x\) gives finite-time
blowup for every \(0<a<1\).
\end{theorem}

\begin{proof}
Let \(\eta\) and \(J=\eta_\theta\) be the flow and Jacobian in
\eqref{sb:gwp:flow-identity}. Strict negativity is preserved in
\((0,\pi)\), and that identity gives
\[
 (-\omega(t,\eta(t,\theta)))^{-a}J(t,\theta)
       =(-\omega_0(\theta))^{-a}.
\]
Changing variables therefore yields
\begin{equation}\label{sb:inverse:conservation}
 \int_0^\pi(-\omega(t,x))^{-a}\,\idiff x=M_{-a}.
\end{equation}
Define the first spatial moment by
\[
 F(t)=\int_0^\pi x(-\omega(t,x))^{-a}\,\idiff x
     =\int_0^\pi\eta(t,\theta)(-\omega_0(\theta))^{-a}\,\idiff\theta.
\]
The velocity is bounded on each compact regular time interval, so the
last integral can be differentiated by dominated convergence. Hence
\begin{equation}\label{sb:inverse:first-moment}
 0<F(t)<\pi M_{-a},\qquad
 F'(t)=a\int_0^\pi u(t,x)(-\omega(t,x))^{-a}\,\idiff x.
\end{equation}

To bound this derivative from below, fix \(I=[\pi/3,2\pi/3]\).
For \(x,y\in I\), \(x\ne y\), the identity
\[
 \sin^2((x+y)/2)-\sin^2((x-y)/2)=\sin x\sin y\geq3/4
\]
and the bound \(\sin^2((x-y)/2)\leq1/4\) give
\[
 \log\left|\frac{\sin((x+y)/2)}{\sin((x-y)/2)}\right|\geq\log2.
\]
The kernel in \eqref{sb:gwp:positive-velocity} is nonnegative on the
whole half-circle. Tonelli's theorem and
\eqref{sb:inverse:first-moment} consequently imply
\begin{equation}\label{sb:inverse:restricted-growth}
 F'(t)\geq\frac{a\log2}{\pi}
 \left(\int_I(-\omega)^{-a}\,\idiff x\right)
 \left(\int_I(-\omega)\,\idiff x\right).
\end{equation}
H\"older's inequality gives
\[
 |I|\leq
 \left(\int_I(-\omega)^{-a}\,\idiff x\right)^{1/(a+1)}
 \left(\int_I(-\omega)\,\idiff x\right)^{a/(a+1)}.
\]
Since \(1-1/a<0\), combining this inequality with
\eqref{sb:inverse:conservation}--\eqref{sb:inverse:restricted-growth}
yields
\begin{equation}\label{sb:inverse:uniform-growth}
 \begin{aligned}
 F'(t)
 &\geq\frac{a\log2}{\pi}|I|^{1+1/a}
       \left(\int_I(-\omega)^{-a}\,\idiff x\right)^{1-1/a}\\
 &\geq\frac{a\log2}{\pi}
       \left(\frac{\pi}{3}\right)^{1+1/a}M_{-a}^{1-1/a}>0.
 \end{aligned}
\end{equation}
Integration contradicts the upper bound in
\eqref{sb:inverse:first-moment} if the solution is global.
Thus \(T_*<\infty\), and Proposition~\ref{sb:gwp:continuation}
gives the asserted vorticity blowup. For \(\omega_0=-\sin x\), the
inverse moment is finite because both endpoint zeros are simple
and \(a<1\).
\end{proof}

\section{Focusing blowup by finite-dimensional selection}
\label{sb:focus:section}

We transfer a strictly focusing whole-line profile to smooth periodic
data. The proof first obtains the coefficient bounds at the origin and
infinity, then constructs the closed normalized generator and its finite
spectral splitting. Moving localization gives exponentially decaying
periodic errors. A lower decaying norm and a higher bounded energy close
the evolution estimates, after which a finite-dimensional exit argument
selects smooth initial data.

\subsection{Profile assumptions}

\begin{assumption}[Strictly focusing profile]\label{sb:focus:profile-input}
For a parameter $a\in(0,1)$ there are odd smooth functions $\bar\Omega,\bar U$ on
$\mathbb R$, constants $\bar c_l>0$, $\bar c_\omega<0$, and $\rho=1-\bar c_\omega/\bar c_l>2$
such that
\begin{equation}\label{sb:focus:profile}
 \bar U(0)=0,\quad \bar U'=H\bar\Omega,\quad
 (\bar c_lx+a\bar U)\bar\Omega'=(\bar c_\omega+\bar U')\bar\Omega,
 \qquad \bar\Omega'(0)=-1,\quad \bar\Omega'''(0)=6.
\end{equation}
Moreover $P=\bar U/x>0$ on $(0,\infty)$, $\bar\Omega(x)\sim-Ax^{1-\rho}$ for some
$A>0$, and every ordinary derivative of $\bar\Omega$ has at most polynomial
growth at infinity.
\end{assumption}

For every selected whole-line fixed point with positive spatial rate,
the reconstruction and smoothness statements follow from
\cite[Theorems~4.4 and~5.8]{Profiles}, and the tail and derivative bounds
from \cite[Proposition~5.5 and Corollary~5.12]{Profiles}.
Thus Assumption~\ref{sb:focus:profile-input} holds. For
\(0<a\leq a_{\mathrm{foc}}\), such profiles exist by
\cite[Theorem~1.3(i)]{Profiles}.

\begin{lemma}[Conormal profile bounds]\label{sb:focus:conormal}
Let $D=x\partial_x$. Under Assumption~\ref{sb:focus:profile-input}, for
every integer $j\geq0$,
\begin{equation}\label{sb:focus:tails}
 |D^j\bar\Omega(x)|\leq C_jx^{1-\rho},\qquad
 |D^jP(x)|+|D^jH\bar\Omega(x)|\leq C_jx^{-1/2}\quad(x\geq1).
\end{equation}
In particular no finite-second-moment assumption is needed.
\end{lemma}
\begin{proof}
For odd $q$ with $|q(x)|\lesssim x^{1-\rho}$ and
$|q'(x)|\lesssim x^M$, pair the two half-lines in the Hilbert integral.
The regions $y<x/2$ and $y>3x/2$ are bounded by
\[
 Cx^{-2}\int_0^{x/2}y(1+y)^{1-\rho}\,\mathrm{d}y,
 \qquad C\int_{3x/2}^\infty y^{-\rho}\,\mathrm{d}y.
\]
On the intervening interval subtract $q(x)$ in the singular kernel.
Its integral is bounded by
$C\int_0^{x/2}\min\{x^{1-\rho},x^Mr\}\,\mathrm{d}r/r$.
Thus $Hq=O(x^{-\min(\rho-1,2)}\log(2+x))$.
The same three-region split of the positive logarithmic formula for
$P$ gives this bound for $P$. That formula follows by integrating $H\bar\Omega$;
its density is $-\bar\Omega(y)>0$, since the stationary scalar ODE and $\bar\Omega'(0)<0$
imply $\bar\Omega<0$ on $(0,\infty)$.
Now use
\[
 D\bar\Omega=\frac{\bar c_\omega+H\bar\Omega}{\bar c_l+aP}\bar\Omega,\qquad DP=H\bar\Omega-P,
 \qquad [D,H]=0.
\]
The denominator is at least $\bar c_l$. Inductively the Hilbert estimate
controls $H(D^j\bar\Omega)$, and these identities give the next derivatives of
$P$ and $\bar\Omega$. The assumed ordinary derivative bounds supply the
polynomial Lipschitz bound at each step. Dilation covariance proves the
commutator distributionally; the conormal derivatives already reached
are integrable because $\rho>2$. This proves the induction.
\end{proof}

The weights must accommodate the two canceled origin derivatives and the
algebraic tail. Smooth tangent perturbations are $O(x^5)$ at the origin,
so a left exponent between four and five both admits these perturbations
and controls four ordinary derivatives. At infinity we choose a small
positive exponent below the vorticity decay exponent $\rho-1$.
Fix an integer $m\geq7$ and put $n=m+2$. Define, for $x=e^z$,
\begin{equation}\label{sb:focus:spaces}
 \nu=\frac92,\quad \sigma=\frac1{16},\quad
 w(z)=e^{-\nu z}+e^{\sigma z},\qquad
 \|h\|_{X_k}^2=\sum_{j=0}^k\int_{\mathbb R}
      w(z)^2|\partial_z^jh(e^z)|^2\,\mathrm{d}z.
\end{equation}
The space $X_k$ is the completion of odd compactly supported smooth
functions with $h'(0)=h'''(0)=0$. Equivalently it is the weighted
Sobolev space on the positive logarithmic half-line, with odd extension.
Completed elements need not be pointwise $O(x^5)$.

\begin{theorem}[Smooth focusing blowup]\label{sb:focus:theorem}
Assume~\ref{sb:focus:profile-input}. There are finitely many odd
$\psi_j\in C_c^\infty(\mathbb R)$ with vanishing first and third derivatives at the origin, and
constants $\rho_{\mathrm{loc}},\varepsilon>0$, with the following property.
For every odd compactly supported smooth seed $h$ with zero first and
third jets and $\|h\|_{X_n}\leq\varepsilon$, and all sufficiently large
$L_0$ depending on its support, one can choose finitely many real $z_j$
so that
\begin{equation}\label{sb:focus:data}
 \omega_0(y)=\left(\chi(x/R_0)\bar\Omega(x)+h(x)+\sum_jz_j\psi_j(x)\right)_
 {x=L_0y/\pi},\qquad R_0=\rho_{\mathrm{loc}} L_0,
\end{equation}
is smooth and odd on $\mathbb T=\mathbb R/(2\pi\mathbb Z)$ and its
solution blows up at a finite time $T$. Here $\chi$ is a fixed even
smooth cutoff, equal to one on $[-1,1]$, positive on $(-2,2)$ and zero
outside that interval. The datum is supported in a proper subarc.

There are $\kappa>0$ and normalized scales as in
\eqref{sb:common:rescaling}, with $C_\omega(0)=1$ and $C_l(0)=L_0/\pi$,
such that, for the zero extension of the central rescaled cell,
\begin{equation}\label{sb:focus:normalized-limit}
 \|\Omega_{\rm ext}(\tau)-\bar\Omega\|_{X_m}
 +|c_l(\tau)-\bar c_l|+|c_\omega(\tau)-\bar c_\omega|
 \leq Ce^{-\kappa\tau}.
\end{equation}
Writing $\lambda=C_l^{-1}$, one has for some $\lambda_*>0$
\begin{equation}\label{sb:focus:physical-limit}
 \lambda(t)=\lambda_*(T-t)^{-\bar c_l/\bar c_\omega}
       \bigl(1+O((T-t)^{-\kappa/\bar c_\omega})\bigr),\qquad
 (T-t)\|\omega(t)\|_\infty\longrightarrow-\frac{\|\bar\Omega\|_\infty}{\bar c_\omega}.
\end{equation}
The central-cell lift of $(T-t)\omega(\lambda(t)x,t)$ converges
uniformly to $-\bar\Omega/\bar c_\omega$, and the actual composition converges in
$C^4_{\rm loc}(\mathbb R)$. Only the coefficients $z_j$ are tuned.
\end{theorem}

\subsection{The closed normalized operator}

Use the gauges of Lemma~\ref{sb:common:gauges}. For a line perturbation
$h$ put
\[
 \begin{gathered}
 U_h(x)=\int_0^xHh(y)\,\mathrm{d}y,\\
 \ell_l(h)=-\frac{3-a}{12}(Hh)''(0)-aHh(0),\qquad
 \ell_\omega(h)=-\frac{3-a}{12}(Hh)''(0)-Hh(0).
 \end{gathered}
\]
Set $\bar B (x)=\bar c_lx+a\bar U(x)$ and $q(x)=\bar c_\omega+H\bar\Omega(x)$.
The exact linear and quadratic expressions are
\begin{align}
 \mathcal Lh&=-\bar B h'+qh+\mathcal Rh,
 \notag\\
 \mathcal Rh&=\bar\Omega Hh-a\bar\Omega'U_h+\ell_\omega(h)\bar\Omega-\ell_l(h)D\bar\Omega,
 \label{sb:focus:linear}\\
 \mathcal N(h)&=(Hh+\ell_\omega(h))h
                  -(aU_h/x+\ell_l(h))Dh.
 \label{sb:focus:quadratic}
\end{align}
The linear nonlocal terms are always grouped as $\mathcal R$.
Their separate summands need not belong to $X_k$.

\begin{lemma}[Weighted calculus]\label{sb:focus:calculus}
For every finite $k$, $\mathcal R$ is bounded on $X_k$.
Its real quadratic form, localized to a bounded logarithmic interval,
has a compact self-adjoint representative on $X_k$, and its source and
target tails tend to zero in operator norm. In addition
\begin{align}
 |Hh(0)|+|(Hh)''(0)|+\|Hh\|_\infty+\|U_h/x\|_\infty
 &\leq C\|h\|_{X_2},\label{sb:focus:traces}\\
 \|\mathcal N(h)\|_{X_k}&\leq C\|h\|_{X_{k+2}}\|h\|_{X_k},
 \quad k\geq2.\label{sb:focus:tame}
\end{align}
Weighted Sobolev embedding controls the ordinary $C^4$ norm by $X_m$.
\end{lemma}

\begin{proof}
Let $F(z)=h(e^z)$, $O(z)=\bar\Omega(e^z)$ and $O_1=\partial_zO$.
The Hilbert and velocity-average kernels are respectively
\begin{equation}\label{sb:focus:kernels}
 K_H(t)=\frac2{\pi(e^{2t}-1)},\qquad
 K_U(t)=\frac{e^{-t}}\pi\log|\tanh(t/2)|,
 \qquad (\partial_z+1)\mathscr U=\mathscr H.
\end{equation}
They behave at zero as $(\pi t)^{-1}+O(1)$ and $O(1+|\log|t||)$.
At the positive end both are $O(e^{-2t})$. At the negative end, with
$c_0=-2/\pi$ and $c_2=-4/\pi$,
\[
 K_H=c_0+\tfrac12c_2e^{2t}+O(e^{4t}),\qquad
 K_U=c_0+\tfrac16c_2e^{2t}+O(e^{4t}).
\]
Also $Hh(0)=c_0\int F$ and $(Hh)''(0)=c_2\int e^{-2z}F(z)\,\mathrm{d}z$.
Both moments are bounded by Cauchy--Schwarz with $w^{-1}$.
The exact grouped kernel of $\mathcal R$ is
\begin{equation}\label{sb:focus:grouped-kernel}
 M(z,\eta)=O(z)K_H(z-\eta)-aO_1(z)K_U(z-\eta)
 +c_0(aO_1-O)(z)-\frac{3-a}{12}c_2(O-O_1)(z)e^{-2\eta}.
\end{equation}
In the cone $z\ll-1$, $\eta\geq z+1$, the constant terms cancel and
\begin{equation}\label{sb:focus:quadratic-cancel}
 -\frac{3-a}{12}(O-O_1)+e^{2z}(O/2-aO_1/6)=O(e^{5z}).
\end{equation}
This follows from $O=-e^z+e^{3z}+O(e^{5z})$ and fixes both signs of
the quadratic subtraction.

For clarity we give envelopes for the off-diagonal kernel after
conjugation by $w$. The same envelopes hold after any fixed number
of derivatives in either variable, once input derivatives are moved
off convolution kernels. At the left end they are
\[
 \begin{cases}
 Ce^z e^{-(\nu-2)(z-\eta)},&z\leq-2,\ \eta\leq z-1,\\
 Ce^{(5-\nu)z}(1+e^{-4\eta})/w(\eta),
                  &z\leq-2,\ \eta\geq z+1.
 \end{cases}
\]
On $z,\eta\geq1$, the convolution pieces have envelopes
\[
 \begin{cases}
 Ce^{-(\rho-1)z}e^{-(2-\sigma)(z-\eta)},&z\geq\eta+1,\\
 Ce^{-(\rho-1)z}e^{-\sigma(\eta-z)},&\eta\geq z+1.
 \end{cases}
\]
The moment pieces have target factor $e^{-(\rho-1-\sigma)z}$ and source
factors $w(\eta)^{-1}$ or $e^{-2\eta}w(\eta)^{-1}$.
When the target is near the right end and the source near the left,
their slower factors are
\[
 e^{-(\rho-1-\sigma)z}e^{\nu\eta},\qquad
 e^{-(\rho-1-\sigma)z}e^{(\nu-2)\eta}.
\]
In the reverse region \eqref{sb:focus:quadratic-cancel} gives
\[
 e^{(5-\nu)z}(e^{-(2+\sigma)\eta}+e^{-(4+\sigma)\eta}).
\]
Every row and column tail is integrable and exponentially small.
One may implement these subtractions with a smooth cutoff equal to
one in the indicated left forward cone. Derivatives of that cutoff
lie in a fixed strip or where $1/2\leq\eta-z\leq1$ and obey the same
adjacent envelopes. Thus no differentiated cutoff terms are omitted.

Near the diagonal the sole principal block is
$O(z)/(\pi(z-\eta))$; all other local kernels are logarithmic or
bounded. The coefficient $O$ and all its derivatives vanish at both
ends, giving small source and target tails for this block as well.
On a compact square its symmetrized principal numerator is
\[
 \frac{w(z)O(z)}{w(\eta)}-\frac{w(\eta)O(\eta)}{w(z)}.
\]
It vanishes on the diagonal, so the symmetric kernel is bounded.
It is therefore Hilbert--Schmidt; the Volterra logarithm is also
Hilbert--Schmidt there. The nonsymmetric Hilbert block is only claimed
to be bounded. After $k$ derivatives the sole equal-order terms are
these principal blocks on $\partial_z^kF$; every other localized
term factors through $H^{k-1}$ and local Rellich compactness. This
proves the asserted compact form. The same statements hold for
$\vartheta w$ when $\vartheta$ is smooth, bounded above and below,
and constant at both ends.

The unweighted Hilbert expression in \eqref{sb:focus:kernels} is the
sum of a local Hilbert kernel, an integrable kernel and a one-sided
constant tail. The first two are controlled by Sobolev embedding and
the last by $\|F\|_1\leq\|wF\|_2\|w^{-1}\|_2$.
This proves \eqref{sb:focus:traces}; Leibniz' rule gives
\eqref{sb:focus:tame}. Finally
$\partial_x^jh=x^{-j}\prod_{r=0}^{j-1}(\partial_z-r)F$.
At zero the weighted bound gains $x^{\nu-j}$, positive for $j\leq4$,
and at infinity it gains $x^{-\sigma-j}$. This proves the embedding.
\end{proof}

\begin{proposition}[Finite spectral splitting]\label{sb:focus:spectrum}
On every $X_k$, $2\leq k\leq n$, the expression \eqref{sb:focus:linear} has
the closed realization
\[
 D(\mathcal L_k)=X_{k+1},\qquad \mathcal L_k=-\bar B \partial_x+q+\mathcal R.
\]
Its semigroup has negative essential growth bound. In particular, for
$k=m$ there are bounded real projections $\Pi_u,\Pi_s$, with
$N=\operatorname{rank}\Pi_u<\infty$, and numbers
\begin{equation}\label{sb:focus:rates}
 0\leq\eta_u<\eta<\kappa<\min\{\eta_s,\mu\},\qquad
 \mu=\frac7{32}\bar c_l,
\end{equation}
such that
\[
 \|e^{\tau\mathcal L_m}\Pi_s\|\leq Me^{-\eta_s\tau},\qquad
 \|e^{-\tau\mathcal L_m|_{X_u}}\|\leq Me^{\eta_u\tau}.
\]
There is a Euclidean norm $|\cdot|_u$ on $X_u=\Pi_uX_m$ with
$\operatorname{Re}\langle\mathcal L_mz,z\rangle_u\geq-\eta_u|z|_u^2$.
\end{proposition}

\begin{proof}
In logarithmic coordinates put $B(z)=\bar B (e^z)/e^z$, $Q(z)=q(e^z)$.
All their derivatives are bounded, $B\geq\bar c_l>0$, and
\[
 (B,Q)(-\infty)=(B_-,B_-),\quad B_-=\bar c_l+aH\bar\Omega(0)>0,
 \qquad (B,Q)(+\infty)=(\bar c_l,\bar c_\omega).
\]
The transport flow of $B$ is global. Its characteristic formula and
the bounded logarithmic derivatives of $w$ generate a strongly
continuous semigroup on each $X_k$. The maximal distributional domain
is $X_{k+1}$, because $F_z=B^{-1}(QF-TF)$.
Logarithmic cutoffs followed by mollification give a graph core of
smooth functions supported away from both endpoints. Bounded
perturbation by $\mathcal R$ proves generation for $\mathcal L_k$.
The realizations and semigroups agree on their common regularity classes.

For $0\leq j\leq n$ the top transport energy coefficient is
\[
 d_j=Q-(j-\tfrac12)B'+B(\log w)',
 \qquad d_j(-\infty)=-(\nu-1)B_-,\quad
 d_j(+\infty)=\bar c_\omega+\sigma\bar c_l<0.
\]
Choose a nonnegative compactly supported $\zeta$ and $M_0<\infty$
so that $\vartheta(z)=\exp\bigl(-M_0\int_{-\infty}^z\zeta(s)\,\idiff s\bigr)$ makes every
coefficient with weight $\vartheta w$ strictly negative. Local lower
commutators and Lemma~\ref{sb:focus:calculus} then give
\begin{equation}\label{sb:focus:compact-energy}
 2\operatorname{Re}\langle\mathcal L_kF,F\rangle_{X_{k,\vartheta}}
 \leq-8\delta\|F\|_{X_{k,\vartheta}}^2
       +C\|F\|_{H^{k-1}([-R,R])}^2
\end{equation}
for some $\delta>0$. The restriction $J:X_k\to H^{k-1}([-R,R])$
is compact, so $J^*J$ is a compact self-adjoint operator on $X_k$.
The generator $\mathcal L_k-(C/2)J^*J$ still generates a semigroup
by bounded perturbation. Applying \eqref{sb:focus:compact-energy}
on the graph core and then using density shows that this semigroup
decays at least as $e^{-4\delta\tau}$. This step uses the compact
quadratic remainder supplied by the symmetric form of $\mathcal R$.
A bounded compact perturbation changes a semigroup by compact operators
at positive times, so its essential growth bound is unchanged
\cite[Proposition~IV.2.12]{EngelNagel2000}.

Choose a line $\operatorname{Re}\lambda=-\eta$ with
$0<\eta<\min(2\delta,\mu/4)$ which misses the spectrum.
The spectrum to its right consists of finitely many poles of finite
algebraic multiplicity; the semigroup growth bound is the maximum of
the essential and spectral bounds \cite[Corollary~IV.2.11]{EngelNagel2000}.
The Riesz decomposition \cite[Proposition~IV.1.16]{EngelNagel2000} gives
the stated stable estimate and finite-dimensional backward estimate,
with strict gaps on either side of the chosen line. It permits
\eqref{sb:focus:rates}. A Lyapunov inner product for
$-(\mathcal L_m|_{X_u}+\eta_u)$ gives the final inequality.
This construction includes zero eigenvalues and neutral Jordan blocks.
\end{proof}

\subsection{Moving localization and the two energy levels}

Let $R(\tau)=R_0e^{\bar c_l\tau}$,
$L_{\rm app}(\tau)=L_0e^{\bar c_l\tau}$ and
$\Phi(x,\tau)=\chi(x/R(\tau))\bar\Omega(x)$. The actual period is
$2L(\tau)=2\pi C_l(\tau)$. While the physical support remains in the
central arc, let $S_\Omega$ be its rescaled outer endpoint and set
\[
 h=\Omega_{\rm ext}-\Phi,\qquad
 S_h=\max\{S_\Omega,2R\},\qquad G_0=(1+\rho_{\mathrm{loc}}^2)R_0^{-7/32}.
\]
Thus $h$ is supported in $[-S_h,S_h]$, and its first and third derivatives at the origin
vanish. Only the physical endpoint $S_\Omega$ is differentiated;
$S_h$ records the support of the actual perturbation input.

\begin{lemma}[Periodic forcing and high energy]\label{sb:focus:nonlinear}
Suppose $S_h/R$ is bounded above and below, $2S_h\leq L$,
$1/2\leq L/L_{\rm app}\leq2$, and $\|h\|_{X_n}\leq1$.
The exact perturbation equation is
$h_\tau=\mathcal Lh+\mathcal F(h,\tau)$, with
\begin{equation}\label{sb:focus:forcing}
 \|\mathcal F\|_{X_m}
 \leq C\{\|h\|_{X_n}\|h\|_{X_m}
               +(S_h/L)^2\|h\|_{X_m}+G_0e^{-\mu\tau}\}.
\end{equation}
There is a fixed equivalent high energy $E_n\asymp\|h\|_{X_n}^2$
satisfying
\begin{equation}\label{sb:focus:high-energy}
 E_n'\leq-\bigl(4\delta-C(\|h\|_{X_n}+(S_h/L)^2+G_0e^{-\mu\tau})\bigr)E_n
       +C\|h\|_{X_m}^2+CG_0^2e^{-2\mu\tau}.
\end{equation}
Moreover
\begin{equation}\label{sb:focus:rate-error}
 |c_l-\bar c_l|+|c_\omega-\bar c_\omega|
 \leq C\{(1+(S_h/L)^2)\|h\|_{X_m}+G_0e^{-\mu\tau}\}.
\end{equation}
The constants are uniform in the period and support sizes under the
stated ratios.
\end{lemma}

\begin{proof}
The analytic kernel difference is
\[
 E_L(s)=\frac1{2L}\cot\frac{\pi s}{2L}-\frac1{\pi s}.
\]
For $|s|\leq L$, its zeroth and first derivatives are bounded by
$C|s|/L^2$ and $C/L^2$; its $j$th derivative, $j\geq2$, is bounded
by $C_jL^{-j-1}$. These bounds persist after any fixed number of
$\partial_{\log L}$ derivatives. Equivalently
$SE_L(S\xi)=(S/L)^2e_{S/L}(\xi)$, where the smooth functions
$e_{\rho_{\mathrm{loc}}}$ have uniform derivative bounds on $|\xi|\leq2$, $\rho_{\mathrm{loc}}\leq1/2$.
The support premise ensures that every kernel use stays in this region.
For an odd supported input and $r\geq0$,
\begin{equation}\label{sb:focus:supported-moment}
 \int_0^S y^r|D^jh(y)|\,\mathrm{d}y
 \leq C_r(1+S^{r+1-\sigma})\|h\|_{X_j}.
\end{equation}
Write $\Delta H=H_L-H$, $\Delta A h=x^{-1}\int_0^x\Delta Hh$,
and use $\Delta\ell$ for their rate differences. At the origin
\[
 \Delta Hh=a_0+\tfrac12a_2x^2+O(x^4),\qquad
 \Delta Ah=a_0+\tfrac16a_2x^2+O(x^4).
\]
Consequently
$\Phi\Delta Hh-aD\Phi\Delta Ah+\Delta\ell_\omega(h)\Phi
-\Delta\ell_l(h)D\Phi$ has the same two cancellations as
\eqref{sb:focus:quadratic-cancel}. Its $X_k$ norm is bounded by
$C(S_h/L)^2\|h\|_{X_k}$ through $k=n$.
The corresponding quadratic expression has bound
$C(S_h/L)^2\|h\|_{X_{k+2}}\|h\|_{X_k}$ for $k\leq m$.
These follow by the scaled kernel bounds and
\eqref{sb:focus:supported-moment}; the origin estimate is for the
complete grouped expression before taking absolute values.

Let $\mathcal E^L$ be minus the residual obtained by substituting
$\Phi$ into the normalized equation with period $2L$ and its own
gauges. Since $R'/R=\bar c_l$, the largest cutoff transport cancels
exactly. The remaining line terms are
\[
 -a\bar U\partial_x\chi(x/R)\bar\Omega
 -aU_{\Phi-\bar\Omega}\Phi'+H(\Phi-\bar\Omega)\Phi,
\]
together with their rate corrections. Lemma~\ref{sb:focus:conormal}
and the grouped calculus give
\[
 \|\Phi-\bar\Omega\|_{X_{n+1}}+\|\mathcal E^{L_{\rm app}}\|_{X_n}
 \leq CG_0e^{-\mu\tau}.
\]
Indeed the line powers are $R^{-(\rho-1-\sigma)}$ and
$R^{-(1/2-\sigma)}$. The periodic profile moments satisfy
\[
 \int_0^{2R}y^r|D^j\Phi|\lesssim
 R^{\max(r+2-\rho,0)}(1+\log R).
\]
Their worst weighted error is
$C\rho_{\mathrm{loc}}^2R^{-(1-\sigma)}(1+\log R)$.
All these powers dominate $R^{-7/32}$.
The same estimate is Lipschitz in $\log L$, hence
\[
 \|\mathcal E^L-\mathcal E^{L_{\rm app}}\|_{X_n}
 \leq CG_0e^{-\mu\tau}|\log(L/L_{\rm app})|.
\]

For an explicit accounting, if $\mathcal B$ polarizes $\mathcal N$,
the forcing consists of $\mathcal N(h)+2\mathcal B(\Phi-\bar\Omega,h)$,
the two periodic linear polarizations, the periodic quadratic term,
and $\mathcal E^L$. The two linear polarizations are
\[
 \Phi\Delta Hh-aD\Phi\Delta Ah
       +\Delta\ell_\omega(h)\Phi-\Delta\ell_l(h)D\Phi,
 \quad
 (\Delta H\Phi+\Delta\ell_\omega(\Phi))h
       -(a\Delta A\Phi+\Delta\ell_l(\Phi))Dh.
\]
The second must be retained: it is the half on which the transform
acts on the profile. Its coefficients and all required derivatives
are bounded by $CG_0e^{-\mu\tau}$. The preceding estimates and
\eqref{sb:focus:tame} prove \eqref{sb:focus:forcing}.

We detail the top energy to specify the derivative count. Put
$F(z)=h(e^z)$ and take
\[
 E_n=\sum_{j=0}^n\|\vartheta w\partial_z^jF\|_2^2.
\]
The linear part obeys \eqref{sb:focus:compact-energy}. Its compact
remainder is absorbed by
\[
 \|F\|_{H^{n-1}([-R,R])}^2
 \leq\epsilon E_n+C_\epsilon\|h\|_{X_m}^2.
\]
The only unknown $(n+1)$st derivative is self-transport. For its
coefficient $A_h=a\mathscr UF+\ell_l(h)$ and $G=\partial_z^nF$,
\[
 -2\operatorname{Re}\int(\vartheta w)^2G A_h\overline{G'}
 =\int(\vartheta w)^2\bigl(A_h'+2A_h(\log(\vartheta w))'\bigr)|G|^2.
\]
For the last Leibniz splits, the one-sided constant tail of $K_H$
is differentiated on the input. Integration by parts gives, for $r\geq1$,
\[
 \|\mathscr H\partial_z^rF\|_2
 \leq C(\|\partial_z^rF\|_2+\|\partial_z^{r-1}F\|_2),\qquad
 \|\partial_z^n\mathscr UF\|_2
 \leq C(\|\partial_z^{n-1}F\|_2+\|\partial_z^{n-2}F\|_2).
\]
Put the other weighted factor with at most two derivatives in
$L^\infty$. This proves the transport bound $C\|h\|_{X_n}E_n$.
For the top Hilbert product, symmetrize as in the linear calculation.
Its local singular numerator cancels. The undifferentiated one-sided
constant kernel is instead Hilbert--Schmidt: with $g=\vartheta wF$
and $b_*=1/(\vartheta w)$,
\[
 \iint_{z<\eta-1}|g(z)b_*(\eta)|^2\,\mathrm{d}z\,\mathrm{d}\eta
 \leq\|g\|_2^2\|b_*\|_2^2.
\]
It is not estimated by a Schur bound. All remaining Hilbert kernels
are integrable or local principal kernels. Thus
$|\operatorname{Re}\langle h,\mathcal N(h)\rangle_{X_{n,\vartheta}}|
\leq C\|h\|_{X_n}E_n$.

For periodic terms the scaled analytic kernel gives coefficient bounds
on $(-\infty,\log S_h]$: $L^\infty$ through order $n-2$, and unweighted $L^2$
at orders $n-1,n$, all with factor $(S_h/L)^2\|h\|_{X_n}$.
These differentiated coefficients decay at the left end; their
undifferentiated constants are only bounded, not put in $L^2$.
The same Leibniz allocation and self-transport integration yield
$C(S_h/L)^2(1+\|h\|_{X_n})E_n$.
The profile polarization uses its available $X_{n+1}$ bound, giving
$CG_0e^{-\mu\tau}E_n$, with no extra derivative of $h$.
Young's inequality for the residual proves \eqref{sb:focus:high-energy}.
All integrations first use logarithmic cutoffs; their boundary terms
vanish by weighted tail convergence. Finally the two exact gauge
functionals give \eqref{sb:focus:rate-error}.
\end{proof}

\subsection{Support, smooth continuation and selection}

\begin{proof}[Proof of Theorem~\ref{sb:focus:theorem}]
We give the choice of constants because the high norms of the smooth
tuning directions need not be small. Choose a large $K$ so that the
geometric forcing divided by $KD$ is smaller than
$(\kappa-\eta_u)/8$. Density of the smooth annular core supplies a
linear map $E:X_u\to C_c^\infty(\mathbb R)$ with vanishing normalized origin derivatives,
$\Pi_uE=I$, and $\|\Pi_sE\|_{X_u\to X_m}$ as small as desired.
Its image has a smooth basis $\psi_1,\ldots,\psi_N$.
Choose this leakage smaller than a fixed multiple of $1/K$.
Next choose a high radius $A$ large enough to dominate both
$K\|E\|_{X_u\to X_n}$ and the compact forcing in
\eqref{sb:focus:high-energy}. Only after fixing $K,E,A$ choose
$D>0$ and $\rho_{\mathrm{loc}}>0$ sufficiently small. The seed threshold
$\varepsilon$ is a fixed sufficiently small multiple of $D$.

For $v\in X_u$, $|v|_u\leq KD$, take
\[
 h(0)=h+E(v-\Pi_uh),\qquad
 \Omega(0)=\chi(x/R_0)\bar\Omega+h(0).
\]
Choose $L_0$ so large that $G_0\leq D$ and the supports of the seed
and $E$ lie in $|x|<R_0/2$. Then $\Pi_uh(0)=v$,
$\|\Pi_sh(0)\|_{X_m}\leq D$ and
$\|h(0)\|_{X_n}\leq c_HAD$, where $c_H$ is a fixed small energy
equivalence constant. The outer support endpoint is exactly $2R_0$
for every $v$: the profile has strict sign and the cutoff remains
positive arbitrarily close to that endpoint. Thus the data are
smooth periodic functions on the physical circle.

Consider the bootstrap conditions
\begin{align*}
 \|\Pi_sh(\tau)\|_{X_m}&\leq2CD e^{-\kappa\tau},&
 |\Pi_uh(\tau)|_u&\leq KD e^{-\kappa\tau},&
 \|h(\tau)\|_{X_n}&\leq2AD,\\
 (2C_S)^{-1}R\leq S_\Omega&\leq2C_SR,&
 S_\Omega/L&\leq8\rho_{\mathrm{loc}},&1/2\leq L/L_{\rm app}&\leq2.
\end{align*}
Here $C,C_S$ are fixed large constants. Take $\rho_{\mathrm{loc}}\leq1/16$;
then $S_h/L\leq8\rho_{\mathrm{loc}}$ and $2S_h\leq L$ already in this bootstrap,
so Lemma~\ref{sb:focus:nonlinear} is applicable before improving it.
The low norm obeys $\|h\|_{X_m}\leq C(K)D e^{-\kappa\tau}$,
independently of $A$. Stable Duhamel and \eqref{sb:focus:forcing} give
\[
 \|\Pi_sh(\tau)\|_{X_m}
 \leq Me^{-\eta_s\tau}D+C\int_0^\tau e^{-\eta_s(\tau-s)}
 \bigl((AD+\rho_{\mathrm{loc}}^2)C(K)De^{-\kappa s}+G_0e^{-\mu s}\bigr)\,\mathrm{d}s.
\]
Since $\kappa<\eta_s,\mu$, enlarging $C$ and reducing $D,\rho_{\mathrm{loc}}$
improves this to $CD e^{-\kappa\tau}$. The high energy then gives
\[
 E_n'\leq-2\delta E_n+C(1+C(K)^2)D^2e^{-2\kappa\tau}.
\]
The earlier choice of $A$, the smaller initial high radius $c_HAD$,
and energy equivalence improve the high norm to $AD$ with a strict
margin. This closes the two derivative levels without an $X_{n+2}$
norm. Equation~\eqref{sb:focus:rate-error} improves both rates.

The support is transported by the smooth physical characteristic flow.
Its central rescaled endpoint satisfies
\[
 S_\Omega'=c_lS_\Omega+aU_\Omega(S_\Omega),\qquad L'=c_lL.
\]
At $S_\Omega\asymp R$, the profile velocity, cutoff tail,
perturbation and periodic moment estimates give
\[
 |aU_\Omega(S_\Omega)/S_\Omega|
 \leq C\bigl(G_0e^{-\mu\tau}+\|h\|_{X_m}\bigr).
\]
The right-hand sides of the equations for
$\log(S_\Omega/R)$ and $\log(S_\Omega/L)$ are integrable and small.
Similarly $\log(L/L_{\rm app})=\int_0^\tau(c_l-\bar c_l)$.
Starting from $S_\Omega(0)=2R_0$, these bounds improve to
$C_S^{-1}R\leq S_\Omega\leq C_SR$, $S_\Omega/L\leq4\rho_{\mathrm{loc}}$ and
$3/4<L/L_{\rm app}<4/3$. In particular the support never reaches
the physical seam, and every use of the central-cell lift is justified.

We record the continuity needed for this bootstrap. The equation is
always solved on the fixed physical circle. On finite smooth intervals
the characteristic flow and solutions depend continuously on the
finite parameter $v$ in arbitrarily high ordinary Sobolev norms.
The normalization scales are recovered from the two nonzero normalized origin
derivatives and hence vary continuously there. Their normalized differences
have vanishing first and third derivatives at the origin. Taylor's theorem, applied also after
the finitely many $D$ derivatives, controls the origin contribution
to $X_n$ by an integrable multiple of $x^{5-\nu}$.
Away from zero, the support and period lie in compact ranges on such
an interval, so ordinary Sobolev continuity implies weighted
continuity. The common initial outer endpoint $2R_0$ is carried by
the continuous characteristic flow, making the support bounds
continuous as well. Finally the weighted $C^1$ bound, the bounded
scales on finite $\tau$ intervals and the fixed-circle continuation
criterion exclude termination of a smooth solution before a bootstrap
boundary. Higher ordinary derivatives propagate by the usual
differentiated transport energy.

All boundaries except the finite-dimensional one now improve strictly.
At a contact $|\Pi_uh|_u=KD e^{-\kappa\tau}$,
\[
 \frac{\mathrm{d}}{\mathrm{d}\tau}\bigl(e^{2\kappa\tau}|\Pi_uh|_u^2\bigr)
 \geq2e^{2\kappa\tau}|\Pi_uh|_u^2
 \left(\kappa-\eta_u-C\frac{\|\mathcal F\|_{X_m}}{|\Pi_uh|_u}\right)>0.
\]
Indeed the last ratio is at most $C(K^{-1}+AD+C(K)\rho_{\mathrm{loc}}^2)$.
The constant order chosen above makes the bracket strictly positive.
If every $v$ in the closed ball exited, its first exit time and
normalized exit direction would be continuous. For the time, this
follows from strict interior margins before contact and strict
outward motion just after contact; it also holds at initial boundary
points. The exit direction would retract the ball onto its boundary,
and equals the identity there, a contradiction. If $N=0$, the strict
improvements directly give a global normalized solution. In either
case at least one initial parameter remains trapped for all $\tau$.

The constructed solution has
$\|h\|_{X_m}\lesssim e^{-\kappa\tau}$ and
$\|\Phi-\bar\Omega\|_{X_m}\lesssim e^{-\mu\tau}$, proving
\eqref{sb:focus:normalized-limit}. The ordinary local convergence
follows from Lemma~\ref{sb:focus:calculus}. The rate errors are
integrable with exponential tails, and the common reconstruction
of Lemma~\ref{sb:common:physical} gives
\[
 C_l=C_{l,*}e^{\bar c_l\tau}(1+O(e^{-\kappa\tau})),\qquad
 C_\omega=C_{\omega,*}e^{-\bar c_\omega\tau}(1+O(e^{-\kappa\tau})),
 \quad (T-t)C_\omega=(-\bar c_\omega)^{-1}(1+O(e^{-\kappa\tau})).
\]
Here $T=\int_0^\infty C_\omega^{-1}\,\mathrm{d}\tau$ is finite.
The global uniform limit of the zero extension gives the supremum
asymptotic and rules out smooth continuation at $T$. Eliminating
$\tau$ proves \eqref{sb:focus:physical-limit}. The whole-line uniform limit concerns the zero extension of the central
cell, while the periodic composition has the stated local convergence.
\end{proof}

\begin{corollary}[An explicit focusing parameter]\label{sb:focus:explicit}
There are smooth odd periodic data with focusing blowup at $a=1/2$
and spatial concentration exponent $1/3$.
\end{corollary}
\begin{proof}
Take
\[
 \bar\Omega(x)=-\frac{x}{(1+x^2/2)^2},\qquad
 \bar U(x)=\frac{x}{\sqrt2(1+x^2/2)},\qquad
 \bar c_l=\frac1{8\sqrt2},\quad\bar c_\omega=-\frac3{8\sqrt2}.
\]
The elementary Hilbert pair for $(x^2+r^2)^{-1}$, differentiated
in $r$ and evaluated at $r=\sqrt2$, gives
$H\bar\Omega=\bar U'=\frac1{\sqrt2}(1-x^2/2)/(1+x^2/2)^2$.
Direct substitution verifies \eqref{sb:focus:profile}. The origin
jets are $(-1,6)$, $P>0$ and $\rho=4$, so all the profile hypotheses
hold. Theorem~\ref{sb:focus:theorem} gives the conclusion.
\end{proof}

\begin{corollary}[The analytic focusing interval]\label{sb:focus:parameter-range}
For every \(0<a\leq a_{\mathrm{foc}}\), the periodic gCLM equation admits
smooth odd initial data with finite-time focusing blowup as in
Theorem~\ref{sb:focus:theorem}.
\end{corollary}
\begin{proof}
By \cite[Proposition~6.2]{Profiles}, every normalized whole-line fixed
point in this parameter range has a positive spatial rate, including at
\(a=a_{\mathrm{foc}}\). The existence, regularity, and tail results in
\cite[Theorems~1.2 and~1.3(i)]{Profiles} supply a profile satisfying
Assumption~\ref{sb:focus:profile-input}, to which
Theorem~\ref{sb:focus:theorem} applies.
\end{proof}

\section{Critical focusing and logarithmic concentration}
\label{sb:paper:critical}

We use the equation and normalization of Section~\ref{sb:common:framework}.
The critical profile has zero stationary spatial rate, but the periodic
kernel produces a positive drift of the spatial scale. Its constant
order-\(L^{-2}\) correction cancels in the normalized shape equation;
the first shape source has order \(L^{-4}\). We prove that these facts,
together with a strictly coercive profile, give a smooth periodic blowup
with a logarithmic concentration law.

\subsection{The selected critical profile}
\begin{assumption}[A critical profile with a coercive estimate]\label{sb:dyn:profile-input}
There are \(a\in(1/2,3/4)\) and an even positive
\(f\in\mathcal S(\mathbb R)\) with
\(f(0)=1\), \(f''(0)=-2\), such that, writing
\[
 \bar\Omega=-xf,\quad \bar U'=H\bar\Omega,\quad \bar U(0)=0,\quad P=\bar U/x,
 \quad b=H\bar\Omega(0),\quad \bar c_\omega=-(1-a)b,
\]
we have \(b>0\), \(P>0\), and
\begin{equation}\label{sb:crit:stationary}
 a\bar U\bar\Omega'=(\bar U'+\bar c_\omega)\bar\Omega,\qquad \bar\Omega'(0)=-1,
       \qquad \bar\Omega'''(0)=6 .
\end{equation}
For the normalized operator \(\mathcal L\) defined in
\eqref{sb:crit:linear}, there are \(\kappa_0>0\) and an even positive
weight \(\rho_{\mathrm E}\) such that
\begin{equation}\label{sb:dyn:energy-input}
 \mathop{\rm Re}\int_{\mathbb R}\rho_{\mathrm E}\,\bar h\,\mathcal Lh\,\idiff x
       \leq-\kappa_0\int_{\mathbb R}\rho_{\mathrm E}|h|^2\,\idiff x
\end{equation}
for smooth odd compactly supported \(h\) with
\(h'(0)=h'''(0)=0\). The weight has the form
\(\rho_{\mathrm E}=1/(x^6 f_{\rm ref})\), where \(f_{\rm ref}>0\) is smooth and even,
bounded above and below near zero, and comparable to
\((1+x^2)^{-p_0}\) at infinity for some finite \(p_0>0\).
For \(j=1,2\), the functions
\((x\partial_x)^j\log\rho_{\mathrm E}(x)\) are bounded.

\end{assumption}

The profile theorem proved in Appendix~\ref{sb:sc:introduction}
supplies this input for one selected \(a\in(0.689,0.690)\), with
\(\kappa_0=1/20\) and \(\rho_{\mathrm E}=\rho_0\), where \(\rho_0\) is
the reference weight defined in \eqref{sb:sc:space}. This quantitatively
selected profile is the input for the dynamical construction below.

We record consequences needed below. Put
\[
 \mathsf b_j=\frac2\pi\int_0^\infty x^jf(x)\,\idiff x.
\]
Polynomial subtraction in the paired Hilbert kernel, using Schwartz
decay near the diagonal and in the remote input region, gives
\begin{equation}\label{sb:crit:tails}
 \bar U(x)=\mathsf b_2x^{-1}+O(x^{-3}),\qquad
 \bar U'(x)=-\mathsf b_2x^{-2}+O(x^{-4}),\qquad
 \lim_{x\to\infty}\frac{\log f(x)}{x^2}
               =\frac{\bar c_\omega}{2a\mathsf b_2}.
\end{equation}
The last limit follows by integrating
\(\bar\Omega'/\bar\Omega=(\bar U'+\bar c_\omega)/(a\bar U)\). The first two expansions hold after every
fixed number of logarithmic derivatives: use dilation covariance of
\(H\), apply the same kernel argument to \((x\partial_x)^j\bar\Omega\), and use
\(xP'=H\bar\Omega-P\).

For an odd perturbation \(h\), set
\[
 U_h(x)=\int_0^xHh(y)\,\idiff y .
\]
Define the two rate functionals by
\[
 B(h)=-\frac{3-a}{12}(Hh)''(0)-aHh(0),\qquad
 G(h)=-\frac{3-a}{12}(Hh)''(0)-Hh(0).
\]
We also set \(\alpha=(3-a)/6\).
The normalized linearization and its quadratic remainder are
\begin{align}
 \mathcal Lh={}&-a\bar U h'+(\bar U'+\bar c_\omega)h
       +(Hh+G(h))\bar\Omega-(aU_h+B(h)x)\bar\Omega',
       \label{sb:crit:linear}\\
 \mathcal N(h)={}&(Hh+G(h))h-(aU_h+B(h)x)h'.
       \label{sb:crit:quadratic}
\end{align}
The prescribed tangent conditions are \(h'(0)=h'''(0)=0\).
We use \(B_L,G_L\) for these functionals with \(H_L\) in place
of \(H\). They equal \(\ell_{l,L},\ell_{\omega,L}\) in
Lemma~\ref{sb:common:gauges}; hence the normalized evolution has
\(c_l=B_L(\Omega)\), \(c_\omega=G_L(\Omega)\).

\subsection{A closed operator on polynomially weighted spaces}

The singular energy is useful for excluding eigenvalues. For the periodic
nonlinear problem we use a fixed space with polynomial absolute weights.
Let \(z=\log x\), \(D=x\partial_x\), fix \(\sigma>0\), and define
\begin{equation}\label{sb:crit:polynomial-space}
 \nu=\frac92,\qquad w_\sigma(z)=e^{-\nu z}+e^{\sigma z},\qquad
 \|h\|_{Y_\sigma^m}
   =\|w_\sigma(z)h(e^z)\|_{H^m(\mathbb R_z)} .
\end{equation}
The odd space is the completion of smooth functions supported on annuli
away from zero; odd extension identifies it with the indicated weighted
space on the positive half-line. For \(m\geq1\), weighted embedding gives
\[
 |D^rh(x)|\leq C_m\|h\|_{Y_\sigma^m}
               (x^{-9/2}+x^\sigma)^{-1},\qquad 0\leq r<m.
\]
In particular \(Y_\sigma^2\) controls \(C^1\), and it controls
ordinary \(H^1\) when \(\sigma>1/2\).
The functionals \(B,G\) are bounded even on \(Y_\sigma^0\), since
\[
 Hh(0)=-\frac2\pi\int_0^\infty h(x)\frac{\idiff x}{x},
 \qquad
 (Hh)''(0)=-\frac4\pi\int_0^\infty h(x)\frac{\idiff x}{x^3}.
\]

\begin{lemma}[Closed realization and essential damping]\label{sb:crit:generator}
For every finite integer \(m\geq0\), \(\mathcal L\) generates a strongly
continuous group on \(Y_\sigma^m\), with domain
\begin{equation}\label{sb:crit:domain}
 \mathcal D(\mathcal L_m)
 =\{h\in Y_\sigma^m:\ (a\bar U/x)Dh\in Y_\sigma^m
                                  \text{ distributionally}\}.
\end{equation}
Its essential semigroup growth bound is strictly negative. For some
\(\delta\in(0,\min\{-\bar c_\omega,ab\})\), the half-plane
\(\operatorname{Re}\lambda>-\delta\) contains only finitely many
spectral points, all poles with finite algebraic multiplicity.
Their generalized eigenvectors are common to all finite \(m\) and
all \(\sigma>0\), are smooth and odd, have vanishing first and third derivatives at the origin, and are Schwartz.
For finitely many chosen levels there are equivalent Hilbert norms
\(\|\cdot\|_{m,*}\) such that
\begin{equation}\label{sb:crit:high-linear-energy}
 2\mathop{\rm Re}\langle \mathcal Lh,h\rangle_{m+1,*}
       \leq-c_m\|h\|_{m+1,*}^2+C_m\|h\|_{Y_\sigma^m}^2,
 \qquad c_m>0.
\end{equation}
\end{lemma}

\begin{proof}
Here are the operator and endpoint details. In logarithmic coordinates put
\[
 \mathsf B=a\bar U/x,\quad Q=\bar U'+\bar c_\omega,\quad
 O(z)=\bar\Omega(e^z),\quad O_1=\partial_zO .
\]
The transport is \(\mathcal T=-\mathsf B\partial_z+Q\).
At the two ends,
\[
 \begin{gathered}
 (\mathsf B,Q)=(ab,ab)+O(e^{2z})\quad(z\to-\infty),\\
 \mathsf B=a\mathsf b_2e^{-2z}+O(e^{-4z}),\quad
 Q=\bar c_\omega+O(e^{-2z})\quad(z\to+\infty).
 \end{gathered}
\]
Every fixed derivative has the corresponding bounds.
The two logarithmic convolution kernels for \(Hh\) and \(U_h/x\) are
\[
 K_H(t)=\frac2{\pi(e^{2t}-1)},\qquad
 K_U(t)=\frac{e^{-t}}{\pi}\log|\tanh(t/2)|,\qquad
 (\partial_t+1)K_U=K_H.
\]
The nonlocal operator is bounded after grouping it as
\begin{equation}\label{sb:crit:grouped}
 Rh=(Hh-Hh(0))\bar\Omega-a(U_h-xHh(0))\bar\Omega'
                 -\frac{3-a}{12}(Hh)''(0)(\bar\Omega-x\bar\Omega').
\end{equation}
For clarity, its origin estimate uses the expansions
\[
 O=-e^z+e^{3z}+O(e^{5z}),\quad
 O_1=-e^z+3e^{3z}+O(e^{5z}),\quad
 -\frac{3-a}{12}(O-O_1)
       +e^{2z}\left(\frac O2-\frac a6O_1\right)=O(e^{5z}).
\]
Subtract the constant and quadratic tails of \(K_H,K_U\) in the
region \(z<-1,\ \eta>z+1/2\), with a fixed smooth transition.
On the two separated negative-end regions the conjugated grouped
kernel is bounded by constant multiples of
\[
 e^z e^{-(\nu-2)(z-\eta)}\quad(\eta<z-1),\qquad
 \frac{e^{(5-\nu)z}(1+e^{-4\eta})}{w_\sigma(\eta)}
                                      \quad(\eta>z+1).
\]
These bounds have finite row and column integrals and vanishing tails
because \(4<\nu<5\). Higher Sobolev estimates use the total kernel
derivative: after integration by parts in the input,
\[
 \partial_z(RF)=R\partial_zF+R_1F,\qquad
 K_{R_1}=(\partial_z+\partial_\eta)K_R .
\]
This total derivative annihilates a kernel depending only on
\(z-\eta\), and preserves the grouped endpoint cancellations.
Iteration leaves the principal operator only on the top input
derivative; the commutators contain fewer input derivatives.
At the positive end, each kernel is a localized principal kernel,
an integrable remainder, and a constant one-sided tail. The tail acts
by an integral of \(h(e^\eta)\) over \(\eta>z\), bounded by
\(Ce^{-\sigma z}\|h\|_{Y_\sigma^0}\). The other parts have the usual
unweighted \(L^2\) bounds. Multiplication by \(O,O_1\), which decay
faster than every \(e^{-Nz}\), permits arbitrary positive \(\sigma\).
For example the right-backward envelope has the form
\(C_Ne^{-Nz}e^{(\sigma-2)(z-\eta)}\), \(z>\eta+1>1\), with
arbitrarily large \(N\). This also proves arbitrary gain of right
polynomial weight for \(Rh\), without claiming derivative gain.
The mixed-end regions follow by combining these bounds. Schur's
inequality, local Hilbert boundedness, and differentiation give boundedness
of \(R\) at every finite Sobolev level.

Its symmetric part is compact. On a compact square the only principal
singularity, after conjugation by any smooth positive weight \(v(z)\),
is \(v(z)O(z)/(v(\eta)\pi(z-\eta))\).
Symmetrization cancels its diagonal singularity. The velocity kernel is
logarithmic and gains one local derivative; lower derivative terms are
compact by local Rellich compactness. The separated source and target
tails just estimated vanish in operator norm. These facts prove compactness of the symmetric part of the grouped
operator on the full weighted space, including its higher derivative
norms. The principal Hilbert block itself remains a bounded operator of
order zero. The same compactness argument applies to \(v=\theta w_\sigma\) when \(\theta,\theta^{-1}\)
are smooth and bounded and constant at both ends.

The characteristic equation \(\dot z=\mathsf B(z)\) has a complete
flow, since \(\mathsf B\) is bounded and globally Lipschitz. Its Jacobian
and every fixed derivative are bounded on bounded time intervals.
The characteristic formula therefore defines a strongly continuous
transport group. Its maximal generator domain is
\(\{\mathsf B\partial_zF\in Y_\sigma^m\}\):
the characteristic fundamental theorem of calculus proves inclusion
in one direction; testing a generator difference quotient against compact
smooth functions proves the converse. Smooth compact functions form a
graph core. To see this, cut off at \(|z|\leq n\), using
\[
 \mathcal T(\chi_nF)-\chi_n\mathcal TF
                       =-\mathsf B\chi_n'F,
\]
then mollify on the resulting compact set. On that set \(\mathsf B>0\)
gives the needed local extra derivative. Globally division by
\(\mathsf B\) is not allowed: \eqref{sb:crit:domain} need not be one
full Sobolev derivative above the underlying space.
Bounded perturbation by \(R\) proves the asserted closed realization.

For the energy, the top transport coefficient at derivative order \(r\)
in the weight \(\theta w_\sigma\) is
\[
 Q-(r-\tfrac12)\mathsf B'
                   +\mathsf B(\log(\theta w_\sigma))'.
\]
Its endpoint limits without \(\theta\) are
\(-(\nu-1)ab\) and \(\bar c_\omega\).
Choose a compactly supported nonnegative function \(\chi\) and then
\(\theta(z)=\exp\bigl(-M\int_{-\infty}^z\chi(\zeta)\,\idiff\zeta\bigr)\) with \(M\) large.
Since \(\mathsf B\) has a positive minimum on the intervening compact
set, all the finitely many top coefficients become strictly negative.
This changes to an equivalent norm. Lower transport commutators have
coefficients tending to zero at both ends. Together with the compact
symmetric part of \(R\), they give a compact quadratic remainder.
Approximate that remainder by finite-rank forms generated by compact
smooth test functions. At level \(m+1\), integration of derivatives onto
these tests bounds the finite-rank part in \(Y_\sigma^m\); the small
approximation error is absorbed. This proves
\eqref{sb:crit:high-linear-energy}.

Subtracting a bounded compact operator representing the positive
quadratic remainder gives a generator with an exponentially decaying
semigroup. Bounded compact perturbations change its semigroup by a
compact operator at each positive time: approximate the perturbation by
finite-rank maps in Duhamel's formula; strong continuity of the two groups
and their Hilbert-space adjoints makes the resulting integrands norm
continuous. Thus the essential semigroup growth is negative.
The spectral decomposition of a semigroup with negative essential growth
then gives the stated finite set of poles in a smaller right half-plane.
Here one may use \cite[Chapter~IV, Corollary~2.11 and
Proposition~2.12]{EngelNagel2000}; the preceding argument
verifies its generator and essential-growth hypotheses.

We next justify the ordinary regularity needed to apply the singular
energy estimate to these spectral modes. Suppose that
\(\mathcal Lh=\lambda h\) with
\(\operatorname{Re}\lambda>-\delta\). On every compact interval of
the logarithmic coordinate, \(\mathsf B>0\), so the eigen-equation
first gives one local Sobolev derivative. The local Hilbert bounds and
the smoothness of the separated kernels then permit iteration. Near
the left endpoint, \(\mathsf B\to ab>0\); the same argument with
cutoffs on a negative half-line gives all logarithmic derivatives in
the origin weight. Thus weighted Sobolev embedding gives the initial
conormal bound \(D^j h=O(x^{9/2})\) at zero for every fixed \(j\).
For the right endpoint, put \(s_x=1/(a\bar U)\).
At the right end \(s\asymp x^2\) and
\[
 \partial_sh+(\lambda-Q)h=Rh,\qquad
 \mathop{\rm Re}(\lambda-Q)>c>0
\]
for large \(x\). The forcing has all polynomial right weights, by the
kernel estimate above. Fix a sufficiently large \(s_0>0\), and express the functions in this
characteristic coordinate in the next two displays. Variation of
constants gives
\[
 h(s)=e^{-\int_{s_0}^{s}(\lambda-Q(\xi))\,\idiff\xi}h(s_0)
 +\int_{s_0}^{s}e^{-\int_r^{s}(\lambda-Q(\xi))\,\idiff\xi}
                       Rh(r)\,\idiff r .
\]
The kernel is bounded by \(e^{-c(s-r)}\mathbf1_{r<s}\). For every
\(N\geq0\), the inequality
\((1+s)^N\leq(1+r)^N(1+s-r)^N\), followed by Young's inequality,
yields
\[
 \|(1+s)^Nh\|_{L^2(s_0,\infty)}
 \leq C_N\bigl(|h(s_0)|+
                 \|(1+s)^NRh\|_{L^2(s_0,\infty)}\bigr).
\]
Since \(s\asymp x^2\), this proves every right polynomial weight
for \(h\). Returning to
\(Dh=\mathsf B^{-1}((Q-\lambda)h+Rh)\) costs only two powers of \(x\).
Iteration gives all conormal derivatives and all right polynomial weights.

To complete the origin regularity argument, let an odd function have a
conormal remainder of noninteger order \(r>0\). Subtracting its
convergent Hilbert moments gives an even Taylor polynomial plus a
conormal remainder of order \(r\). This follows by splitting at \(x/2\)
and \(3x/2\), with principal-value subtraction in the middle interval.
Multiplication by \(\bar\Omega=O(x)\), or integration from zero followed by
multiplication by \(\bar\Omega'\), gives one further order in \(Rh\).
The constant Hilbert coefficient is absent from
\eqref{sb:crit:grouped}. If \(c_2x^2\) denotes the quadratic Hilbert
term, its contribution to \(Rh\) is
\[
 c_2\bigl(x^2\bar\Omega-ax^3\bar\Omega'/3-\alpha(\bar\Omega-x\bar\Omega')\bigr).
\]
The coefficient of \(x^3\) is
\(c_2(-1+a/3+2\alpha)=0\), since \(2\alpha=1-a/3\).
Thus the smooth odd Taylor part starts at order five, while the
noninteger remainder gains one order. The local eigen-equation is
\[
 Dh-c(x)h=g(x),\quad
 c(0)=1-\lambda/(ab),\quad g=Rh/\mathsf B .
\]
Here \(\mathop{\rm Re}c(0)<2<9/2\). Multiplying by a smooth even
integrating factor \(\mu\) gives
\[
 \mu(x)h(x)
   =x^{c(0)}\int_0^x t^{-c(0)-1}\mu(t)g(t)\,\idiff t .
\]
The homogeneous term is excluded by the origin weight. The grouped
linear and cubic Taylor coefficients vanish. Each higher odd monomial
is divided by \(k-c(0)\neq0\), and the remainder gains the next order.
Iteration through half-integer orders proves ordinary smoothness and
the vanishing first and third derivatives at the origin. The conormal right estimates imply Schwartz regularity.
Induction treats generalized eigenvectors and smooth odd Schwartz
inhomogeneous terms with vanishing first and third derivatives at the origin in the same way.
For completeness, the spectral realizations are compatible because
the characteristic groups and their bounded-perturbation expansions
agree on common data. For \(\operatorname{Re}\lambda\) sufficiently
large, their resolvents therefore agree by
\[
 (\lambda-\mathcal L)^{-1}h
   =\int_0^\infty e^{-\lambda\tau}e^{\tau\mathcal L}h\,\idiff\tau .
\]
The regularity just proved places every generalized eigenspace in the
intersection of all the stated spaces. Meromorphic continuation of the
resolvents, and integration over a common resolvent contour, then
identify the Riesz projections. This completes the proof.
\end{proof}

\begin{proposition}[From scalar coercivity to a stable semigroup]
\label{sb:crit:spectral-gap}
Under Assumption~\ref{sb:dyn:profile-input}, for every finite \(m\) and
\(\sigma>0\) there are \(C_m,\eta_m>0\) such that
\begin{equation}\label{sb:crit:semigroup}
 \|e^{\tau\mathcal L}h\|_{Y_\sigma^m}
                      \leq C_me^{-\eta_m\tau}\|h\|_{Y_\sigma^m},
 \qquad \tau\geq0 .
\end{equation}
In particular zero is in the resolvent of this normalized realization.
\end{proposition}

\begin{proof}
By Lemma~\ref{sb:crit:generator}, a spectral point with nonnegative
real part would have a smooth Schwartz eigenvector \(h\) with two zero
jets. For the polynomial reference weight in
Assumption~\ref{sb:dyn:profile-input}, these facts immediately give
\(h\in L^2(\rho_{\mathrm E})\). More explicitly, choose
\(\sigma>\max\{0,p_0-5/2\}\). Then \(Y_\sigma^0\) embeds into
\(L^2(\rho_{\mathrm E})\): its measure has density comparable to \(x^{-10}\)
at zero and \(x^{2\sigma-1}\) at infinity, whereas \(\rho_{\mathrm E}\) has
orders \(x^{-6}\) and \(x^{2p_0-6}\), respectively.
The eigenvector belongs to this stronger space by the preceding lemma.
Choose smooth even cutoffs that remove \(|x|<\varepsilon\)
and \(|x|>R\). The corresponding truncated eigenvectors converge in
\(Y_\sigma^0\). Their transport commutators have the form
\(-a\bar U\chi'h\). At zero, \(\bar U=O(x)\) and \(h=O(x^5)\);
at infinity, \(\bar U=O(x^{-1})\) and \(h\) is Schwartz. These bounds
make the commutators tend to zero in \(Y_\sigma^0\) as
\(\varepsilon\downarrow0\) and \(R\to\infty\).
The boundedness of the grouped operator \(R\) proves convergence
of the nonlocal terms in the same space. The embedding into
\(L^2(\rho_{\mathrm E})\) therefore gives convergence of both the functions
and their images under \(\mathcal L\) in the energy space. Hence \eqref{sb:dyn:energy-input} applies to \(h\)
and contradicts a nonnegative real part of its eigenvalue.

The finite collection of poles in a right half-plane
contains none with nonnegative real part. Its distance from the imaginary
axis, together with the negative essential growth bound, gives
\eqref{sb:crit:semigroup}, with any sufficiently small positive
\(\eta_m\).
\end{proof}

\subsection{Periodic geometry and a uniform two-level estimate}

In the remaining subsections write \(Y_m=Y_\sigma^m\), with a fixed
\(\sigma>6\). Let \(\chi\) be even and smooth, with \(0\leq\chi\leq1\), equal to one on
\([-1/4,1/4]\), and supported in \((-1/2,1/2)\). Set
\[
 \Phi_L(x)=\chi(x/L)\bar\Omega(x),\qquad
 q_0=\frac{\pi^2\mathsf b_2}{12},\quad k=aq_0,\quad
 r_0=\frac{\pi^4\mathsf b_4}{720},\quad
 r_2=\frac{\pi^4\mathsf b_2}{240}.
\]
For an odd smooth function supported strictly inside \((-L,L)\),
write
\[
 \mathcal F_L(v)=(G_L(v)+H_Lv)v-(B_L(v)x+aU_L(v))v',
 \qquad U_L(v)(x)=\int_0^xH_Lv(y)\,\idiff y .
\]
Its periodic lift and its zero extension agree inside that cell.
When a product containing \(v\) or its derivative is extended to the line
below, it is extended by zero outside the cell.

\begin{lemma}[The critical source in the fixed weighted space]
\label{sb:crit:periodic-source}
\label{sb:prop:critical-correction}
Put \(k_4=ar_0+\alpha r_2\) and
\[
 S=r_2\{(\alpha-x^2)\bar\Omega-(\alpha x-ax^3/3)\bar\Omega'\}.
\]
For every fixed \(m\), as \(L\to\infty\),
\begin{align}
 B_L(\Phi_L)&=kL^{-2}+k_4L^{-4}+O(L^{-6}),
                                  \label{sb:crit:rate-l}\\
 G_L(\Phi_L)&=\bar c_\omega+q_0L^{-2}
                    +(r_0+\alpha r_2)L^{-4}+O(L^{-6}),
                    \label{sb:crit:rate-w}\\
 \|\mathcal F_L(\Phi_L)-L^{-4}S\|_{Y_m}
                      &\leq C_mL^{-6}.
                      \label{sb:crit:source}
\end{align}
Also \(\|\Phi_L-\bar\Omega\|_{Y_m}\leq C_{m,N}L^{-N}\) for every \(N\).
\end{lemma}

\begin{proof}
The cotangent expansion, with odd cancellation in the input, gives
\begin{equation}\label{sb:crit:kernel-expansion}
 H_L\Phi_L-H\bar\Omega
       =-q_0L^{-2}-(r_0+r_2x^2)L^{-4}
                       +O\bigl(L^{-6}(1+|x|^4)\bigr)
\end{equation}
on the core \(|x|\leq L/4\). The estimate holds with every fixed number
of logarithmic derivatives and with the corresponding ordinary Taylor
remainders at zero. Indeed,
\[
 \frac1{2L}\cot\frac{\pi z}{2L}
 =\frac1{\pi z}-\frac{\pi z}{12L^2}
                      -\frac{\pi^3z^3}{720L^4}+O(z^5L^{-6}),
\]
and the remainder after pairing the odd input is bounded by
\[
 CL^{-6}\int_{\mathbb R}|y|(|x|+|y|)^4|\bar\Omega(y)|\,\idiff y.
\]
The moments of the removed tails have arbitrary polynomial smallness.
The value and second derivative at zero give
\eqref{sb:crit:rate-l}--\eqref{sb:crit:rate-w}.

There is an exact cancellation that must precede the shape estimate.
Adding a spatial constant \(d\) to \(H_Lv\) adds \(d x\) to \(U_L(v)\)
and subtracts \(ad,d\) from \(B_L(v),G_L(v)\), respectively. In fact,
\[
 \bigl(B_L(v)-ad\bigr)x+a\bigl(U_L(v)+d x\bigr)
       =B_L(v)x+aU_L(v),\qquad
 G_L(v)-d+H_Lv+d=G_L(v)+H_Lv.
\]
Both brackets defining \(\mathcal F_L(v)\) remain unchanged, so the order
\(L^{-2}\) correction in \eqref{sb:crit:kernel-expansion} contributes only
to the rates.
The remaining bracket corrections at order \(L^{-4}\) give exactly \(S\).
Their combined Taylor expansion at zero starts at order \(x^5\).
After multiplication by \(\bar\Omega,\bar\Omega'\), the differentiated remainder is in
\(Y_m\) with bound \(C_mL^{-6}\).
The cutoff annulus and all its derivatives contribute \(O(L^{-N})\)
for every \(N\), by the Schwartz bounds. The same argument shows
\(\mathcal F_\infty(\Phi_L)=O_{Y_m}(L^{-N})\) and
the final assertion. This proves the globally weighted statement.
\end{proof}

\begin{lemma}[Uniform estimates on the entire growing cell]
\label{sb:crit:periodic-estimates}
Let a smooth odd normalized periodic solution have support strictly
inside its cell \((-L,L)\), and put
\[
 h=\Omega_{\rm ext}-\bar\Omega.
\]
For each fixed \(m\geq4\), its exact equation on the line has the form
\begin{equation}\label{sb:crit:line-evolution}
 h_\tau=\mathcal Lh+\mathcal N(h)
                +\bigl(\mathcal F_L-\mathcal F_\infty\bigr)
                                      (\Omega_{\rm ext}) .
\end{equation}
There is no derivative of a moving cutoff in this equation.
While \(L\geq L_0>1\) and \(\|h\|_{Y_{m+1}}\) is small, its forcing
relative to \(\mathcal L\) satisfies
\begin{equation}\label{sb:crit:base-forcing}
 \|h_\tau-\mathcal Lh\|_{Y_m}
 \leq C_mL^{-4}
       +C_m\bigl(L^{-2}+\|h\|_{Y_{m+1}}\bigr)\|h\|_{Y_{m+1}} .
\end{equation}
In the equivalent high norm of Lemma~\ref{sb:crit:generator},
\begin{equation}\label{sb:crit:nonlinear-high}
 \frac{\idiff}{\idiff\tau}\|h\|_{m+1,*}^2
 \leq-c\|h\|_{m+1,*}^2+C\|h\|_{Y_m}^2+CL^{-8}
       +C\bigl(L^{-2}+\|h\|_{Y_{m+1}}\bigr)\|h\|_{m+1,*}^2 .
\end{equation}
Constants are independent of the distance from the support to the seam.
The rate bounds are
\begin{equation}\label{sb:crit:rate-errors}
 |B_L(\Omega)-kL^{-2}|
 +|G_L(\Omega)-\bar c_\omega-q_0L^{-2}|
             \leq C\bigl(\|h\|_{Y_0}+L^{-4}\bigr).
\end{equation}
\end{lemma}

\begin{proof}
We supply the uniform kernel and derivative estimates, including the
region near the seam. For a compact cell input \(v\), the paired kernels are
\[
 H_Lv(x)=\frac1{2L}\int_0^L v(y)
 \left\{\cot\frac{\pi(x-y)}{2L}
                   -\cot\frac{\pi(x+y)}{2L}\right\}\,\idiff y ,
\]
\[
 U_L(v)(x)=\frac1\pi\int_0^L v(y)
       \log\left|\frac{\sin(\pi(x-y)/(2L))}
                       {\sin(\pi(x+y)/(2L))}\right|\,\idiff y .
\]
Cauchy--Schwarz in \(Y_0\) controls
\[
 \int_0^\infty x|v(x)|\,\idiff x
       +\int_0^\infty |v(x)|x^{-3}\,\idiff x
                         \leq C\|v\|_{Y_0}.
\]
The two rate kernels at zero are bounded by \(C/x\) and \(C/x^3\)
on \((0,L)\). Hence \(|B_L(v)|+|G_L(v)|\leq C\|v\|_{Y_0}\).

On the core \(0<x\leq L/2\), the difference
\(\Delta_Lv=H_Lv-Hv\) has a smooth kernel away from the non-diagonal
image poles. Oddness in the input and evenness in the output give
\begin{equation}\label{sb:crit:core-difference}
 |\Delta_Lv(0)|\leq CL^{-2}M_1^*(v),\qquad
 |D^r(\Delta_Lv(x)-\Delta_Lv(0))|
                   \leq C_rL^{-4}x^2M_1^*(v),
 \quad M_1^*(v)=\int_0^L y|v(y)|\,\idiff y .
\end{equation}
For the second inequality subtract the value at zero from the paired
smooth kernel. Taylor's formula supplies \(x^2yL^{-4}\); logarithmic
derivatives preserve this bound. Integration from zero gives the matching
primitive bound. The line convolution kernels in the preceding proof
give unweighted Hilbert and primitive bounds on the \(z\) axis, including
their bounded one-sided tail integrals.

For \(L/2\leq|x|\leq L\), split the input using a fixed cutoff into
\(v_{\rm in}\), supported in \(|x|<3L/8\), and \(v_{\rm out}\),
supported in \(|x|>L/4\). The separated inner input satisfies
\[
 |D^rH_Lv_{\rm in}(x)|
       +|D^r(U_L(v_{\rm in})(x)/x)|
                              \leq C_rL^{-2}M_1^*(v).
\]
For the outer input use \(\theta=x/L\) on the fixed torus of length two.
There the weight is comparable to \(L^\sigma\), and \(D\) derivatives
are equivalent to ordinary \(\theta\) derivatives on its support. Thus
\[
 \|v_{\rm out}(L\,\cdot)\|_{H^r(\mathbb T_1)}
                                \leq C_rL^{-\sigma}\|v\|_{Y_r}.
\]
Ordinary periodic derivatives commute with \(H_L\), whose Fourier
multiplier has modulus at most one; its mean-zero primitive gains a
derivative. Division by \(\theta\) on the output arc is harmless.
The core and outer estimates together yield, for \(r\geq1\),
\begin{equation}\label{sb:crit:coefficient-bounds}
 \begin{split}
 &\|H_Lv\|_\infty+\|DH_Lv\|_{H^{r-1}(\mathrm{d}z)}
   +\|U_L(v)/x\|_\infty\\
 &\hspace{35mm}
       +\|D(U_L(v)/x)\|_{H^r(\mathrm{d}z)}
                           \leq C_r\|v\|_{Y_r}.
 \end{split}
\end{equation}
The \(z\) norms in this display are on \(z<\log L\).
This proof uses no lower bound on the support gap, and uses ordinary
periodic derivatives on the outer arc rather than commuting a
nonperiodic logarithmic derivative through the Hilbert transform.

Here is the normalized geometric comparison. Put
\[
 \delta d(v)=-\frac{3-a}{12}(\Delta_Lv)''(0),\qquad
 A_L(v)=\delta d(v)+\frac a x\int_0^x
                  \bigl(\Delta_Lv(y)-\Delta_Lv(0)\bigr)\,\idiff y ,
\]
\[
 C_L(v)=\delta d(v)+\Delta_Lv(x)-\Delta_Lv(0).
\]
Then exactly
\begin{equation}\label{sb:crit:geometric-field}
 (\mathcal F_L-\mathcal F_\infty)(v)
                             =C_L(v)v-A_L(v)Dv.
\end{equation}
Because the input \(\Phi_L\) is supported in half the cell,
\eqref{sb:crit:core-difference} holds for it on the full output cell.
It gives
\begin{equation}\label{sb:crit:geometric-transport}
 |D^rA_L(\Phi_L)(x)|+|D^rC_L(\Phi_L)(x)|
             \leq C_rL^{-4}(1+x^2)\leq C_rL^{-2}.
\end{equation}
The other polarized term satisfies the stronger estimate
\begin{equation}\label{sb:crit:geometric-smoothing}
 \|C_L(v)\Phi_L-A_L(v)D\Phi_L\|_{Y_r}
                                      \leq C_rL^{-4}\|v\|_{Y_0}.
\end{equation}
Its output stays in half the cell. At zero, write
\(\Delta_Lv=d_0+d_2x^2/2+O(L^{-6}M_1^*(v)x^4)\).
The constant coefficients of \(A_L,C_L\) are both
\(-(3-a)d_2/12\); their quadratic coefficients are \(ad_2/6,d_2/2\).
Substitution of \(\Phi_L=-x+x^3+O(x^5)\) cancels the cubic terms in
\eqref{sb:crit:geometric-smoothing}. The result is
\(O(L^{-4}M_1^*(v)x^5)\), which is admissible for \(\nu<5\).
Away from zero the Schwartz profile absorbs every polynomial factor.

Write \(\mathcal N_L\) for \eqref{sb:crit:quadratic} with the periodic
operators and rate functionals. The coefficient bounds
\eqref{sb:crit:coefficient-bounds} also prove
the tame nonlinear estimates
\[
 \|\mathcal N_L(v)\|_{Y_r}
       \leq C_r\|v\|_{Y_r}\|v\|_{Y_{r+1}},\qquad r\geq2,
\]
\[
 |\mathop{\rm Re}\langle\mathcal N_L(v),v\rangle_{r,*}|
                     \leq C_r\|v\|_{Y_2}\|v\|_{Y_r}^2 .
\]
For the latter, conjugate by the fixed weight and integrate the top
transport derivative once by parts. The commutator inequality
\[
 \|[\partial_z^r,A]\partial_zZ\|_2
 \leq C_r\{\|A_z\|_\infty\|Z\|_{H^r}
                +\|A_z\|_{H^{r-1}}\|Z_z\|_\infty\}
\]
controls the remaining derivatives. Every logarithmic derivative of the
weight divided by that weight is bounded. The same estimates hold for
the line operator. A Sobolev extension on a translated interval near
\(\log L\) has uniform constants; the compact cell input vanishes near
that endpoint, so integration by parts has no boundary contribution.

Now \(\mathcal F_\infty(\bar\Omega)=0\), so
\(\mathcal F_\infty(\bar\Omega+h)=\mathcal Lh+\mathcal N(h)\).
Outside the cell \(\bar\Omega+h=\Omega_{\rm ext}=0\), making this sum zero.
The zero extension remains smooth across the moving seam since it has
a support gap. This proves the exact identity
\eqref{sb:crit:line-evolution}.
To estimate its last term, polarize
\eqref{sb:crit:geometric-field} at
\(\Omega_{\rm ext}=\Phi_L+v\), where
\[
 v=h+(\bar\Omega-\Phi_L).
\]
The profile source is \(L^{-4}S+O_{Y_r}(L^{-6})\) by
Lemma~\ref{sb:crit:periodic-source}. The first polarized term is the
small transport \eqref{sb:crit:geometric-transport}; the other is
\eqref{sb:crit:geometric-smoothing}. The quadratic terms have the just
proved tame bounds.

One may extend the profile transport coefficients smoothly past the cell
and cut them off before \(5L/4\), preserving all \(O(L^{-2})\) bounds:
the profile input is separated from their next kernel pole. Replacing
\(v\) by \(h\) gives rapidly small sources and rapidly small coefficients
times \(Dh\). The latter are retained as transport coefficients when
taking high energy; they are not treated as derivative-free forcing.
For the quadratic energy pairing, put \(\eta_L=\bar\Omega-\Phi_L\), so
\(h=v-\eta_L\), and use
\[
 \begin{split}
 \langle\mathcal N_L(v)-\mathcal N_\infty(v),h\rangle_{r,*}
 ={}&\langle\mathcal N_L(v)-\mathcal N_\infty(v),v\rangle_{r,*}\\
 &-\langle\mathcal N_L(v)-\mathcal N_\infty(v),\eta_L\rangle_{r,*}.
 \end{split}
\]
The first term has the compact-input estimate already proved. In the
second integrate the top derivative onto the known Schwartz tail
\(\eta_L\). Since \(v\) vanishes near the seam, there is no boundary
term, and the result is bounded by
\(C_{r,N}L^{-N}(\|v\|_{Y_r}^2+\|v\|_{Y_r})\).
Combining these observations proves \eqref{sb:crit:base-forcing}.
For high energy use \eqref{sb:crit:high-linear-energy}, integrate each
small transport at its top derivative, and estimate the fixed source by
Young's inequality. This gives \eqref{sb:crit:nonlinear-high}.
Finally linearity of \(B_L,G_L\), their uniform moment bounds, the rapid
cutoff error, and \eqref{sb:crit:rate-l}--\eqref{sb:crit:rate-w}
give \eqref{sb:crit:rate-errors}.
\end{proof}

\subsection{Smooth periodic blowup and its concentration law}

\begin{theorem}[Smooth critical periodic blowup]
\label{sb:crit:transfer}
Let \(\bar\Omega\) and \(a\) satisfy Assumption~\ref{sb:dyn:profile-input},
and fix a finite integer \(m\geq4\).
For every sufficiently large \(L_0\), depending on \(m\),
the smooth odd periodic initial data
\begin{equation}\label{sb:crit:physical-data}
 \omega_0(y)=\chi(y/\pi)\,
                    \bar\Omega\!\left(\frac{L_0}{\pi}y\right),
                       \qquad -\pi\leq y\leq\pi ,
\end{equation}
produce a smooth solution up to a finite maximal time \(T\).
There are positive rescaling factors with
\[
 C_l(0)=L_0/\pi,\quad C_\omega(0)=1,\quad
 \Omega(x,\tau)=C_\omega(\tau)^{-1}
             \omega(C_l(\tau)^{-1}x,t(\tau)),\qquad L=\pi C_l,
\]
such that
\begin{equation}\label{sb:crit:convergence}
 \|\Omega_{\rm ext}(\tau)-\bar\Omega\|_{Y_m}
                         \leq C_m L(\tau)^{-4},\qquad
 L(\tau)^2\asymp L_0^2+\tau .
\end{equation}
In particular \(c_l\to0\) and \(C_l\to+\infty\). As \(t\uparrow T\),
\begin{equation}\label{sb:crit:leading-law}
 C_l(t)^2=\frac{a\mathsf b_2}{6(1-a)b}
               \log\frac1{T-t}
                    +O\!\left(\log\log\frac1{T-t}\right),
\end{equation}
and
\begin{equation}\label{sb:crit:physical-limits}
 (T-t)\|\omega(t)\|_\infty\longrightarrow-\frac{\|\bar\Omega\|_\infty}{\bar c_\omega},
 \qquad
 (T-t)\omega(C_l(t)^{-1}x,t)\longrightarrow(-\bar c_\omega)^{-1}\bar\Omega(x)
\end{equation}
locally uniformly in \(x\). The rescaled zero extension converges
globally in \(C^1(\mathbb R)\) as in \eqref{sb:crit:convergence}.
\end{theorem}

\begin{proof}
Start from the smooth periodic local solution of
Proposition~\ref{sb:gwp:local} and use the two
origin gauges above, together with
\begin{equation}\label{sb:crit:scales}
 L'/L=c_l=B_L(\Omega),\qquad
 C_\omega'/C_\omega=-c_\omega=-G_L(\Omega),\qquad
 t'=C_\omega^{-1}.
\end{equation}
Initially \(\Omega=\Phi_{L_0}\), with the exact prescribed jets.
Equivalently these equations define a coupled normalized periodic problem
on the fixed physical torus. The scales remain positive and finite on
every interval where their rates are bounded.

We close the two derivative levels explicitly. Write
\(X=\|h\|_{Y_m}\), \(Z=\|h\|_{m+1,*}\), where \(h=\Omega_{\rm ext}-\bar\Omega\).
Use a common fixed equivalence constant for the high norms. Bootstrap
\begin{equation}\label{sb:crit:bootstrap}
 X\leq K_0L^{-4},\qquad Z\leq K_1L^{-4},\qquad L\geq L_0 .
\end{equation}
Equation \eqref{sb:crit:rate-errors}, with large \(L_0\), gives
\begin{equation}\label{sb:crit:positive-drift}
 \frac{k}{2L^2}\leq c_l\leq\frac{2k}{L^2},\qquad
 (L^2)'=2k+O(L^{-2}),\qquad
 c_\omega=\bar c_\omega+q_0L^{-2}+O(L^{-4}).
\end{equation}
In particular \(0\leq(\log L)'\leq2kL_0^{-2}\).
For every fixed \(p>0\) and \(\eta>0\), increasing \(L_0\) if needed,
\begin{equation}\label{sb:crit:convolution}
 \int_0^\tau e^{-\eta(\tau-s)}L(s)^{-p}\,\idiff s
       \leq\frac{2}{\eta}L(\tau)^{-p},\qquad
 e^{-\eta\tau}L_0^{-p}\leq L(\tau)^{-p}.
\end{equation}
Indeed \(L(s)^{-p}\leq L(\tau)^{-p}
\exp(2pkL_0^{-2}(\tau-s))\).

Duhamel's formula, \eqref{sb:crit:semigroup}, and
\eqref{sb:crit:base-forcing} therefore imply
\begin{equation}\label{sb:crit:base-improvement}
 X(\tau)\leq
 \left(\epsilon_0+C+
        C K_1L_0^{-2}+C K_1^2L_0^{-4}\right)L(\tau)^{-4},
\end{equation}
where \(\epsilon_0\) can be made arbitrarily small since
\(\|\Phi_{L_0}-\bar\Omega\|_{Y_{m+1}}=O(L_0^{-N})\) for every \(N\).
After absorbing the last term of
\eqref{sb:crit:nonlinear-high}, that inequality and
\eqref{sb:crit:convolution} give
\begin{equation}\label{sb:crit:high-improvement}
 Z(\tau)^2\leq
             \bigl(\epsilon_0^2+C(K_0^2+1)\bigr)L(\tau)^{-8}.
\end{equation}
Choose first \(K_0\) larger than four times the fixed source constant
in \eqref{sb:crit:base-improvement}; next choose \(K_1\) so that
the right side of \eqref{sb:crit:high-improvement} is less than
\(K_1^2L^{-8}/4\). Finally take \(L_0\) large, with those constants fixed.
This makes the last two terms in
\eqref{sb:crit:base-improvement} small, absorbs the high nonlinear and
transport coefficients, and verifies \eqref{sb:crit:positive-drift}.
Both norm bounds improve strictly. This order of choices avoids a
circular dependence between the two regularity levels. Continuity and
the integral inequalities prevent exit from the bootstrap region.

The support condition used in the estimates persists for every finite
rescaled time. On the fixed physical circle the solution obeys a
transport equation with multiplication on the right. Hence a zero
neighborhood is carried to a zero neighborhood by its characteristic
flow. Oddness makes the antipodal seam a fixed characteristic for the
physical velocity. A Lipschitz flow cannot close the initial positive
gap in finite time.

For completeness, the norm bounds also give continuation of the smooth
solution at each finite rescaled time. The fixed exponent \(9/2\) at the
origin and \(\sigma>6\) at infinity give bounded ordinary \(H^3\) norm
from \(Y_m\), \(m\geq4\), after adding the fixed smooth profile.
The finite positive scales transfer this bound to the physical circle.
Differentiating the periodic equation, using commutation of ordinary
derivatives with the Hilbert transform and the product commutator
inequality, gives for each higher integer \(r\)
\[
 \frac{\idiff}{\idiff t}\|\omega\|_{H^r}
          \leq C_r\|\omega\|_{H^3}\|\omega\|_{H^r}.
\]
Thus the higher norms of smooth initial data remain finite on each
bounded rescaled-time interval, and local existence extends the solution.
This proves global existence in \(\tau\) and
\eqref{sb:crit:convergence}.

Integration of \eqref{sb:crit:positive-drift} gives, for fixed \(L_0\),
\begin{equation}\label{sb:crit:time-asymptotics}
 L^2=2k\tau+O(\log(2+\tau)),\qquad
 C_\omega=A_\omega e^{-\bar c_\omega\tau}\tau^{-1/(2a)}
        \left(1+O\!\left(\frac{\log\tau}{\tau}\right)\right).
\end{equation}
The second formula uses \(q_0/(2k)=1/(2a)\), with \(A_\omega>0\).
It follows that the physical terminal time is finite:
\(T=\int_0^\infty C_\omega^{-1}\,\idiff\tau<\infty\). Moreover
\[
 (T-t)C_\omega=(-\bar c_\omega)^{-1}(1+O(\tau^{-1})),\qquad
 -\bar c_\omega\tau=\log\frac1{T-t}
                 +\frac1{2a}\log\log\frac1{T-t}+O(1).
\]
For the first identity, divide the tail integral
\(\int_\tau^\infty C_\omega(s)^{-1}\,\idiff s\) by
\(C_\omega(\tau)^{-1}\); the integrand ratio at \(s=\tau+r\)
is \(e^{\bar c_\omega r}\exp(O(r/\tau))\), allowing exponential domination.
Substitution in \(C_l^2=L^2/\pi^2\) proves
\eqref{sb:crit:leading-law}. The embedding \(Y_m\hookrightarrow C^1\),
the exact rescaling, and \((T-t)C_\omega\to(-\bar c_\omega)^{-1}\) give
\eqref{sb:crit:physical-limits}. The limiting supremum is positive,
so the physical solution cannot continue smoothly through \(T\).
\end{proof}

\begin{proposition}[The next logarithmic coefficient]\label{sb:crit:second-law}
Under the hypothesis of Theorem~\ref{sb:crit:transfer}, define the unique
normalized line corrector
\begin{equation}\label{sb:crit:corrector}
 Q_4=-\mathcal L^{-1}S,\qquad
 \varkappa=k_4+B(Q_4).
\end{equation}
It is smooth, odd, Schwartz, and has vanishing first and third derivatives at the origin.
For fixed sufficiently large \(L_0\), as \(\tau\to\infty\),
\begin{equation}\label{sb:crit:shape-correction}
 \Omega_{\rm ext}-\bar\Omega=L^{-4}Q_4+O_{Y_m}(L^{-6}).
\end{equation}
The assertion holds at the prescribed finite level \(m\) in
Theorem~\ref{sb:crit:transfer}. With \(s=\log(1/(T-t))\),
\begin{equation}\label{sb:crit:second-scale}
 C_l(t)^2=-\frac{2k}{\pi^2\bar c_\omega}s
       +\frac1{\pi^2}\left(-\frac{k}{a\bar c_\omega}
                                  +\frac{\varkappa}{k}\right)\log s
       +O(1).
\end{equation}
\end{proposition}

\begin{proof}
The stable semigroup makes zero a resolvent point. Lemma
\ref{sb:crit:generator}, applied to the smooth source \(S\), gives
the stated regularity of \(Q_4\).
The high norm in the two-level construction gives
\(\|h\|_{Y_{m+1}}=O(L^{-4})\). The sharper source expansion
in Lemma~\ref{sb:crit:periodic-source}, the transport bound
\eqref{sb:crit:geometric-transport}, and the quadratic estimates give
\[
 h_\tau=\mathcal Lh+L^{-4}S+O_{Y_m}(L^{-6}).
\]
Since \((L^{-4})'=-4c_lL^{-4}=O(L^{-6})\), subtraction of
\(L^{-4}Q_4\) leaves an exponentially stable equation with
\(O(L^{-6})\) forcing. Duhamel's formula gives more explicitly
\[
 \|h-L^{-4}Q_4\|_{Y_m}
 \leq C e^{-\eta'\tau}L_0^{-4}+C L(\tau)^{-6}
\]
for some \(\eta'>0\). The initial layer is necessary for the cutoff
data \eqref{sb:crit:physical-data}; it is negligible in
\eqref{sb:crit:shape-correction}.

The rate functional on a compact input satisfies
\[
 |B_L(v)-B(v)|\leq CL^{-2}\|v\|_{Y_0}.
\]
This follows from the same origin kernel expansion used above.
Apply this comparison to the compact input
\(v=\Omega_{\rm ext}-\Phi_L=h+(\bar\Omega-\Phi_L)\), whose two jets vanish.
The difference \(v-h\) is rapidly small in every fixed \(Y_m\).
Combining the comparison with \eqref{sb:crit:rate-l} and
\eqref{sb:crit:shape-correction} yields
\[
 c_l=kL^{-2}+\varkappa L^{-4}+O(L^{-6}).
\]
Thus \(y=L^2\) satisfies
\[
 y'=2k+\frac{2\varkappa}{y}+O(y^{-2}),\qquad
 y=2k\tau+\frac{\varkappa}{k}\log\tau+O(1).
\]
In the second equality the derivative of the remaining error is
integrable, using first \(y=2k\tau+O(\log\tau)\).
Finally
\(-\bar c_\omega\tau=s+(2a)^{-1}\log s+O(1)\) from
\eqref{sb:crit:time-asymptotics}. Substitution gives
\eqref{sb:crit:second-scale}.
\end{proof}

\section{Smooth blowup with a fixed spatial scale}
\label{sb:finite:section}

We realize finite-period profiles by solutions with smooth
periodic initial data. An adapted characteristic coordinate reduces the
principal normalized operator to constant transport with damping. We show
that the symmetric nonlocal part is compact in the resulting weighted
spaces, then combine two Sobolev energy levels with finite-dimensional
selection. The initial profile is smoothed only near its seam, and the
tuning perturbations have fixed interior support.

\begin{proposition}[Finite-period profile properties]
\label{sb:finite:profile-input}
For a fixed \(a\in(0,1)\), let \(\bar\Omega\) be any finite-period profile
provided by \cite[Theorem~4.4]{Profiles}. Then, for its selected half-period
\(L=\bar L\in(0,\infty)\), the following properties hold on
\(\mathbb T_L=\mathbb R/(2L\mathbb Z)\).
The functions \(\bar\Omega,\bar U\) are odd and periodic,
\[
 \bar\Omega\in C^1(\mathbb T_L)\cap C^\infty((-L,L)),\qquad
 \bar U'=H_L\bar\Omega,\qquad \bar U(0)=\bar U(L)=0,
\]
\[
 \bar\Omega<0,\quad \bar U>0\quad(0<x<L),\qquad
 \bar\Omega(x)=-x+O(x^3)\quad(x\to0).
\]
There are \(b,\kappa_L,A_0>0\) such that, with
\(\bar c_\omega=-(1-a)b\),
\begin{equation}\label{sb:finite:stationary-input}
 b=\bar U'(0),\qquad \kappa_L=-\bar U'(L),\qquad
 a\bar U\bar\Omega'=(\bar c_\omega+H_L\bar\Omega)\bar\Omega,
\end{equation}
and, for \(p_L=(\kappa_L-\bar c_\omega)/(a\kappa_L)>1\),
\begin{equation}\label{sb:finite:seam-input}
 \bar\Omega(L-y)=-A_0y^{p_L}(1+o(1)),\qquad
 |\bar\Omega'(L-y)|\leq Cy^{p_L-1}\quad(y\downarrow0).
\end{equation}
At least one such finite-period profile exists for every
\begin{equation}\label{sb:finite:parameter-range}
 a\in[a_{\mathrm{per}},1),\qquad
 a_{\mathrm{per}}=\frac{64}{176-9\pi^2}.
\end{equation}
\end{proposition}
\begin{proof}
The reconstruction and regularity statements follow from
\cite[Theorems~4.4 and~5.8]{Profiles}. The endpoint estimate in
\cite[Proposition~5.4]{Profiles} gives \eqref{sb:finite:seam-input} and
\(p_L\geq(3-a)/(2a)>1\). Existence for each \(a\in(0,1)\) follows from
\cite[Theorem~4.3]{Profiles}, and \cite[Corollary~5.3]{Profiles} forces
the selected period to be finite when \(a\geq a_{\mathrm{per}}\).
\end{proof}

All the arguments below use only the displayed profile properties.
Set \(v=a\bar U\) and \(Q=\bar c_\omega+H_L\bar\Omega\), so that \(v\bar\Omega'=Q\bar\Omega\).
For the evolving normalized solution, write
\(U(x,\tau)=\int_0^xH_L\Omega(y,\tau)\,\idiff y\) for its velocity.
We use the one-derivative normalization of
Lemma~\ref{sb:common:gauges}:
\begin{equation}\label{sb:finite:normalized}
 \Omega_\tau+aU\Omega_x=(c_\omega+H_L\Omega)\Omega,
 \qquad c_l=0,\qquad c_\omega=-(1-a)H_L\Omega(0).
\end{equation}
It preserves \(\Omega_x(0)=-1\) and the fixed half-period \(L\).
The cubic origin coefficient is free to evolve.
Constants below may depend on the fixed profile, the chosen weight, and
the indicated derivative order.

\begin{theorem}[Finite-period smooth blowup]
\label{sb:finite:theorem}
Let \(\bar\Omega\) satisfy the profile properties in
Proposition~\ref{sb:finite:profile-input}, and fix an integer \(m\geq2\).
There are \(\mu>0\) and the spaces \(X_j\) defined in
\eqref{sb:finite:spaces} such that, for arbitrarily small \(D>0\),
there are smooth odd \(2L\)-periodic data \(\Omega_0\), with
\(\Omega_0'(0)=-1\) and
\(\|\Omega_0-\bar\Omega\|_{X_{m+2}}\leq CD\), whose solutions of
\eqref{sb:finite:normalized} exist smoothly for all \(\tau\geq0\) and
satisfy
\begin{equation}\label{sb:finite:conclusion}
 \|\Omega(\tau)-\bar\Omega\|_{X_m}+|c_\omega(\tau)-\bar c_\omega|
 \leq CD e^{-\mu\tau},\qquad
 \sup_{\tau\geq0}\|\Omega(\tau)-\bar\Omega\|_{X_{m+2}}\leq CD.
\end{equation}
The data are nonpositive on \((0,L)\) and vanish identically in a
neighborhood of the seam. With
\(C_l=L/\pi\) and \(C_\omega(0)=1\), they give smooth odd data
\(\omega_0(y)=\Omega_0((L/\pi)y)\) on the physical circle of period
\(2\pi\). The physical solution has a finite maximal time \(T\) and
\begin{equation}\label{sb:finite:physical-limit}
 (T-t)\omega(y,t)\longrightarrow
 -\frac1{\bar c_\omega} \bar\Omega\left(\frac L\pi y\right)
 \quad\hbox{in }C^1(\mathbb T),\qquad
 (T-t)\|\omega(t)\|_\infty\longrightarrow-\frac{\|\bar\Omega\|_\infty}{\bar c_\omega}.
\end{equation}
Thus every profile in Proposition~\ref{sb:finite:profile-input} can be
realized by smooth blowup data. In particular this gives
constant-spatial-scale blowup for
every parameter in \eqref{sb:finite:parameter-range}.
\end{theorem}

\subsection{Coordinates and profile coefficients}

The weight is chosen to admit the cubic tangent term at the origin
and to control one ordinary derivative at the seam. These requirements
become $\alpha<3$ and $r>1$ in the endpoint powers below.
Fix
\begin{equation}\label{sb:finite:delta}
 \begin{gathered}
 0<\delta<\min\{2ab,a\kappa_L(p_L-1)\},\\
 \beta=\frac\delta{ab},\quad \eta=\frac\delta{a\kappa_L},\quad
 \alpha=1+\beta\in(1,3),\quad r=p_L-\eta\in(1,p_L).
 \end{gathered}
\end{equation}
Choose an additive constant in
\(s_x=1/v\), and set
\begin{equation}\label{sb:finite:spaces}
 z_h(s)=e^{-\delta s}\frac{h(x(s))}{\bar\Omega(x(s))},
 \qquad \|h\|_{X_j}=\|z_h\|_{H^j(\mathbb R_s)}.
\end{equation}
The positive half-cell determines an odd periodic representative.
The definition of the spaces is by completion of smooth functions supported
in \((0,L)\), using their odd periodic extensions. In particular it does
not require high ordinary derivatives at the profile seam.
The solution belongs to the affine space \(\bar\Omega+X_j\); the profile
itself need not belong to \(X_j\). At level zero the weight is
\begin{equation}\label{sb:finite:weight}
 \rho_L(x)=\frac{e^{-2\delta s(x)}}{v(x)\bar\Omega(x)^2}
 \asymp
 \begin{cases}
 x^{-2\alpha-1},&x\downarrow0,\\
 (L-x)^{-2r-1},&x\uparrow L.
 \end{cases}
\end{equation}
We use the auxiliary coordinate
\[
 k_0=\pi/L,\qquad \zeta=\log\tan(k_0x/2),\qquad
 \sin(k_0x)=\operatorname{sech}\zeta,\quad
 \cos(k_0x)=-\tanh\zeta.
\]

\begin{lemma}\label{sb:finite:profile-conormal}
Let \(D_0=\sin(k_0x)\partial_x\) and
\(\mathsf U=\bar U/\sin(k_0x)\). For some \(\sigma\in(0,1)\), and every
integer \(j\geq0\),
\[
 D_0^j\bar\Omega,\quad D_0^j\mathsf U,\quad
 D_0^j\left(\frac{Q}{a\mathsf U}\right)
 \in C^{1,\sigma}(\mathbb T_L),\qquad \inf \mathsf U>0.
\]
Moreover \(|D_0^j\bar\Omega|\leq C_j|\bar\Omega|\).
The change between \(s\) and \(\zeta\) has Jacobian bounded above and
below and bounded derivatives of every order. If
\(M(\zeta)=\bar\Omega(x(\zeta))e^{\delta s(x(\zeta))}\), then
\begin{equation}\label{sb:finite:multiplier-ends}
 |\partial_\zeta^jM|\leq C_j
 \begin{cases}e^{\alpha\zeta},&\zeta<0,\\
 e^{-r\zeta},&\zeta>0.
 \end{cases}
\end{equation}
Consequently, for \(j\geq0\),
\begin{equation}\label{sb:finite:input-embedding}
 \|h(x(\cdot))\|_{H^j_\zeta}
 +\|h(x(\cdot))\|_{L^1_\zeta}
 \leq C_j\|h\|_{X_j},
 \qquad \|h\|_{C^1(\mathbb T_L)}\leq C\|h\|_{X_2}.
\end{equation}
\end{lemma}

\begin{proof}
Equation~\eqref{sb:finite:seam-input} includes
\(|\bar\Omega _x(L-y)|\leq Cy^{p_L-1}\).
Initially \(\bar U'\in C^\theta\) for every \(\theta<1\).
The profile equation writes \(\bar\Omega _x=(\bar\Omega/y)A(y)\)
near the seam, with \(A\in C^\theta\) bounded, since
\(\bar U/y\) has a nonzero limit. Also
\(|\partial_y(\bar\Omega/y)|\leq Cy^{p_L-2}\).
Thus \(\bar\Omega\in C^{1,\sigma}\) for any
\(0<\sigma<\min\{p_L-1,1\}\), after decreasing \(\sigma\) below
\(\theta\). This argument uses a derivative bound, rather than
differentiating an unspecified little-oh remainder.
Hilbert boundedness on this H\"older space gives
\(\bar U\in C^{2,\sigma}\), hence \(\mathsf U\in C^{1,\sigma}\).
Its endpoint values are \(b/k_0\) and \(\kappa_L/k_0\); positivity of
\(\bar U\) on \((0,L)\) gives \(\inf \mathsf U>0\).

For an odd periodic function \(h\), the Fourier sine series gives
\begin{equation}\label{sb:finite:rank-commutator}
 [D_0,H_L]h=\frac{k_0}{2}\widehat h_{\sin,1}.
\end{equation}
Indeed, multiplication by \(\sin(k_0x)\) shifts a Fourier index by
one, so only the first sine mode crosses the zero frequency; direct
evaluation on that mode gives the displayed constant. Set
\(B=Q/(a\mathsf U)\).
The identities
\begin{equation}\label{sb:finite:coefficient-induction}
 D_0\mathsf U=H_L\bar\Omega-k_0\cos(k_0x)\mathsf U,
 \qquad D_0\bar\Omega=B\bar\Omega
\end{equation}
and \eqref{sb:finite:rank-commutator} prove the asserted
\(C^{1,\sigma}\) regularity by simultaneous induction: to obtain the
next derivative of \(B\), differentiate its quotient, commute \(D_0\)
past \(H_L\), and use the positive lower bound for \(\mathsf U\).
Furthermore \(D_0^j\bar\Omega\) is \(\bar\Omega\) times a polynomial
in \(B,D_0B,\ldots,D_0^{j-1}B\), proving its relative bound.

Since \(\partial_\zeta=D_0/k_0\) and
\(\partial_s=ak_0\mathsf U\partial_\zeta\), the coordinate assertion follows.
At the origin smooth parity gives
\(s=\zeta/(ab)+c_-+O(e^{2\zeta})\), with the same differentiated
bounds. At the other end,
\(s=\zeta/(a\kappa_L)+c_++O(e^{-\sigma_1\zeta})\) for some
\(\sigma_1>0\), again conormally. These statements follow by integrating
\(s_\zeta=(ak_0\mathsf U)^{-1}\); the already proved regularity makes each
positive-order \(\zeta\) derivative of \(\mathsf U\) exponentially small at the
two ends. This proves \eqref{sb:finite:multiplier-ends} and the norm
equivalence under the coordinate change. Multiplication by \(M\) proves
the Sobolev estimate in \eqref{sb:finite:input-embedding}; its \(L^1\)
estimate follows from Cauchy--Schwarz, since \(M\in L^2\).
Finally, differentiation of \(h=Mz_h\) with respect to \(x\) has
bounded prefactors tending to zero like \(x^{\alpha-1}\) and
\((L-x)^{r-1}\). Apply \(H^2(\mathbb R)\subset W^{1,\infty}\)
to obtain the last assertion.
\end{proof}

\subsection{Exact kernels and higher order estimates}

For \(h=\Omega-\bar\Omega\), write
\(\partial_xU_h=H_Lh\), \({U_h}(0)=0\), and introduce
\begin{equation}\label{sb:finite:FV}
 \mathcal V=\frac {U_h}{\bar U},\qquad
 F=H_Lh-(1-a)H_Lh(0)-\mathcal VQ.
\end{equation}
The normalized equation, with \(\Omega_x(0)=-1\), is exactly
\begin{equation}\label{sb:finite:z-equation}
 \begin{gathered}
 z_\tau=\mathcal Lz+\mathcal R(z),\qquad
 \mathcal L=-\partial_s-\delta+\mathcal K,\\
 \mathcal Kz=e^{-\delta s}F,
 \quad \mathcal R(z)=-\mathcal Vz_s+(F-\delta \mathcal V)z.
 \end{gathered}
\end{equation}
In particular \(\mathcal K\) is linear. The stationary identity \(v\bar\Omega _x=Q\bar\Omega\) explains the choice of unknown:
writing \(h=\bar\Omega e^{\delta s}z\) gives
\[
 \bar\Omega^{-1}e^{-\delta s}(-v h_x+Qh)=-z_s-\delta z.
\]
Thus the principal transport has constant speed and constant damping
in this coordinate; the spatial scale remains fixed.

\begin{lemma}\label{sb:finite:tame-coefficients}
For every integer \(j\geq1\),
\begin{equation}\label{sb:finite:coefficient-estimate}
 \|\mathcal V\|_\infty+\|\mathcal V_s\|_{H^j_s}
 +\|F\|_\infty+\|F_s\|_{H^{j-1}_s}
 \leq C_j\|z\|_{H^j_s}.
\end{equation}
The functions \(\mathcal V,F\) themselves can have nonzero endpoint limits;
no unweighted \(L^2(\mathrm{d}s)\) assertion about them is made.
For \(n\geq2\), the nonlinear estimates are
\begin{align}
 \|\mathcal R(z)\|_{H^n}
 &\leq C_n\|z\|_{H^n}\|z\|_{H^{n+1}},
 \label{sb:finite:nonlinear-forcing}\\
 |\langle\mathcal R(z),z\rangle_{H^n}|
 &\leq C_n\|z\|_{H^2}\|z\|_{H^n}^2.
 \label{sb:finite:nonlinear-energy}
\end{align}
The second estimate is an energy estimate on smooth functions and extends
by approximation. Its derivative allocation is given explicitly in the
proof.
\end{lemma}

\begin{proof}
Put \(h_\zeta(\zeta)=h(x(\zeta))\), and use a different variable
\(\upsilon\) in the integrals below. Pairing the positive and negative
halves of the cotangent kernel gives the exact formulas
\begin{align}
 (H_Lh)(x(\zeta))
 &=\frac1\pi\operatorname{p.v.}\int_{\mathbb R}
 h_\zeta(\upsilon)
 \{\coth(\zeta-\upsilon)+\tanh\upsilon\}\,\idiff\upsilon,
 \label{sb:finite:H-logtan}\\
 R(\zeta):=\frac{{U_h}(x(\zeta))}{\sin(k_0x(\zeta))}
 &=\frac1{\pi k_0}\int_{\mathbb R}h_\zeta(\upsilon)
 \{A(\zeta-\upsilon)+B(\zeta-\upsilon)\tanh\upsilon\}
 \,\idiff\upsilon,
 \label{sb:finite:u-logtan}
\end{align}
where
\[
 J(w)=\log|\tanh(w/2)|,\qquad A(w)=\cosh w\,J(w),
 \qquad B(w)=\sinh w\,J(w).
\]
For verification, \(\cos(k_0x)=-\tanh\zeta\) converts the paired
Hilbert kernel into
\(\coth(\zeta-\upsilon)+\tanh\upsilon\). The velocity kernel uses
\[
 \left|\frac{\sin(k_0(x-y)/2)}{\sin(k_0(x+y)/2)}\right|
 =|\tanh((\zeta-\upsilon)/2)|,
\]
and \(\mathrm{d}y=k_0^{-1}\operatorname{sech}\upsilon\,\mathrm{d}\upsilon\).

Here are the precise convolution bounds needed from these formulas.
The kernel \(C(w)=\coth w-\operatorname{sgn}w\) is exponentially
decaying away from zero and is \(1/w+O(1)\) near zero. Subtracting an
even smooth cutoff times \(1/w\) leaves an \(L^1\) kernel. Thus its
principal-value convolution is bounded on every \(H^j\), by the
Fourier multiplier bound for the Hilbert kernel and Young's inequality.
For the cutoff kernel, that multiplier bound follows directly from the
uniform boundedness of \(\int_0^R\sin(\xi w)\,\idiff w/w\), followed
by integration against the derivative of the smooth cutoff.
Since the derivative of \(\operatorname{sgn}*h_\zeta\) is
\(2h_\zeta\), \eqref{sb:finite:H-logtan} implies
\begin{equation}\label{sb:finite:H-convolution-bound}
 \|H_Lh\|_\infty
 +\|\partial_\zeta(H_Lh)\|_{H^{j-1}_\zeta}
 \leq C_j(\|h_\zeta\|_{H^j}+\|h_\zeta\|_1).
\end{equation}
For the supremum norm use \(\|\operatorname{sgn}*h_\zeta\|_\infty
\leq\|h_\zeta\|_1\) and \(H^1\subset L^\infty\) for
\(C*h_\zeta\). The last, spatially constant term in
\eqref{sb:finite:H-logtan} is bounded by \(\|h_\zeta\|_1\).

Similarly \(A+1,B+\operatorname{sgn}\in L^2(\mathbb R)\).
The distribution \(A'\) has singular part \(\operatorname{p.v.}(1/w)\)
and an integrable remainder after a smooth local cutoff, while
\(B'\in L^1(\mathbb R)\). Indeed near zero
\(A'=1/w+O(w\log|w|)\), \(B'=O(1+|\log|w||)\); at infinity
both derivatives decay exponentially. Differentiating
\eqref{sb:finite:u-logtan} therefore gives
\begin{equation}\label{sb:finite:u-convolution-bound}
 \|R\|_\infty+\|R_\zeta\|_{H^j}
 \leq C_j(\|h_\zeta\|_{H^j}+\|h_\zeta\|_1).
\end{equation}
Its supremum bound follows from the displayed \(L^2\) kernels and
Cauchy--Schwarz, together with the constant and sign kernels bounded
against \(\|h_\zeta\|_1\).

Now \(\mathcal V=R/\mathsf U\) and
\(F=H_Lh-(1-a)H_Lh(0)-\mathcal VQ\).
Every positive-order \(\zeta\) derivative of \(\mathsf U,Q\) is exponentially
decaying and bounded in all Sobolev orders, by
Lemma~\ref{sb:finite:profile-conormal}. Leibniz' rule,
\(\inf \mathsf U>0\), \eqref{sb:finite:H-convolution-bound}, and
\eqref{sb:finite:u-convolution-bound} prove
\eqref{sb:finite:coefficient-estimate}, first in \(\zeta\) and then in
\(s\). The right side is bounded by \(\|z\|_{H^j_s}\) using
\eqref{sb:finite:input-embedding}.

For completeness, the energy derivative in the transport term is treated
before applying a norm estimate:
\[
 -\int \mathcal V\,\partial_s^{j+1}z\,\partial_s^jz
 =\frac12\int \mathcal V_s|\partial_s^jz|^2.
\]
The other transport terms obey, for \(j\geq1\),
\[
 \|[\partial_s^j,\mathcal V]z_s\|_2
 \leq C_j\bigl(\|\mathcal V_s\|_\infty\|z\|_{H^j}
       +\|\mathcal V_s\|_{H^{j-1}}\|z_s\|_\infty\bigr).
\]
This follows by expanding the commutator; its intermediate products are
bounded by H\"older and the one-dimensional interpolation inequality
between the displayed endpoint norms. For a bounded function \(G\)
with \(G_s\in H^{j-1}\), the same product calculation gives
\[
 \|Gz\|_{H^j}
 \leq C_j\bigl(\|G\|_\infty\|z\|_{H^j}
             +\|G_s\|_{H^{j-1}}\|z\|_\infty\bigr).
\]
Use these inequalities with \(G=F-\delta \mathcal V\), and
\(H^2\subset W^{1,\infty}\), to obtain
\eqref{sb:finite:nonlinear-energy}. Without integrating the first transport
term by parts its bound is
\(C\|z\|_{H^n}\|z\|_{H^{n+1}}\), which proves
\eqref{sb:finite:nonlinear-forcing}.
\end{proof}

\begin{lemma}\label{sb:finite:linear-bounds}
For every integer \(j\geq0\), \(\mathcal K\) is bounded on
\(H^j(\mathbb R_s)\), and its symmetric part in that Hilbert space is
compact. In particular, for \(j\geq1\) and every \(\epsilon>0\),
\begin{equation}\label{sb:finite:compact-energy}
 \Re\langle\mathcal Lz,z\rangle_{H^j}
 \leq-(\delta-\epsilon)\|z\|_{H^j}^2
       +C_{j,\epsilon}\|z\|_{H^{j-1}}^2.
\end{equation}
\end{lemma}

\begin{proof}
We give endpoint estimates explicitly, since ordinary weighted Hilbert
boundedness without the normalization cancellation would be false here.
Before conjugating, the nonlocal operator is
\[
 Nh=\bar\Omega\{H_Lh-(1-a)H_Lh(0)\}-a{U_h}\bar\Omega _x.
\]
Its paired kernels on \((0,L)\) are
\begin{align*}
 K_H(x,y)&=\frac{\cot(\pi(x-y)/(2L))-
                        \cot(\pi(x+y)/(2L))}{2L},\\
 K_U(x,y)&=\frac1\pi\log\left|
              \frac{\sin(\pi(x-y)/(2L))}
                   {\sin(\pi(x+y)/(2L))}\right|,\\
 K_N(x,y)&=\bar\Omega(x)\{K_H(x,y)-(1-a)K_H(0,y)\}
                     -a\bar\Omega _x(x)K_U(x,y).
\end{align*}
At zero the exact cancellation can be displayed as
\begin{equation}\label{sb:finite:origin-cancellation}
 \frac{Nh}{\bar\Omega}
 =H_Lh-H_Lh(0)-a\left(\frac {U_h}x-H_Lh(0)\right)
       -a{U_h}\left(\frac{\bar\Omega _x}{\bar\Omega}-\frac1x\right).
\end{equation}
The last quotient difference is \(O(x)\), with differentiated smooth
bounds. Expansion of the elementary kernels gives, for small \(x,y\),
\begin{equation}\label{sb:finite:origin-kernel}
 |K_N(x,y)|\leq C
 \begin{cases}
 x/y,&2y<x,\\
 x^3(y^{-3}+y^{-1}),&2x<y.
 \end{cases}
\end{equation}
For example \(K_H(x,y)-K_H(0,y)=O(x^2/y^3)\) when \(2x<y\),
and \(K_U(x,y)-xK_H(0,y)=O(x^3/y^3)\); these are precisely the
two canceled expressions in \eqref{sb:finite:origin-cancellation}.
On comparable dyadic intervals of scale \(R\), the weighted operator
norm is \(O(R)\). Rescaling such an interval gives a uniformly bounded
Hilbert transform multiplied by \(\bar\Omega=O(R)\); the logarithmic
and normalization kernels have the same bound by direct integration.

Normalize an input dyadic interval of scale \(Y\) by its \(X_0\) norm;
its typical size is \(Y^\alpha\). After output normalization at scale
\(X\), the off-diagonal matrices in \eqref{sb:finite:origin-kernel} are
bounded by
\begin{equation}\label{sb:finite:origin-matrix}
 C X(Y/X)^\alpha\quad(Y<X/2),\qquad
 C Y(X/Y)^{3-\alpha}\quad(X<Y/2).
\end{equation}
The additional \(y^{-1}\) term has a smaller bound.

At the seam write \(X=L-x\), \(Y=L-y\). The same kernels give
\begin{equation}\label{sb:finite:seam-kernel}
 |K_N(x,y)|\leq C
 \begin{cases}
 X^{p_L}/Y,&2X<Y,\\
 X^{p_L-2}Y+X^{p_L}Y,&2Y<X.
 \end{cases}
\end{equation}
Here we used \(|\bar\Omega|\leq CX^{p_L}\) and
\(|\bar\Omega _x|\leq CX^{p_L-1}\).
Comparable intervals have norm \(O(X^{p_L})\), and the two normalized
off-diagonal matrices are
\begin{equation}\label{sb:finite:seam-matrix}
 C Y^{p_L}(X/Y)^{p_L-r}\quad(X<Y/2),\qquad
 C X^{p_L}(Y/X)^{r+2}\quad(Y<X/2).
\end{equation}
All four ratio exponents are positive. Their geometric row and column
sums prove boundedness by the Schur test. They also prove that an input
or output restricted to a shrinking endpoint region has operator norm
tending to zero. From seam input to origin output the separated kernel
is \(O(x^3(L-y))\), and in the reverse direction it is
\(O((L-x)^{p_L}/y)\); Cauchy--Schwarz with
\eqref{sb:finite:weight} proves the same assertion for these cross terms.

We verify that these estimates survive arbitrary conormal differentiation.
This also supplies the higher-order boundedness not implied by the
zeroth-order Schur test alone. In the \(\zeta\) coordinate the kernel
of \(\mathcal K\), against \(\mathrm{d}\upsilon\), is exactly
\begin{equation}\label{sb:finite:conjugated-kernel}
 \begin{gathered}
 \frac{e^{-\delta s(\zeta)}M(\upsilon)}\pi
 \Bigl(
 \coth(\zeta-\upsilon)+a\tanh\upsilon+(1-a)
 \\
 \hspace{2em}
 -d(\zeta)\{A(\zeta-\upsilon)
               +B(\zeta-\upsilon)\tanh\upsilon\}\Bigr),
 \qquad d=\frac{Q}{k_0\mathsf U}.
 \end{gathered}
\end{equation}
In particular \(d(-\infty)=a\), \(d(+\infty)=-a p_L\), and
\(\partial_\zeta^\ell(d-a)=O(e^{2\zeta})\) at the negative end,
for every fixed \(\ell\). Apply
\(\mathfrak D=\partial_\zeta+\partial_\upsilon\) to this kernel.
It annihilates each function of \(\zeta-\upsilon\), so the only new
factors are derivatives of the displayed coefficients. Their endpoint
powers are unchanged by Lemma~\ref{sb:finite:profile-conormal}.
The leading negative-end bracket cancels identically for every
\(\upsilon\):
\[
 -1+a\tanh\upsilon+(1-a)-a(-1+\tanh\upsilon)=0.
\]
Its derivatives also cancel. Thus the origin third-order bound, as well
as the seam powers, remain valid after each \(\mathfrak D\)
derivative. Equivalently each differentiated kernel has the matrices
\eqref{sb:finite:origin-matrix} and \eqref{sb:finite:seam-matrix}, with a
constant depending on its derivative order. On comparable intervals it
still consists of a Hilbert kernel and a logarithmic kernel with uniformly
bounded smooth coefficients after rescaling. This proves all needed
local bounds without differentiating a Hilbert singularity in isolation.

For compactly supported inputs integration by parts now gives
\[
 \partial_\zeta^j(\mathcal Kz)
 =\sum_{\ell=0}^j\binom j\ell
             \mathcal K_\ell\partial_\zeta^{j-\ell}z,
 \qquad K_\ell=\mathfrak D^\ell K_0.
\]
There are no endpoint boundary terms; extension follows by density.
The established Schur bounds prove boundedness on every \(H^j_\zeta\),
hence every \(H^j_s\). If \(\chi_R\) is a smooth cutoff equal to one
on \([-R,R]\) in \(\zeta\), and zero outside \([-R-1,R+1]\),
the same argument, including the uniformly bounded derivatives of its
transitions, proves
\begin{equation}\label{sb:finite:endpoint-norm-small}
 \|\mathcal K-\chi_R\mathcal K\chi_R\|_{H^j_s\to H^j_s}
 \longrightarrow0.
\end{equation}
For instance every sufficiently small positive rate below the ratio
exponents in \eqref{sb:finite:origin-matrix} and
\eqref{sb:finite:seam-matrix} yields an exponential bound in \(R\).

It remains to identify the symmetric part on a fixed compact interval.
The principal kernel in \(s\) is
\(\bar\Omega(x(s))/(\pi(s-s'))\): changing coordinates in
\eqref{sb:finite:conjugated-kernel} and subtracting this term leaves a
smooth kernel for the Hilbert contribution. The velocity contribution
gains one derivative by \eqref{sb:finite:u-convolution-bound}, after
localizing, and the gauge contribution has finite rank. Consequently
\(\chi_R\mathcal K\chi_R\), modulo a compact operator on \(H^j_s\),
is \(\chi_R\bar\Omega H\chi_R\).
In the \(j\)-th derivative of its quadratic form the terms with a
derivative on a smooth coefficient have one fewer derivative on the input
and are compact by the local Rellich theorem. The remaining symmetric
form is the localized commutator of \(H\) with the real
smooth multiplier \(\bar\Omega\); its kernel has the removable
singularity \((\bar\Omega(x(s))-\bar\Omega(x(s')))/(s-s')\), and is
compact. Summing the derivative levels proves localized symmetric
compactness. Equation~\eqref{sb:finite:endpoint-norm-small} proves the
global assertion.

Finally \(-\partial_s-\delta\) has symmetric part \(-\delta I\)
on every \(H^j\). Approximate the compact symmetric part just proved
in operator norm by a finite-rank form generated by smooth compactly
supported functions. Each of its linear functionals, initially an
\(H^j\) pairing, is bounded in \(H^{j-1}\) by integrating one
derivative onto its smooth test function. This proves
\eqref{sb:finite:compact-energy}.
\end{proof}

\subsection{Closed generators and a common spectral splitting}

\begin{lemma}\label{sb:finite:spectrum}
On each \(H^j_s\), the closed operator \(\mathcal L\) has domain
\(H^{j+1}_s\) and core \(C_c^\infty\). Its essential semigroup growth
bound is at most \(-\delta\). The set of eigenvalues with
\(\Re\lambda\geq0\) is finite, including algebraic multiplicities,
and is independent of \(j\). Its real Riesz projection \(P\) is
consistent on all Sobolev levels; its range \(E\) consists of vectors
in \(\bigcap_jH^j\).
There exist \(\mu>0\), a norm \(|\cdot|_E\) from an inner product on
\(E\), and constants \(M_j,c_E>0\), such that
\begin{align}
 \|e^{\tau\mathcal L}(I-P)\|_{H^j\to H^j}
 &\leq M_j e^{-2\mu\tau},\label{sb:finite:stable-semigroup}\\
 \Re\langle(\mathcal L+\mu)e,e\rangle_E
 &\geq c_E|e|_E^2\qquad(e\in E).
 \label{sb:finite:outgoing-linear}
\end{align}
The empty-space convention applies when \(E=\{0\}\).
\end{lemma}

\begin{proof}
Bounded perturbation of the translation generator gives the stated
domain, core, and semigroup. On a fixed Sobolev level decompose
\(\mathcal K=K_a+K_c\), where \(K_a^*=-K_a\) and \(K_c\) is
compact self-adjoint. Then
\(A=-\partial_s-\delta+K_a\) is a skew-adjoint operator minus
\(\delta I\). Its semigroup has norm \(e^{-\delta\tau}\).
Duhamel's formula shows that the difference between its semigroup and
that of \(A+K_c\) is compact at each positive time. To justify this
claim in operator norm, approximate \(K_c\) by finite-rank operators;
strong continuity of both semigroups and their Hilbert-space adjoints
gives norm continuity of the resulting finite-rank integrands on each
compact time interval. Their integrals are compact, and the approximation
error is bounded by the uniform semigroup bounds on that interval.

For \(\Re\lambda\geq-\delta+\epsilon\),
\((\lambda-A)^{-1}K_c\) is compact and tends to zero in operator norm
as \(|\Im\lambda|\to\infty\), uniformly in this half-plane. This
follows from the spectral theorem for the skew-adjoint part of \(A\):
the resolvent tends strongly to zero on fixed vectors, is bounded by
\(1/\epsilon\), and the range of \(K_c\) is precompact on bounded
sets. The norm also tends to zero when \(\Re\lambda\to+\infty\).
Analytic Fredholm theory consequently leaves only finitely many poles in
each such half-plane: use the larger open half-plane
\(\Re\lambda>-\delta+\epsilon/2\) to exclude accumulation at
\(\Re\lambda=-\delta+\epsilon\). Together with the compact semigroup difference,
the spectral decomposition for a semigroup with negative essential
growth bound gives exponential decay on the complementary stable
subspace \cite[Chapter~IV, Corollary~2.11 and Propositions~1.16, 2.12]{EngelNagel2000}.

These spectral statements are compatible across levels. An eigenvector
in \(H^j\) satisfies
\(z_s=-(\lambda+\delta)z+\mathcal Kz\), hence belongs to
\(H^{j+1}\); induction gives every higher order. The same argument
applies successively to a generalized eigenvector chain. The resolvents
coincide on common vectors for sufficiently large positive \(\lambda\)
by uniqueness, then by analytic continuation in the half-plane
\(\Re\lambda>-\delta\) away from its isolated poles. Thus their
Riesz projections coincide. Choose \(\mu>0\) so small that the stable
spectrum and the essential growth bound lie strictly to the left of
\(-2\mu\); enlarging the finitely many constants gives
\eqref{sb:finite:stable-semigroup} at every level used below.

On \(E\), all eigenvalues of \(\mathcal L+\mu\) have positive real
part. In any initial Euclidean coordinates the positive matrix
\[
 G=\int_0^\infty
 e^{-t(\mathcal L+\mu)^*}e^{-t(\mathcal L+\mu)}\,\idiff t
\]
defines an inner product and satisfies
\((\mathcal L+\mu)^*G+G(\mathcal L+\mu)=I\).
Finite-dimensional norm equivalence gives
\eqref{sb:finite:outgoing-linear}.
\end{proof}

\subsection{The two Sobolev levels and the shrinking region}

Fix \(m\geq2\), put \(n=m+2\), and let
\(M\geq\max\{1,M_m\}\) in
\eqref{sb:finite:stable-semigroup}. The following estimates close the
derivative loss in \eqref{sb:finite:nonlinear-forcing} by using a bounded
high norm and a decaying lower norm, and integrating the highest transport derivative by parts.

\begin{lemma}\label{sb:finite:trapping}
There are \(K_0\geq1\) and, for each \(K\geq K_0\), a
\(D_K>0\), with the following property. Suppose \(D\leq D_K\),
\[
 \|(I-P)z(0)\|_{H^m}\leq D,\qquad
 \|z(0)\|_{H^n}\leq KD/2.
\]
On any interval on which
\begin{equation}\label{sb:finite:tube}
 \|(I-P)z(\tau)\|_{H^m}\leq2MD e^{-\mu\tau},\qquad
 |Pz(\tau)|_E\leq D e^{-\mu\tau},\qquad
 \|z(\tau)\|_{H^n}\leq2KD,
\end{equation}
the first and third bounds improve to
\begin{equation}\label{sb:finite:tube-improvement}
 \|(I-P)z(\tau)\|_{H^m}\leq\tfrac32 MD e^{-\mu\tau},
 \qquad \|z(\tau)\|_{H^n}\leq KD.
\end{equation}
At a contact with the middle boundary in
\eqref{sb:finite:tube},
\begin{equation}\label{sb:finite:transverse}
 \frac{\mathrm{d}}{\mathrm{d}\tau}\{e^{2\mu\tau}|Pz(\tau)|_E^2\}
 \geq c_E D^2>0.
\end{equation}
The constants in the lower-norm bound are independent of the later
choice of \(K\).
\end{lemma}

\begin{proof}
Finite-dimensional norm equivalence and the tube give
\(\|z(\tau)\|_{H^m}\leq C_0D e^{-\mu\tau}\), where
\(C_0=2M+\sup_{|e|_E=1}\|e\|_{H^m}\) is independent of \(K\).
Interpolation and \eqref{sb:finite:nonlinear-forcing} imply
\begin{equation}\label{sb:finite:decaying-forcing}
 \|\mathcal R(z(\tau))\|_{H^m}
 \leq C\|z\|_{H^m}^{3/2}\|z\|_{H^n}^{1/2}
 \leq C\sqrt K D^2 e^{-3\mu\tau/2}.
\end{equation}
The stable Duhamel formula is valid in \(H^m\), since the actual
solution has the bounded \(H^n\) regularity in the tube. Consequently
\[
 \begin{aligned}
 \|(I-P)z(\tau)\|_{H^m}
 &\leq MD e^{-2\mu\tau}
   +C\sqrt K D^2\int_0^\tau
       e^{-2\mu(\tau-t)}e^{-3\mu t/2}\,\idiff t
 \\
 &\leq MD e^{-2\mu\tau}
       +\frac{2C\sqrt K}{\mu}D^2 e^{-3\mu\tau/2}.
 \end{aligned}
\]
Taking \(D_K\) small proves the first strict improvement.

For the high norm, use \eqref{sb:finite:compact-energy} with
\(\epsilon=\delta/4\), and
\eqref{sb:finite:nonlinear-energy}. Young's inequality and interpolation
\(\|z\|_{H^{n-1}}^2\leq C\|z\|_{H^m}\|z\|_{H^n}\)
give, after reducing \(D_K\),
\begin{equation}\label{sb:finite:high-differential}
 \frac{\mathrm{d}}{\mathrm{d}\tau}\|z\|_{H^n}^2
 \leq-\delta\|z\|_{H^n}^2+C_1D^2e^{-2\mu\tau}.
\end{equation}
Here is the absorption explicitly: the linear compact term is at most
\(\delta\|z\|_{H^n}^2/8+C\|z\|_{H^m}^2\), and the nonlinear
coefficient is at most \(CC_0D\), chosen smaller than \(\delta/8\),
in the half-energy identity. The resulting negative coefficient there
is at least \(\delta/2\). Thus \(C_1\) is independent of \(K\).
Integration gives
\[
 \|z(\tau)\|_{H^n}^2
 \leq\frac14K^2D^2e^{-\delta\tau}+\frac{C_1}{\delta}D^2.
\]
Choose \(K_0^2\geq4C_1/\delta\). This is at most
\(K^2D^2/2\), proving the second strict improvement.

Finally let \(e=Pz\). Its finite-dimensional equation and
\eqref{sb:finite:outgoing-linear} imply
\[
 \frac{\mathrm{d}}{\mathrm{d}\tau}(e^{2\mu\tau}|e|_E^2)
 \geq2e^{2\mu\tau}
       \{c_E|e|_E^2-C\|\mathcal R(z)\|_{H^m}|e|_E\}.
\]
At \(|e|_E=D e^{-\mu\tau}\), the ratio of the second term to
\(|e|_E^2\) is at most \(C\sqrt K D e^{-\mu\tau/2}\).
Decrease \(D_K\) once more to prove \eqref{sb:finite:transverse}.
This chooses \(K\) before \(D_K\), with no circular dependence.
\end{proof}

\subsection{Smooth data, continuity, and the finite-dimensional choice}

\begin{proof}[Completion of the proof of Theorem~\ref{sb:finite:theorem}]
We record the smooth local facts used in the exit argument. For integer
\(q\geq3\), differentiating the normalized periodic equation yields
\begin{equation}\label{sb:finite:ordinary-continuation}
 \frac{\mathrm{d}}{\mathrm{d}\tau}\|\Omega\|_{H^q}
 \leq C_q\bigl(1+|c_\omega|+\|H_L\Omega\|_\infty
                      +\|\Omega\|_{C^1}\bigr)\|\Omega\|_{H^q}.
\end{equation}
Indeed the transport commutator is bounded by
\(C_q(\|U_x\|_\infty\|\Omega\|_{H^q}
+\|U\|_{H^q}\|\Omega_x\|_\infty)\), the stretching product by
\(C_q(\|\Omega\|_\infty+\|H_L\Omega\|_\infty)
\|\Omega\|_{H^q}\), and
\(\|U\|_{H^q}\leq C\|\Omega\|_{H^{q-1}}\).
Fourier truncation and this inequality give smooth local solutions:
at a fixed high order the truncated equations have a common existence
interval and uniform bounds, their differences converge one order lower
by the corresponding difference energy inequality, and interpolation
with the higher bounds gives strong convergence at each prescribed
lower order. Running this construction with arbitrarily high orders
gives a smooth solution and its continuity in the initial smooth data.
The interval does not shrink when the derivative order is raised:
construct first at order three, which bounds the coefficient in
\eqref{sb:finite:ordinary-continuation}, and apply that estimate at
each higher order on the same interval.
The same difference inequality gives uniqueness. Restarting this local
construction and using \eqref{sb:finite:ordinary-continuation} proves
continuation across any finite time on which the displayed coefficient
remains bounded.

For the solutions in \eqref{sb:finite:tube}, this coefficient is uniformly
bounded: \eqref{sb:finite:input-embedding} bounds \(\Omega\) in \(C^1\),
Hilbert boundedness from \(C^{1/2}\) to itself bounds \(H_L\Omega\),
and the gauge bounds \(c_\omega\). Thus no finite-time failure of smooth
existence occurs before an exit.
Parity is preserved, and Lemma~\ref{sb:common:gauges} preserves the origin
slope. Along characteristics \(X_\tau=aU(X,\tau)\),
\[
 \Omega(X(\tau),\tau)=\Omega_0(X(0))
 \exp\left(\int_0^\tau(c_\omega+U_x)(X(t),t)\,\idiff t\right).
\]
Hence a seam neighborhood on which the datum vanishes remains a
vanishing neighborhood at every finite time; its width on a compact time
interval is bounded below by its initial width times
\(\exp(-a\int\|U_x\|_\infty)\).

These facts also give the particular weighted continuity needed below.
On a fixed compact time interval and for a compact smooth parameter
family, the seam gap is uniformly positive and the ordinary smooth norms
are uniformly bounded. Differences of two normalized solutions have zero
first derivative at zero and are odd. Taylor's theorem therefore bounds
every fixed conormal derivative of their difference by \(Cx^3\), with
\(C\) tending to zero as the ordinary \(C^{n+3}\) difference tends to
zero. This bound is integrable in \eqref{sb:finite:weight}, since
\(\alpha<3\); the seam gap handles the other end. Thus their
\(X_n\) distance tends to zero as well. The same reasoning applies to
time differences. It proves joint weighted continuity, including on a
short common continuation beyond each finite contact with the tube.

Choose a smooth cutoff \(0\leq\vartheta_\epsilon\leq1\) equal to zero when
\(L-|x|\leq\epsilon\), equal to one when
\(L-|x|\geq2\epsilon\), and even and periodic. The function
\(\Omega_\epsilon=\vartheta_\epsilon \bar\Omega\) is smooth on the
circle, since the profile is smooth on every compact subset of
\((-L,L)\) and the cutoff vanishes on a seam neighborhood.
It agrees with the profile near the origin. For
\(z_\epsilon=z_{\Omega_\epsilon-\bar\Omega}\),
Lemma~\ref{sb:finite:profile-conormal} gives
\begin{equation}\label{sb:finite:cutoff-smallness}
 \|z_\epsilon\|_{H^j_s}\leq C_j\epsilon^\eta
 \qquad\text{for each fixed }j.
\end{equation}
Indeed \(z_\epsilon=e^{-\delta s}(\vartheta_\epsilon-1)\);
the transition has bounded derivatives in \(s\), its location satisfies
\(s\geq(a\kappa_L)^{-1}\log(1/\epsilon)-C\), and the integral of
\(e^{-2\delta s}\) on this half-line is \(O(\epsilon^{2\eta})\).

Let \(P,E\) be as in Lemma~\ref{sb:finite:spectrum}. Density of
\(C_c^\infty(\mathbb R_s)\) and finite dimensionality give a linear map
\(J:E\to C_c^\infty(\mathbb R_s)\) such that
\begin{equation}\label{sb:finite:bump-right-inverse}
 PJ=I_E,\qquad
 \|(I-P)Je\|_{H^m}\leq\varepsilon_J|e|_E,
 \qquad \|Je\|_{H^n}\leq C_J|e|_E.
\end{equation}
Here \(\varepsilon_J>0\) can be prescribed first. To construct \(J\),
approximate a basis of \(E\subset H^n\) by compactly supported smooth
functions in \(H^n\), then invert their projected finite matrix.
This proves both the exact right-inverse property and the small stable
component. Their corresponding physical perturbations
\(\bar\Omega e^{\delta s}Je\) are smooth functions supported in a
fixed compact subset of \((0,L)\); extend them oddly.

For \(\xi\in E\), \(|\xi|_E\leq D\), set
\begin{equation}\label{sb:finite:shooting-family}
 z_0(\xi)=z_\epsilon+J(\xi-Pz_\epsilon),\qquad
 \Omega_0(\xi)=\Omega_\epsilon+
 \bar\Omega e^{\delta s}J(\xi-Pz_\epsilon).
\end{equation}
Then \(Pz_0(\xi)=\xi\) exactly. First fix
\(\varepsilon_J<1/4\), then choose a sufficiently large constant
\(C_*\) and set \(D=C_*\epsilon^\eta\).
Equations~\eqref{sb:finite:cutoff-smallness}--
\eqref{sb:finite:bump-right-inverse} give, uniformly on this parameter ball,
\[
 \|(I-P)z_0(\xi)\|_{H^m}\leq D,\qquad
 \|z_0(\xi)\|_{H^n}\leq C_{\rm init}D.
\]
Choose \(K\geq\max\{K_0,2C_{\rm init}\}\), and only then take
\(\epsilon\) small enough that \(D\leq D_K\).
All constants and all smooth correction functions are now fixed before
the last smallness choice. The family consists of smooth odd data,
has a common seam gap, and preserves the origin slope.
It can also be made nonpositive on \((0,L)\). Indeed the correction
bumps have fixed interior support, where \(\vartheta_\epsilon=1\)
for small \(\epsilon\), and their ratio to \(\bar\Omega\) is uniformly
\(O(D)\). On that support \(\Omega_0=\bar\Omega(1+O(D))<0\), while
outside it \(\Omega_0=\vartheta_\epsilon \bar\Omega\leq0\).

Suppose \(E\neq\{0\}\). For each member of
\eqref{sb:finite:shooting-family}, let \(\tau_*(\xi)\) be the first
contact with a boundary of \eqref{sb:finite:tube}. The strict improvements
in Lemma~\ref{sb:finite:trapping} leave only the boundary
\(|Pz|_E=D e^{-\mu\tau}\), and
\eqref{sb:finite:transverse} makes every such contact outward. If every
\(\tau_*(\xi)\) were finite, the map
\begin{equation}\label{sb:finite:exit-map}
 \xi\longmapsto e^{\mu\tau_*(\xi)}Pz(\tau_*(\xi);\xi)
\end{equation}
would be continuous from the closed radius-\(D\) ball to its sphere.
Here is the required continuity argument. Before a first contact the
continuous boundary functions are strictly negative on every compact
subinterval; the weighted continuity just proved preserves this property
for nearby parameters. At the contact the scalar function
\(e^{2\mu\tau}|Pz|_E^2-D^2\) has strictly positive derivative,
continuous in parameter and time by the finite-dimensional projected
equation and \eqref{sb:finite:nonlinear-forcing}. It therefore changes sign
in any sufficiently short neighborhood of the contact, uniformly for
nearby parameters. These two observations prove continuity of the first
contact time and of \eqref{sb:finite:exit-map}, including contacts at time
zero. A short common continuation exists by
\eqref{sb:finite:ordinary-continuation}. For \(|\xi|_E=D\), the contact
time is zero and \eqref{sb:finite:exit-map} equals \(\xi\).
This would retract a closed Euclidean ball onto its boundary, which is
impossible. Hence at least one parameter never exits. If \(E=\{0\}\),
there is no possible exit boundary, so the single member never exits.

The smooth continuation argument gives this solution for every finite
\(\tau\). The tube yields \eqref{sb:finite:conclusion}, since
\(H_Lh(0)\) is bounded by \(\|z_h\|_{H^m}\) through
\eqref{sb:finite:H-convolution-bound}. The embedding in
\eqref{sb:finite:input-embedding} gives uniform, indeed \(C^1\),
convergence on the fixed circle. Also
\(c_\omega-\bar c_\omega\) is exponentially integrable. With
\(C_l=L/\pi\) and \(C_\omega(0)=1\), physical reconstruction therefore
gives a finite \(T\),
\((T-t)C_\omega\to-1/\bar c_\omega\), and the uniform profile and type-I
limits of Lemma~\ref{sb:common:physical}. The same reconstruction and
the \(C^1\) convergence just proved give \eqref{sb:finite:physical-limit}.
The original solution is
smooth for every \(t<T\) by the continuation already established;
its nonzero limiting profile excludes smooth continuation through \(T\).
Finally \(D\to0\) with \(\epsilon\to0\), so the initial distance in
\(X_{m+2}\) is arbitrarily small, as asserted.
\end{proof}

\vspace{0.2in}
\noindent{\bf Acknowledgments.}
This work was supported in part by NSF Grant DMS-2508463,
the Choi Family Gift Fund, and the Mike Yan Gift Fund.

\vspace{0.2in}
\noindent{\bf Disclosure of AI assistance.} The code for the computer-assisted analysis of the critical profile was developed with the assistance of OpenAI’s Codex. The authors directed and verified this work and take full responsibility for its correctness.

\appendix
\section{A coercive critical profile from scalar estimates}
\label{sb:sc:introduction}

We construct one critical whole-line profile together with the energy
estimate needed to produce smooth periodic blowup. The proof uses one
rational reference, one-dimensional scalar bounds, and an exact
four-by-four inertia calculation.

Our convention is
\[
 Hv(x)=\frac1\pi\operatorname{p.v.}\int_{\mathbb R}
                      \frac{v(y)}{x-y}\,\mathrm{d}y,\qquad
 U_W(x)=\int_0^xHW(y)\,\mathrm{d}y .
\]
In this appendix $W=-xf$ denotes a variable vorticity profile and
$U_W$ its velocity. At the exact critical profile, $W_*=\bar\Omega$
and $U_*=\bar U$ in the notation of Section~\ref{sb:paper:critical}.
For even $f$ with $f(0)=1$ set $P_f=U_W/x$,
$b_f=HW(0)$ and $c_f=-\tfrac12(HW)''(0)$. Equivalently,
\[
 b_f=\frac2\pi\int_0^\infty f(y)\,\idiff y,\qquad
 c_f=\frac2\pi\int_0^\infty\frac{1-f(y)}{y^2}\,\mathrm{d}y,
 \qquad a_f=\frac{3c_f}{6b_f+c_f}.
\]
These moment formulas define the affine extensions on the weighted
space below. The positive potential kernel is
\[
 P_f(x)=\int_0^\infty \frac{y}{\pi x}
             \log\left|\frac{x+y}{x-y}\right|f(y)\,\mathrm{d}y;
\]
it follows by integrating the paired Hilbert kernel. Let $f_0$ denote
the exact rational function in \eqref{sb:sc:reference}, and put
\begin{equation}\label{sb:sc:space}
 W_0=-xf_0,\quad P_0=P_{f_0},\quad a_0=a_{f_0},\quad
 \rho_0=(x^6f_0)^{-1},\quad
 \|v\|_X^2=\int_{\mathbb R}\frac{|v|^2}{x^4f_0}\,\idiff x.
\end{equation}
The map $v\mapsto -xv$ is an isometry from the even space $X$
onto the odd space $L^2(\rho_0)$. We use it when passing between
profile and vorticity variables.

For a profile $W=-xf$ and its parameter $a$, write $b=HW(0)$,
$U=U_W$, and
\[
 B(h)=-\frac{3-a}{12}(Hh)''(0)-aHh(0),\quad
 G(h)=B(h)+(a-1)Hh(0),\quad U_h(x)=\int_0^xHh(y)\,\idiff y.
\]
The normalized fixed-parameter linearization is
\begin{equation}\label{sb:sc:linear}
 Lh=-aUh'+(U'-(1-a)b)h+(Hh+G(h))W-(aU_h+B(h)x)W'.
\end{equation}
Subscript zero means evaluation at the reference. For constructing the
profile, the parameter also varies:
\[
 \mathcal F_*(W)=\bigl(HW-(1-a_f)b_f\bigr)W-a_fU_WW',\qquad
 M_0=D\mathcal F_*(W_0).
\]
These derivatives are closed transport operators with bounded parts,
as justified below; no smoothness of the nonlinear differential
vector field on the rough $X$ ball is being assumed.
Differentiation of the moments gives
\[
 Da(W_0)[h]=\frac{B_0(h)}{b_0+c_0/6},\qquad
 M_0=L_0+e_0B_0,\quad
 e_0=\frac{b_0W_0-U_0W_0'}{b_0+c_0/6}-(W_0-xW_0').
\]
Here $a(W)=a_{-W/x}$ and the derivative is in the odd variable;
in the even variable substitute $h=-xv$.

\begin{theorem}\label{sb:sc:profile}
There are $a_{\mathrm c}\in(689/1000,690/1000)$ and an even positive
$f_*\in\mathcal S(\mathbb R)$ with $f_*(0)=1$, $f_*''(0)=-2$,
such that $W_*=-xf_*$ and $U_*'=HW_*$, $U_*(0)=0$, satisfy
\[
 a_{\mathrm c}U_*W_*'=(U_*'-(1-a_{\mathrm c})b_*)W_*,\quad
 P_*=U_*/x>0,\quad b_*=HW_*(0)>0.
\]
The fixed-parameter operator $L_*$ in \eqref{sb:sc:linear} obeys
\begin{equation}\label{sb:sc:profile-coercivity}
 \operatorname{Re}\int_{\mathbb R}\rho_0\bar h L_*h\,\idiff x
       \leq-\frac1{20}\int_{\mathbb R}\rho_0|h|^2\,\idiff x
\end{equation}
for smooth odd compactly supported $h$ with $h'(0)=h'''(0)=0$.
Moreover $\|f_*-f_0\|_X<10^{-15}$.
\end{theorem}

The proof is given in the following sections. It constructs the parameter
simultaneously with the profile and concerns this single local fixed point.
The certificate is part of the proof; the result is computer assisted.
The proof has the following order. The exact rational reference gives
verified residual bounds and separate coercivity gaps for \(M_0\) and
\(L_0\). These bounds yield the analytic inverse and Newton construction
of the exact profile. A form comparison then transfers fixed-parameter
coercivity to that profile. Proposition~\ref{sb:crit:spectral-gap} converts
this coercivity into semigroup decay on the dynamical spaces, and
Theorem~\ref{sb:crit:transfer} gives periodic blowup. The reference
bounds use the finite certificate; the inverse, existence, regularity,
and transfer steps are analytic arguments with explicit constants.

\subsection{The finite scalar certificate}\label{sb:sc:subsection}\label{sb:sc:certificate}

The reference is an exact rational function, specified by
\begin{equation}\label{sb:sc:reference}
 \begin{gathered}
 z=(1+x^2/25)^{-1},\qquad u=1-z,\qquad f_0(x)=z^5 A(u),\\
 A(u)=1-20u+u^2Q(1-u),\qquad Q(z)=\sum_{j=0}^{96}q_jz^j.
 \end{gathered}
\end{equation}
The file \texttt{reference.json} in the repository described in
Appendix~\ref{sb:repro:appendix} gives the rational numbers \(q_j\)
as its \texttt{q\_coefficients} array, ordered by \(j\).
This finite list is part of the definition, rather than rounded data.
Formula \eqref{sb:sc:reference} gives \(f_0(0)=1\) and
\(f_0''(0)=-2\) identically. The verification below proves
\(f_0>0\), \(P_0>0\), and \(f_0(x)\sim c_\infty x^{-10}\)
with \(c_\infty>0\). The stationary defect of this smooth reference is
defined in \eqref{sb:sc:defect} and retained in the energy calculation.

All Hilbert transforms used to check the reference are polynomial
calculations. For \(S=5\), \(z=(1+x^2/S^2)^{-1}\), and
\(\beta_p=4^{1-p}\binom{2p-2}{p-1}\), one has
\begin{equation}\label{sb:sc:hilbert}
 H(z^p)=\frac{x}{S}\sum_{j=1}^p\beta_{p-j+1}z^j,
 \qquad b_{z^p}=S\beta_p,\qquad
 c_{z^p}=\frac{2p-1}{S}\beta_p.
\end{equation}
These identities follow by residues, or induction from
\(H(z)=xz/S\). The potential is recovered by the polynomial identity
\(P_0-2z(1-z)\partial_zP_0=H(-xf_0)\), with \(P_0(0)=0\)
in the \(z\) coordinate. Thus \(P_0\), its derivatives, and
\(a_0=3c_0/(6b_0+c_0)\)
are computed exactly before any interval evaluation.

Here is the energy criterion used in the verification. Write
\(W_0=-xf_0\), \(j=-xf_0'/f_0\), \(\alpha=(3-a_0)/6\), and
\begin{align}
 r_{\rm ref}&=U_0'-(1-a_0)b_0-a_0U_0W_0'/W_0,
                                                        \label{sb:sc:defect}\\
 D&=(1-a_0)b_0+\frac{5a_0-2}{2}P_0
          -\frac{2+a_0}{2}xP_0'-\frac{a_0}{2}P_0j,
                                                        \label{sb:sc:damping}\\
 E&=(1-a_0/3)(1-P_0/b_0)
       -\bigl(1-(1-a_0/3)P_0/b_0\bigr)j.             \label{sb:sc:E}
\end{align}
The damping \(D\) includes the nonzero reference defect:
it equals the stationary expression
\((1-a_0)b_0/2+(4a_0-1)P_0/2-(1+a_0)xP_0'/2\), minus \(r_{\rm ref}/2\).
For a real odd smooth test function \(h\) with vanishing first and third derivatives at the origin, put \(g=h/x^2\), \(J_3v=\int_0^1t^2v(tx)\,\mathrm{d}t\), and
\[
 \mathbf t=\left(\int_{\mathbb R}\frac g x\,\idiff x,
       \int_{\mathbb R}\frac{jg}{x^3}\,\idiff x,
       \int_{\mathbb R}xg\,\idiff x,
       \int_{\mathbb R}\frac{Eg}{x^3}\,\idiff x\right)^{\!T},\qquad
 C=\begin{pmatrix}
 1&-\alpha&0&-\alpha\\ -\alpha&0&0&0\\
 0&0&0&-a_0\\ -\alpha&0&-a_0&0
 \end{pmatrix}.
\]
Polynomial subtraction gives
$Hh=Hh(0)+x^2Hg$, $U_h=xHh(0)+x^3J_3Hg$, and
$c_1:=\tfrac12(Hh)''(0)=Hg(0)$. The elementary Hilbert identity
$H(gHg)=((Hg)^2-g^2)/2$ follows by expanding the positive and
negative Fourier projections. Evaluation at zero yields
$-\int_{\mathbb R}(gHg/x)\,\idiff x=(\pi/2)c_1^2$. Together with $M_0=L_0+e_0B_0$
and integration by parts, this proves the exact identity
\begin{equation}\label{sb:sc:energy}
 \langle M_0h,h\rangle_{\rho_0}
 =-\int_{\mathbb R}D\rho_0h^2\,\idiff x
   +a_0\int_{\mathbb R}g\frac{W_0'}{W_0}J_3Hg\,\idiff x
   +\frac1{2\pi}\mathbf t^TC\mathbf t.
\end{equation}
For \(L_0\), only the upper left two-by-two block of \(C\) remains.
The identities initially hold for compactly supported tests and extend
on the closed operator domains established below.

Define \(m_0=\sup x^2f_0/D\),
\(m_1=\sup f_0(1-j)^2/D\), and the Gram matrix
\begin{equation}\label{sb:sc:gram}
 G=\int_{\mathbb R}\mathbf v\mathbf v^T\,\idiff x,\qquad
 \mathbf v=\sqrt{f_0/D}\,(1,j/x^2,x^2,E/x^2)^T.
\end{equation}
The Hardy bound \(\|J_3H\|_{L^2\to L^2}\leq2/5\) shows that the
primitive term is at most
\[
 \lambda_p\int_{\mathbb R}D\rho_0h^2\,\idiff x,\qquad
 \lambda_p=(2a_0/5)\sqrt{m_0m_1}.
\]
The matrix \(C\) has inertia \((2,2)\). If
\(2\pi C^{-1}-G/\lambda_r\) also has inertia \((2,2)\),
the rank term in \eqref{sb:sc:energy} is at most
\(\lambda_r\int_{\mathbb R}D\rho_0h^2\,\idiff x\), by the Schur complement.
Thus \(D\geq d>0\) and \(\lambda_p+\lambda_r<1\) imply an
energy gap \(d(1-\lambda_p-\lambda_r)\).

We also bound the bounded part in
\(M_0=-a_0U_0\partial_x+(U_0'-(1-a_0)b_0)+R_M\).
Let \(G^{\rm raw}\) denote \eqref{sb:sc:gram} with \(D\) replaced
by \(1\), and let \(m_i^{\rm raw}\) have the same meaning. Using the certified bound $f_0\leq1$, we have
\begin{align}\label{sb:sc:raw}
 \|R_M\|\leq{}&\sqrt{m_0^{\rm raw}}
  +\frac{2a_0}{5}\sqrt{m_0^{\rm raw}m_1^{\rm raw}}
  +\frac\alpha\pi\sqrt{G^{\rm raw}_{00}G^{\rm raw}_{11}}\notag\\
 &+\frac1\pi\sqrt{G^{\rm raw}_{33}
   (\alpha^2G^{\rm raw}_{00}+2\alpha a_0G^{\rm raw}_{02}
                                    +a_0^2G^{\rm raw}_{22})}.
\end{align}
The last summand is omitted for the bounded part of \(L_0\).
This follows directly from Cauchy--Schwarz and the same Hardy bound.

For completeness, the finite arithmetic covers the entire half-line.
Use \(x=5\tan\theta\), \(0\leq\theta\leq\pi/2\), so
\(u=\sin^2\theta\). Write the exact polynomial identities
\(P_0=(1-u)V\), \(j=uJ/A\), \(E=uB/A\),
\(r_{\rm ref}=u^2R/A\), and \(D=D_* /A\). Then
\begin{equation}\label{sb:sc:integrals}
\begin{gathered}
 \mathbf Z=((1-u)A,(1-u)^2J/25,25uA,(1-u)^2B/25)^T,\\
 G=10\int_0^{\pi/2}\frac{(1-u)^2}{D_*}\mathbf Z\mathbf Z^T\,\mathrm{d}\theta,\\
 \|W_0r_{\rm ref}\|_{\rho_0}^2
 =\frac2{125}\int_0^{\pi/2}
          \frac{u^2(1-u)^6R^2}{A}\,\mathrm{d}\theta.
\end{gathered}
\end{equation}
Replacing the power \(6\) by \(3\) gives the squared residual norm
in \(\rho_0K\), where \(K=(1+x^2/25)^3\).

The program uses 4096 full angle cells and 96 exact dyadic polynomial
charts. On each chart, positive center and radius make the affine
coefficient transformation monotone. Separate 180-digit Decimal
rounding downward and upward therefore encloses each transformed
coefficient. Terms above degree 32 are bounded by their absolute
coefficient sum on the chart interval \([-1,1]\), also for the first
two derivatives. A second-order interval Taylor formula encloses each
whole cell. Outward cell integration, including both endpoints,
then bounds \eqref{sb:sc:integrals}; interval LDL factorization proves
the required four-by-four inertia.

\begin{table}[ht]
\centering
\caption{Certified bounds for the rational reference. Decimal endpoints are
rounded outward, and all weighted norms are taken over \(\mathbb R\).}
\label{sb:sc:table}
\small
\begin{tabular}{@{}lll@{}}
\toprule
Estimate & Quantity & Certified bound \\
\midrule
Damping & \(D,\ D-3a_0P_0u\) & \(\geq 0.09\) \\
Scalar bounds & \(m_0\) & \(<0.636\) \\
 & \(m_1\) & \(<1.604\) \\
Operator bounds & \(\lambda_r\) & \(\leq 0.6\) \\
 & \(\lambda_p\) & \(<0.278294\) \\
 & \(\|R_M\|\) & \(<1.542\) \\
 & \(\|R_L\|\) & \(<1.252\) \\
Coercivity & Gap of \(M_0\) & \(>0.01095\) \\
 & Gap of \(L_0\) & \(>0.06302\) \\
Residuals & \(\|W_0r_{\rm ref}\|_{\rho_0}\) & \(<1.721\times 10^{-18}\) \\
 & \(\|W_0r_{\rm ref}\|_{\rho_0K}\) & \(<7.885\times 10^{-17}\) \\
Reference parameters & \(a_0\) & \(\in(0.689,0.690)\) \\
 & \(b_0\) & \(\in(0.613,0.614)\) \\
 & \(c_0\) & \(\in(1.097,1.098)\) \\
\bottomrule
\end{tabular}
\end{table}
The same check proves
\begin{gather*}
 0<f_0\leq\min\{1,4(1+x^2/4)^{-5}\},\quad
 f_0\geq e^{-x^2/2}/10,\quad
 \frac{0.4}{1+x^2}\leq P_0\leq\frac{0.75}{1+x^2},\\
 |xP_0'|\leq\frac{2x^2}{1+x^2},\quad
 b_0-P_0\geq\frac{x^2}{5(1+x^2)},\quad
 |P_0-b_0+c_0x^2/3|\leq x^4\quad(|x|\leq1/4),\\
 x^2f_0<0.418,\qquad x^2P_0<0.44\quad(|x|\geq6).
\end{gather*}
For the Gaussian lower bound, a binomial comparison is checked on
\([0,3]\), and the polynomial inequality
\(25A-(1-u)J\geq0\) proves that \(\log f_0+x^2/2\) increases
on \([3,\infty)\).
The directory \texttt{critical\_scalar} in the verification repository
(see Appendix~\ref{sb:repro:appendix}) contains
\texttt{reference.json}, \texttt{scalar.json}, \texttt{verify.py},
\texttt{interval.py}, \texttt{check.py}, and \texttt{manifest.json}.
The script \texttt{check.py} checks the supplied inequalities and the
consistency of the proof files. This certificate establishes reference
coercivity, residuals, and scalar envelopes. Exact-profile existence
and stability follow from the separate Newton and transfer arguments.

\subsection{Elementary bounds for the evaluation kernel}

\begin{lemma}\label{sb:sc:kernel}
The moment duals satisfy $C_b,C_c\leq2$. For every even $h\in X$,
\[
 |P_h(x)|\leq400P_0(x)\|h\|_X,\quad
 |P_h(x)-b_h|\leq304\frac{x^2}{1+x^2}\|h\|_X,
\]
and
\[
 |P_h(x)-b_h+c_hx^2/3|\leq13x^{5/2}\|h\|_X\qquad(x>0).
\]
\end{lemma}
\begin{proof}
The dual moment norms are
$C_b=\pi^{-1}(\int_{\mathbb R}x^4f_0\,\idiff x)^{1/2}$ and
$C_c=\pi^{-1}(\int_{\mathbb R}f_0\,\idiff x)^{1/2}$.
For $k(r)=\pi^{-1}r\log|(1+r)/(1-r)|$, the squared potential
dual norm is $k_P(x)^2=\tfrac12\int_0^\infty k(y/x)^2y^4f_0(y)\,\mathrm{d}y$.
The factor $1/2$ reflects the full-line norm of even inputs.
One has $k(r)\leq8r^2/(3\pi)$ for $r\leq1/2$, and
$k(r)\leq8/(3\pi)$ for $r\geq2$. On the intervening interval,
\[
 \int_{1/2}^{2}\log^2\frac{1+r}{|1-r|}\,\mathrm{d}r
 \leq2\{(\log3)^2+2\log3+2\}<\frac{541}{50}.
\]
The reference majorant gives the exact beta-moment bounds
\[
 \begin{gathered}
 M_4:=\int_0^\infty y^4f_0(y)\,\idiff y\leq3\pi/2,\quad
 M_8:=\int_0^\infty y^8f_0(y)\,\idiff y\leq280\pi,\\
 \sup_{y\geq0} y^9f_0(y)\leq4\cdot6^9/10^5<404.
 \end{gathered}
\]
Splitting at $y=x/2,2x$, and using $f_0\leq1$ in the middle
when $x\leq1$, gives
\[
 k_P(x)^2\leq\frac{x^5}{1296\pi^2}
      +\frac{32x^5}{\pi^2}\frac{541}{50}
      +\frac{32M_4}{9\pi^2}<41.
\]
For $x\geq1$ use $y^4f_0\leq404/y^5$ in the middle and obtain
\[
 x^4k_P(x)^2\leq\frac{34M_8}{9\pi^2}
                  +\frac{4\cdot404}{\pi^2}\frac{541}{50}<2300.
\]
Division by $P_0\geq 0.3/(1+x^2)$ bounds $k_P/P_0$ by $320$;
we retain $400$ as a convenient constant. These estimates use only
$\pi>3$ and $\log3<11/10$.

For the canceled kernel put
$t_2(r)=k(r)-2/\pi-2/(3\pi r^2)$. Scaling gives squared dual
norm at most $\tfrac12x^5\int_0^\infty r^4t_2(r)^2\,\mathrm{d}r$.
On $(0,1/2)$ bound $|r^2t_2|$ by $(2r^2+2/3)/\pi$;
on $(2,\infty)$ bound $|t_2|$ by $8/(15\pi r^4)$.
In the middle use
$|r^2t_2|\leq\bigl(8\log(3/|1-r|)+26/3\bigr)/\pi$ and the preceding
log-square integral. The resulting constant is less than $13$.
Thus the dual norm $D_0$ of $P_h-b_h$ satisfies
\[
 D_0\leq\tfrac23x^2+13x^{5/2},\qquad D_0\leq400P_0+2.
\]
For $x\leq1$ the first inequality gives
$D_0(1+x^2)/x^2<28$; for $x\geq1$ the second gives a bound
$300/x^2+2(1+x^2)/x^2\leq304$.
This proves all assertions.
\end{proof}

\subsection{A scalar Newton construction}
\label{sb:sc:newton}

Let $f_0$ be the positive even rational reference of
Section~\ref{sb:sc:certificate}. All Hilbert norms below are full-line norms.
Set
\[
 X=\left\{h\text{ even}:\ \|h\|_X^2
 =\int_{\mathbb R}\frac{|h(x)|^2}{x^4f_0(x)}\,\mathrm{d}x<\infty\right\},
 \qquad r_B=10^{-5}.
\]
The affine ball $f_0+B_X(r_B)$ may contain signed rough functions.
For its elements define
\begin{align*}
 b_f&=\frac2\pi\int_0^\infty f(y)\,\mathrm{d}y,&
 c_f&=\frac2\pi\int_0^\infty\frac{1-f(y)}{y^2}\,\mathrm{d}y,\\
 P_f(x)&=\int_0^\infty\frac{y}{\pi x}
       \log\left|\frac{x+y}{x-y}\right|f(y)\,\mathrm{d}y,&
 a_f&=\frac{3c_f}{6b_f+c_f},\quad g_f=\frac{P_f}{b_f}.
\end{align*}
Write $A_f=a_f^{-1}$ and $C_f=A_f-1$ only in the reconstruction
formulas, and put
\begin{equation}\label{sb:sc:map}
 J_f(x)=\int_0^x\frac{g_f(s)-1}{s g_f(s)}\,\mathrm{d}s,
 \qquad \Phi(f)=\exp\{A_f\log g_f+C_fJ_f\}.
\end{equation}
In particular $A_f=1/3+2b_f/c_f$.

\begin{lemma}\label{sb:sc:mapbounds}
The map $\Phi$ is twice continuously differentiable and compact from
an open neighborhood of the closed enclosing ball into $f_0+X$.
Its derivative is compact. On that ball,
\begin{equation}\label{sb:sc:hessian}
 \|D^2\Phi(f)\|_{X\times X\to X}\leq10^{10}.
\end{equation}
Moreover, with $m_3=(1+x^2/25)^{3/2}$,
\begin{equation}\label{sb:sc:image}
 \left\|\frac{\Phi(f)-\Phi(f_0)}{\sqrt{f_0}}\right\|_\infty
 +\|m_3a_0xP_0\partial_x\bigl(\Phi(f)-\Phi(f_0)\bigr)\|_X
 \leq4\times 10^7\|f-f_0\|_X.
\end{equation}
\end{lemma}

\begin{proof}
Write $C_b,C_c$ for the dual norms of $b$ and $c$ and $D_0(x)$ for
the dual norm of $P_h(x)-b_h$. Lemma~\ref{sb:sc:kernel} and the
reference inequalities give $C_b,C_c\leq2$, $|P_h|\leq400P_0\|h\|_X$,
and $D_0\leq304x^2/(1+x^2)$. In particular
\[
 P_f\geq(1-400r_B)P_0>0,\qquad
 b_f-P_f\geq(1/5-304r_B)\frac{x^2}{1+x^2}>0.
\]
The same reference bounds imply, throughout this signed ball,
\[
 0.6\leq b_f\leq 0.7,\quad 1\leq c_f\leq1.2,\quad 0<a_f<0.7,
 \qquad \frac{0.4}{1+x^2}\leq g_f\leq
       \min\left\{1,\frac{1.25}{1+x^2}\right\}.
\]
For the last upper bound one may retain $b_f\geq 0.61$ and
$P_f\leq1.004\cdot 0.75/(1+x^2)$. Also $C_f\geq 0.43$ and
$P_f\leq 0.47/x^2$ for $x\geq6$.

These inequalities yield $\Phi(f)\leq2e^{-0.27x^2}$. Indeed $\Phi\leq1$
on $[0,1]$. On $[1,6]$, integration in \eqref{sb:sc:map} gives
\[
 \log\Phi\leq1.43\log\frac{1.25}{1+x^2}
                 +0.43\{0.2\log x-0.4(x^2-1)\}.
\]
After adding $0.27x^2$ this is convex in $x^2$ and negative at
both endpoints. For $x\geq6$, use
$J_f(x)\leq J_f(6)+\log(x/6)-(0.6/0.94)(x^2-36)$.
The resulting upper bound for $\log\Phi+0.27x^2$ decreases.
Since the reference satisfies $f_0\geq e^{-x^2/2}/10$, we obtain
\begin{equation}\label{sb:sc:relativeimage}
 \Phi(f)/\sqrt{f_0}\leq G(x):=\sqrt{40}\,e^{-x^2/50}.
\end{equation}

Here are explicit derivative estimates, including the origin
cancellation. Put $d=c_f/(3b_f)$, $v=g_f-1+d x^2$, $u=1-g_f$, and
\[
 E=\log g_f+d x^2=v+\log(1-u)+u,\qquad
 K=J_f+d x^2/2=\int_0^x\frac{v-u^2/g_f}{s}\,\mathrm{d}s.
\]
Exactly $\log\Phi=-x^2+A_fE+C_fK$. The canceled kernel estimate
in Lemma~\ref{sb:sc:kernel} is $\|P_h-b_h+c_hx^2/3\|_{X^*}
\leq13x^{5/2}$. Together with the certified reference remainder
$|P_0-b_0+c_0x^2/3|\leq x^4$, it gives the following outward
bounds on $0<x\leq1/4$; derivatives are in unit $X$ directions:
\[
\begin{array}{c|rrr}
 &|\cdot|/x^{5/2}&|D\cdot|/x^{5/2}&|D^2\cdot|/x^{5/2}\\ \hline
 v&0.22&23&154\\ E&0.27&25&202\\ K&0.11&11&88
\end{array}
\]
For verification use $|Dd|\leq4$, $|D^2d|\leq24$,
$|DA_f|\leq7$, $|D^2A_f|\leq28$, $|u|\leq 0.81x^2$,
$|Du|\leq16x^2$, $|D^2u|\leq102x^2$, and $g_f\geq 0.94$.
The derivatives of $\log(1-u)+u$ are $-u/g_f,-g_f^{-2}$;
those of $u^2/g_f$ are $u(2-u)/g_f^2,2/g_f^3$.
Thus
\[
 |D\log\Phi|\leq64x^{5/2},\quad
 |D^2\log\Phi|\leq1007x^{5/2},\quad
 |D^2\Phi|\leq1135x^{5/2}.
\]
The corresponding $J,DJ,D^2J$ prefixes at $1/4$ are bounded
by $0.03,0.5,4$. Off this interval direct integration gives
\[
\begin{array}{c|r}
 |D\log g_f|&421\\ |D^2\log g_f|&174000\\
 |DJ_f|/(1+x^2)&2106\\ |D^2J_f|/(1+x^2)&1755574\\
 |D\log\Phi|/(1+x^2)&3000\\
 |D^2\log\Phi|/(1+x^2)&2200000
\end{array}
\]
For instance $|P_h/P_f|\leq417$, $|b_h/b_f|\leq4$;
the $D^2J$ integrand is bounded by
$(2\cdot417^2+2\cdot417\cdot4)/g_f$.
Consequently $|D^2\Phi|\leq11200000\Phi(1+x^2)^2$ there.
Using \eqref{sb:sc:relativeimage}, its $X$ norm factor is below $440$:
expand $(1+x^2)^4/x^4=x^{-4}+4x^{-2}+6+4x^2+x^4$ and integrate
the last three terms against $e^{-x^2/25}$. The origin contributes
less than $1000$, proving \eqref{sb:sc:hessian}.

For completeness the high derivative estimate involves no pointwise
Hilbert bound. The exact identity
\[
 P_h'/x=-H(h/x)+\int_0^1t^2H(h/x)(tx)\,\mathrm{d}t
 \quad\hbox{gives}\quad \|P_h'/x\|_2\leq\tfrac75\|h/x\|_2
 \leq\tfrac75\|h\|_X.
\]
Put $\|z\|_{X_3}=\|m_3z\|_X$, $R=(1-400r_B)^{-1}\leq25/24$,
$R_p=400R\leq1250/3$, $L=D\log\Phi[h]$, and $Q=|g_f-1|/g_f$.
Use $Q\leq27x^2$, $|D\log g_f[h]|\leq7157x^2/(1+x^2)$ and
\[
 B(V)=\|xP_0'V\|_{X_3}+\tfrac75r_B\|Vm_3/\sqrt{f_0}\|_\infty.
\]
The rational derivative bound $|xP_0'|\leq2x^2/(1+x^2)$ and
Gaussian moments give $\|Gm_3\|_2<500$, $\|Gm_3\|_\infty<700$,
$B(\Phi)<1001$, $B(\Phi|L|)<3070000$,
$\|P_0\Phi|L|Q\|_{X_3}<40500000$,
$\|P_0\Phi Q\|_{X_3}<13500$ and
$\|\Phi D\log g_f[h]\|_{X_3}<3578500$.
Differentiating $D\Phi[h]=\Phi L$ gives seven terms bounded by
\[
\begin{split}
 &\tfrac34\,2R B(\Phi|L|),\quad
 \tfrac34\|P_0\Phi|L|Q\|_{X_3},\quad
 \tfrac34\,7R B(\Phi),\quad \tfrac34\,7\|P_0\Phi Q\|_{X_3},\\
 &\tfrac34\,2R\tfrac75\|Gm_3\|_\infty,\quad
 \tfrac34\,2RR_p B(\Phi),\quad
 \tfrac34\,0.7R\|\Phi D\log g_f[h]\|_{X_3}.
\end{split}
\]
Their sum is below $3.7858441\times10^7$. The point norm is below $10^6$,
by \eqref{sb:sc:relativeimage} and the bound for $L$.
Integrating these derivative bounds proves \eqref{sb:sc:image}.

The same formulas justify differentiating twice on the rough ball.
Locally $P_h'=-Hh+\int_0^1tHh(tx)\,\mathrm{d}t$ supplies $H^1$ bounds.
The uniform $O(x^{5/2})$ origin gain, local Rellich compactness and
Gaussian tails prove compactness into $X$ and compactness of $D\Phi$.
They also provide dominators for continuity of the first two derivatives.
\end{proof}

\begin{proposition}\label{sb:sc:constructed}
There is a positive smooth Schwartz fixed point $f_*$ of
\eqref{sb:sc:map} satisfying $f_*(0)=1$, $f_*''(0)=-2$, and
$\|f_*-f_0\|_X<10^{-15}$. Its potential is positive and
$a_{\mathrm c}=3c_{f_*}/(6b_{f_*}+c_{f_*})$ lies in $(0.689,0.690)$.
\end{proposition}

\begin{proof}
Define the frozen transport and residual by
\[
 T_0z=-a_0xP_0z'+\mathcal C_0z,
\qquad \mathcal C_0=(1-a_0)(P_0-b_0)+xP_0',\qquad
 F_0=-a_0xP_0f_0'+\mathcal C_0f_0.
\]
Its domain is $\mathcal D=\{z\in X:a_0xP_0z'\in X\}$,
with distributional derivative. Put $\mathrm{d}s/\mathrm{d}x=(a_0xP_0)^{-1}$.
The unitary change $z\mapsto\sqrt{2a_0xP_0/(x^4f_0)}z$ turns
$T_0$ into $-\partial_s-D_T$ on $H^1(\mathbb R)$, where
\[
 D_T=\frac{(1-a_0)b_0}{2}+\frac{4a_0-1}{2}P_0
 -\frac{1+a_0}{2}xP_0'-\frac{F_0}{2f_0}\geq 0.09.
\]
Thus $T_0$ is closed, densely defined and invertible with inverse
norm at most $1/0.09$. The displayed residual term must be retained
because $f_0$ is not stationary.

For a smooth multiplier $\phi$ let $K_\phi h=DT_{f_0}[h]\phi$,
where $T_f=-a_fxP_f\partial_x+(1-a_f)(P_f-b_f)+xP_f'$.
Writing $a_h=Da_{f_0}[h]$, explicitly
\[
 K_\phi h=-(a_hxP_0+a_0xP_h)\phi'
 +\bigl(-a_h(P_0-b_0)+(1-a_0)(P_h-b_h)+xP_h'\bigr)\phi.
\]
If $C_a=\|Da_{f_0}\|\leq4$, its bounded norm is at most
\begin{equation}\label{sb:sc:boundedK}
 (C_a+400a_0)\|xP_0\phi'\|_X
 +C_a\|(P_0-b_0)\phi\|_X+(1-a_0)\|D_0\phi\|_X
 +\tfrac75\|\phi/\sqrt{f_0}\|_\infty.
\end{equation}
This applies to $f_0$, $\Phi_0=\Phi(f_0)$, and $q=\Phi_0-f_0$.
The certificate gives, on the smooth annular core and hence on
$\mathcal D$, $\operatorname{Re}\langle M_0z,z\rangle_X
\leq-\eta_0\|z\|_X^2$, $\eta_0>0.0109535$, where $M_0=T_0+K_{f_0}$, and
$\|K_{f_0}\|<1.542$. It also bounds
$\|q\|_X<1.913\times 10^{-17}$,
$\|xP_0q'\|_X<1.553\times 10^{-17}$, and
$\|q/\sqrt{f_0}\|_\infty<7.711\times 10^{-15}$.
These follow from $T_0q=-F_0$, the native and high transport energies,
and their flux estimate. Formula~\eqref{sb:sc:boundedK}, using
$D_0\leq402$, yields $\delta_K:=\|K_q\|<3\times 10^{-14}$.

The image estimate \eqref{sb:sc:image} implies that
\(D\Phi(f_0):X\to\mathcal D(T_0)\) is bounded in the graph norm.
Differentiating the identity \(T_f\Phi(f)=0\) therefore gives the
operator identity
\begin{equation}\label{eq:critical-derivative-factorization}
 T_0D\Phi(f_0)=-K_{\Phi(f_0)}.
\end{equation}
Let \(\widetilde M=T_0+K_{\Phi(f_0)}\). The reference energy estimate
and the bound for \(K_{\Phi(f_0)-f_0}\) give
\begin{equation}\label{eq:critical-corrected-coercivity}
 \operatorname{Re}\langle\widetilde M z,z\rangle_X
 \leq-\kappa\|z\|_X^2,
 \qquad z\in\mathcal D(T_0),
 \qquad \kappa>0.01095.
\end{equation}
The compactness of \(D\Phi(f_0)\), proved in
Lemma~\ref{sb:sc:mapbounds}, now supplies the required Fredholm
argument directly. Indeed, if \(z=D\Phi(f_0)z\), then
\(z\in\mathcal D(T_0)\), and
\eqref{eq:critical-derivative-factorization} gives
\(\widetilde Mz=0\). Equation~\eqref{eq:critical-corrected-coercivity}
implies \(z=0\). The Fredholm alternative shows that
\(I-D\Phi(f_0)\) is invertible on \(X\).

The same identities show that \(\widetilde M\) is onto \(X\).
For \(y\in X\), set
\(z=(I-D\Phi(f_0))^{-1}T_0^{-1}y\). Since
\(z=T_0^{-1}y+D\Phi(f_0)z\), both summands belong to
\(\mathcal D(T_0)\), and \(\widetilde Mz=y\).
Consequently,
\begin{equation}\label{eq:critical-inverse-bound}
 \|\widetilde M^{-1}\|_{X\to X}\leq\kappa^{-1},
 \qquad
 (I-D\Phi(f_0))^{-1}
   =I-\widetilde M^{-1}K_{\Phi(f_0)}.
\end{equation}
In particular, the bound
\(\beta:=\|(I-D\Phi(f_0))^{-1}\|<142\) follows from the certified
bound for \(K_{\Phi(f_0)}\) and \eqref{eq:critical-inverse-bound}.

The preconditioned defect is exactly
$ (I-D\Phi_0)^{-1}q=-\widetilde M^{-1}F_0$ and has norm
$\eta<1.58\times 10^{-16}$, since the certified residual obeys
$\|F_0\|_X<1.721\times 10^{-18}$.

The Newton map $f\mapsto f+(I-D\Phi_0)^{-1}(\Phi(f)-f)$
therefore maps the closed radius-$r$ ball into itself and contracts,
where $r=10^{-15}$: indeed
\[
 \beta10^{10}r<0.00142,
 \qquad \eta+\tfrac12\beta10^{10}r^2
 <1.588\times 10^{-16}<r.
\]
Its fixed point is a fixed point of $\Phi$ and is positive.
To enclose its parameter, the moment bounds on the larger radius-$r_B$
ball give
\[
 b_0-2r_B\leq b_f\leq b_0+2r_B,\qquad
 c_0-2r_B\leq c_f\leq c_0+2r_B.
\]
The exact rational check verifies
\[
 \frac{689}{1000}
 <\frac{3(c_0-2r_B)}{6(b_0+2r_B)+c_0-2r_B}
 \leq a_f\leq
 \frac{3(c_0+2r_B)}{6(b_0-2r_B)+c_0+2r_B}
 <\frac{690}{1000}.
\]
This enclosure applies in particular to \(a_{\mathrm c}=a_{f_*}\).

Finally the canceled formula gives $f_*=1-x^2+O(x^{5/2})$.
To justify ordinary smoothness, put $g=P_{f_*}/b_{f_*}$ and
$F_N=f_*/g^N$. All $F_N$ have Gaussian envelopes, and
\[
 F_N'=(A-N)F_{N+1}g'+CF_{N+1}(g-1)/x.
\]
Initially $g'\in L^2$ by the Hilbert identity, so every $F_N\in H^1$.
If all are in $H^m$, the same identity gives $g'\in H^m$, while
$(g-1)/x=\int_0^1g'(tx)\,\mathrm{d}t\in H^m$ by the dilation bound
$\int_0^1t^{j-1/2}\mathrm{d}t=2/(2j+1)$. The Sobolev algebra property
proves all $F_N\in H^{m+1}$. Since $g\asymp(1+x^2)^{-1}$,
differentiating $f_*=g^NF_N$ and choosing $N$ arbitrarily large
proves Schwartz decay of every derivative. The normalized jets now
follow classically.
\end{proof}

\subsection{Coercivity of the exact profile}
\label{sb:sc:transfer}

The reference weight is useful here: it has polynomial growth at infinity,
so no uniform comparison of the rational and Gaussian logarithmic slopes
is required. We give the quantitative form estimate used to conclude the
proof of Theorem~\ref{sb:sc:profile}.

Put $v=f_*-f_0$, $K=(1+x^2/25)^3$, and
\[
 e=\|v\|_X,\quad A_v=\|v/\sqrt{f_0}\|_\infty,\quad
 H_v^{(3)}=\|\sqrt K\,a_0xP_0v'\|_X.
\]
The frozen high transport estimate and its flux identity give
\[
 H_q^{(3)}\leq(1+\|\mathcal C_0\|_\infty/0.09)
                   \|F_0\|_{X_3}<4.901\times 10^{-16},\qquad
 \|q/\sqrt{f_0}\|_\infty<7.711\times 10^{-15}.
\]
To check the flux constant, the conjugated transport equation
$T_0q=-F_0$ yields, after integration from zero,
\[
 \|q/\sqrt{f_0}\|_\infty
 \leq\frac{\|F_0\|_{X_3}}{2\sqrt{d_3}}
          \sup_{x>0}\frac{x^2}{\sqrt{a_0xP_0K}},\qquad d_3=0.09.
\]
The full-line convention accounts for the factor $2$ in the denominator.
The supremum is one of the scalar rational checks.
Applying \eqref{sb:sc:image} to
$v=q+\Phi(f_*)-\Phi(f_0)$ gives
\begin{equation}\label{sb:sc:strong}
 e<10^{-15},\qquad A_v,H_v^{(3)}<4.1\times 10^{-8}.
\end{equation}

Here are all multipliers in the form comparison:
\[
\begin{gathered}
 S_j=\sup |x|^j\sqrt{f_0}\ (j=1,2,3),\quad
 S_d=\sup\frac{x^2\sqrt{f_0}}{a_0P_0\sqrt K},\\
 C_p=\sup\frac{x^2}{a_0P_0(1+x^2)\sqrt K},\quad
 C_r=\sup\frac1{a_0P_0\sqrt K},\\
 N_p^2=\int_{\mathbb R}\frac{(W_0')^2}{(1+x^2)^2f_0}\,\idiff x,\quad
 N_r^2=\int_{\mathbb R}\frac{(f_0')^2}{x^2f_0}\,\idiff x,\quad
 J_0=\sup P_0|j-6|.
\end{gathered}
\]
Their certified upper bounds are, respectively,
\[
\begin{array}{c|rrrrr}
 &S_1&S_2&S_3&S_d&C_p\\ \hline
 &0.646&0.946&2.015&10.509&31.531\\[2pt]
 &C_r&N_p&N_r&J_0&\|U_0'\|_\infty\\ \hline
 &32.096&0.887&2.509&3.682&0.614
\end{array}
\]
All bounds concern whole intervals, including both endpoints.

For clarity, the complete perturbation estimate is given next. Set
$C_a=4$, $C_b=2$, $C_P=400$, $C_D=304$ and
$a_+=0.7$, $\alpha_+=0.4$, $d_+=0.4$. These bound
$|Da_f|$, the moment and kernel norms, $a_{\mathrm c}$,
$(3-a_{\mathrm c})/6$, and $1-a_{\mathrm c}$ on the enclosing ball. Define
\[
 M_v=(S_2+S_3)e+S_dH_v^{(3)},
\]
\begin{align*}
 \epsilon_{\rm loc}={}&(1+a_+/2)M_v
 +C_ae\{b_0+\|U_0'\|_\infty/2+J_0/2\}
 +d_+C_be+(a_+/2)C_PeJ_0,\\
 \epsilon_{\rm nl}={}&S_1A_v
 +C_D\{C_aeN_p+a_+(e+C_pH_v^{(3)})\}\\*
 &+\sqrt{b_0/\pi}\{(C_ae/6)N_r+\alpha_+C_rH_v^{(3)}\}.
\end{align*}
Then
\begin{equation}\label{sb:sc:formdifference}
 |\operatorname{Re}\langle (L_*-L_0)h,h\rangle_{\rho_0}|
 \leq(\epsilon_{\rm loc}+\epsilon_{\rm nl})\|h\|_{\rho_0}^2.
\end{equation}

To verify this inequality, integration by parts in the local transport
difference gives the multiplier
\[
 \delta U'-\delta\bigl((1-a)b\bigr)+\tfrac12\delta(aU)'
                          +\tfrac12\delta(aP)(j-6).
\]
Weighted embeddings give $\|xv\|_{H^1}\leq M_v$, and hence
$\|\delta U'\|_\infty=\|H(-xv)\|_\infty\leq M_v$.
The moment bounds and $|P_*-P_0|\leq C_PeP_0$ yield
$\epsilon_{\rm loc}$.
For the nonlocal part set $g=h/x^2$ and use exactly
\[
 Rh=Wx^2Hg-aW'x^3J_3Hg-\alpha c_1(W-xW'),\quad
 c_1=Hg(0),\quad \alpha=(3-a)/6.
\]
Here $\|g\|_2\leq S_1\|h\|_{\rho_0}$ and
$|c_1|\leq\sqrt{b_0/\pi}\|h\|_{\rho_0}$.
For $w=-h/x$, the identity
$J_3Hg=(P_w-b_w)/x^2$ gives the point bound
$|J_3Hg|\leq C_D\|h\|_{\rho_0}/(1+x^2)$.
Finally
\[
 \left\|\frac{v+xv'}{(1+x^2)\sqrt{f_0}}\right\|_2
 \leq e+C_pH_v^{(3)},\qquad
 \left\|\frac{v'}{x\sqrt{f_0}}\right\|_2\leq C_rH_v^{(3)}.
\]
These prove $\epsilon_{\rm nl}$, including the change in
$\alpha$ and the primitive coefficient $aW'$. No pointwise
$L^2$ Hilbert-transform estimate is used.

Substitution of \eqref{sb:sc:strong} and the displayed scalar bounds
gives $\epsilon_{\rm loc}+\epsilon_{\rm nl}<10^{-3}$.
The exact rational assembly in \texttt{check.py} gives the sharper
bound $0.000270$, a Newton contraction factor below $0.001418$, and
a Newton image radius below $1.579\times 10^{-16}$.
Since the reference gap of $L_0$ exceeds $0.06302$,
\eqref{sb:sc:formdifference} proves
\eqref{sb:sc:profile-coercivity}, with reserve.
Together with Proposition~\ref{sb:sc:constructed}, this completes
the proof of Theorem~\ref{sb:sc:profile}.

The weight $\rho_0$ is comparable to $x^{-6}$ at zero and to $x^4$
at infinity. Thus smooth Schwartz perturbations with vanishing first and third derivatives at the origin
belong to its energy space. This is precisely the profile input needed
by the polynomial-space stability and periodic cutoff argument in Section~\ref{sb:paper:critical}. The parameter remains fixed during that
physical evolution.

\section{Finite data, interval arithmetic, and reproduction}
\label{sb:repro:appendix}

The computer-assisted component is the scalar verification supporting
Theorem~\ref{sb:sc:profile}. This appendix specifies the exact data,
directed arithmetic, and reproduction procedure. The verification code
and data are available at
\begin{center}
\url{https://github.com/HouGroup2026/periodic-gclm-verification}.
\end{center}
These electronic files form part of the proof. All remaining existence,
regularity, spectral, and nonlinear arguments are given in the paper.

\subsection{Exact reference and proof data}

The reference in \eqref{sb:sc:reference} is specified by 97 rational
coefficients, ordered by ascending powers of \(z\), in
\texttt{critical\_scalar/reference.json}.
Each numerator and denominator is an integer written exactly in base ten.
The spatial scale is \(5\), and the tail order is \(5\).
Reproduction uses these exact coefficients and does not require repeating
the numerical search that produced them.

The program first forms the Hilbert transforms, potential, moments,
rate coefficients, residual polynomial, damping polynomials, and Gram
integrands by rational arithmetic. All identities in
\eqref{sb:sc:hilbert}--\eqref{sb:sc:integrals} therefore hold exactly for
the reference used by the program. The report \texttt{scalar.json}
records these polynomials, the exact reference coefficients, and every
required scalar bound. The separate \texttt{check.py} verifies the
explicit scalar hypotheses used in the analytic Newton argument and
the final form comparison, using exact rational arithmetic throughout.
Stored binary floating-point endpoints are interpreted as their exact
binary rational values in this last check.

\subsection{Directed arithmetic and full-interval bounds}

The scalar interval implementation assumes IEEE 754 binary64 arithmetic.
For addition, subtraction, multiplication, and reciprocal, compute each
endpoint and move it by one representable value toward the appropriate
infinity. Rational-to-binary conversion is checked against the exact
rational value and enlarged when needed. Division by an interval
containing zero is rejected. Thus every successful arithmetic operation
contains its exact real result, including at underflow.

The endpoints for \(\pi\) follow from Machin's identity
\[
 \pi=16\arctan(1/5)-4\arctan(1/239),
\]
\[
 \arctan(q^{-1})=\sum_{j=0}^{n-1}
       \frac{(-1)^j}{(2j+1)q^{2j+1}}+R_n,
 \quad |R_n|\leq\frac1{(2n+1)q^{2n+1}},
\]
with the alternating sign of the remainder and \(n=50,15\), respectively.
Sine is enclosed by its alternating Taylor series on the required
bounded angle interval. Logarithms and exponentials use correctly
rounded Decimal operations at 30 digits, enlarged by the neighboring
Decimal values before outward binary conversion. Square roots are
treated analogously at 45 digits. No library trigonometric evaluation
is used as an enclosure of a continuum integral.

For polynomial evaluation write \(u=\sin^2\theta\),
\(x=5\tan\theta\), and partition \([0,1]\) into 96 intervals with
exact dyadic endpoints. Floating-point sine values only propose this
partition; the resulting rational endpoints are checked to be strictly
increasing and to cover \([0,1]\). No accuracy of that proposal is
needed. On a chart \(u=c+rt\), \(|t|\leq1\), both \(c,r\) are
positive dyadic rationals. The transformed polynomial coefficients are
computed twice, with Decimal rounding downward and upward at 180 digits.
The affine Horner recurrence involves nonnegative multipliers \(c,r\),
so the directed recurrences enclose every coefficient even when the
original coefficients have different signs.

If \(p(t)=\sum_{j=0}^n p_jt^j\), retain the coefficients through
degree 32 and bound the omitted part by
\[
 \left|\sum_{j>32}p_jt^j\right|
      \leq\sum_{j>32}|p_j|,\qquad |t|\leq1.
\]
Apply the same procedure to the first and second derivatives. For a
subinterval \(I\subset[-1,1]\) and a fixed midpoint \(t_0\in I\),
use the interval Taylor enclosure
\[
 p(I)\subset p(t_0)+p'(t_0)(I-t_0)
                      +\tfrac12p''(I)(I-t_0)^2.
\]
If a cell intersects more than one chart, take the hull of the chart
enclosures. Every discarded coefficient is included in the remainder.

Divide the full angle interval \([0,\pi/2]\) into 4096 equal cells.
The outward sine-square endpoints enclose their images in \(u\).
For each complete cell, evaluate the polynomial comparisons and all
integrands in \eqref{sb:sc:integrals} and in the multiplier estimates.
An interval enclosure for an integrand, multiplied by the cell-width
interval and summed outward, bounds its integral. The compactified
formulas are regular at both endpoints; hence the first and last cells
include the origin and the entire infinite spatial tail. No tail is
discarded and no sign condition is inferred from sampled values.

The four-dimensional energy comparison uses the Gram matrix in
\eqref{sb:sc:gram}. Rational matrix identities reduce the desired
comparison to the inertia of a symmetric interval matrix; interval
\(LDL^{T}\) elimination encloses each nonzero pivot and verifies its
sign. All pivots are separated from zero. The program also checks
the required denominator and damping inequalities on every cell.
It aborts if any of these checks is unresolved.

\subsection{Files and commands}

The verification files are in the directory \texttt{critical\_scalar}
of the repository:
\begin{center}
\begin{tabular}{ll}
\texttt{reference.json}&exact rational reference\\
\texttt{verify.py}&polynomial and continuum interval verification\\
\texttt{interval.py}&outward scalar arithmetic\\
\texttt{scalar.json}&complete recorded scalar report\\
\texttt{check.py}&exact rational assembly of the analytic inequalities\\
\texttt{nonlinear.json}&readable summary of that assembly\\
\texttt{manifest.json}&integrity metadata for the proof files\\
\texttt{README.md}&certificate documentation
\end{tabular}
\end{center}
The repository README gives further instructions. With Python 3.10 or
later and its standard library, check the recorded certificate from the
repository root:
\begin{verbatim}
python3 reproduce.py --check-only
\end{verbatim}
This checks the integrity of the proof files and verifies the scalar
hypotheses for the analytic arguments using the recorded bounds.
For a full recomputation, run:
\begin{verbatim}
python3 reproduce.py \
    --output-dir /tmp/gclm-scalar-reproduced
\end{verbatim}
The output directory must be absent or empty. The command copies the
frozen files, recomputes all continuum bounds from the rational reference,
and runs the exact checker on the newly produced report. The supplied
frozen files remain unchanged. The log file \texttt{continuum.log}
records the computation.

The critical existence assertion uses the finite electronic data and
interval verification specified in this appendix. The reproduction
program runs with the Python standard library and recomputes the complete
continuum bounds from the exact rational reference.

\bibliographystyle{myalpha}
\bibliography{references}
\end{document}